\documentclass[preprint,12pt]{elsarticle}

\usepackage{amssymb}
\usepackage{amsmath}
\usepackage{tikz}
\usepackage{verbatim}
\usepackage{mathtools}
\usetikzlibrary{shapes}
\usetikzlibrary{plotmarks}

\usepackage{subcaption,siunitx,booktabs}
\usepackage{pgfplotstable}

\journal{Computer Methods in Mechanics and Engineering}

\begin{document}

\begin{frontmatter}


 \title{A $p$-multigrid method enhanced with an ILUT smoother and its comparison to $h$-multigrid methods within Isogeometric Analysis.}
 \author{R. Tielen \corref{cor1} \fnref{label1}}
 \author{M. M\"{o}ller  \fnref{label1}}
 \author{D. G\"{o}ddeke \fnref{label2}}
 \author{C. Vuik \fnref{label1}}
 \cortext[cor1]{r.p.w.m.tielen@tudelft.nl}

\title{}

 \address[label1]{Delft Institute of Applied Mathematics, Delft University of Technology, \\ Van Mourik Broekmanweg 6, 2628XE, Delft}
 \address[label2]{Institute for Applied Analysis and Numerical Simulation, University of Stuttgart  \\ Allmandring 5b, 70569, Stuttgart}

\begin{abstract}
Over the years, Isogeometric Analysis has shown to be a successful alternative to the Finite Element Method (FEM). However, solving the resulting linear systems of equations efficiently remains a challenging task. In this paper, we consider a $p$-multigrid method, in which coarsening is applied in the approximation order $p$ instead of the mesh width $h$. Since the use of classical smoothers (e.g. Gauss-Seidel) results in a $p$-multigrid method with deteriorating performance for higher values of $p$, the use of an ILUT smoother is investigated. Numerical results and a spectral analysis indicate that the resulting $p$-multigrid method exhibits convergence rates independent of $h$ and $p$. In particular, we compare both coarsening strategies (e.g. coarsening in $h$ or $p$) adopting both smoothers for a variety of two and threedimensional benchmarks.   
\end{abstract}


\end{frontmatter}


\section{Introduction}
\label{intro}

Isogeometric Analysis (IgA) \cite{hugh} has become widely accepted over the years as an alternative to the Finite Element Method (FEM). The use of B-spline basis functions or Non-Uniform Rational B-splines (NURBS) allows for a highly accurate representation of complex geometries and establishes the link between computer-aided design (CAD) and computer-aided engineering (CAE) tools.  Furthermore, the $C^{p-1}$ continuity of the basis functions offers a higher accuracy per degree of freedom compared to standard FEM \cite{sang}.  

\noindent IgA has been applied with success in a wide range of engineering fields, such as structural mechanics \cite{cot}, fluid dynamics \cite{baz} and shape optimization \cite{wall}. Solving the resulting linear systems efficiently is, however, still a challenging task. The condition numbers of the mass and stiffness matrices increase exponentially with the approximation order $p$, making the use of (standard) iterative solvers inefficient. On the other hand, the use of (sparse) direct solvers is not straightforward due to the increasing stencil of the basis functions and increasing bandwidth of matrices for higher values of $p$. Furthermore, direct solvers may not be practical for large problem sizes due to memory constraints, which is a common problem in high-order methods in general. 

\noindent Recently, various solution techniques have been developed for discretizations arising in Isogeometric Analysis. For example, preconditioners have been developed based on fast solvers for the Sylvester equation \cite{tani} and overlapping Schwarz methods \cite{beirao}.     

\noindent As an alternative, geometric multigrid methods have been investigated, as they are considered among the most efficient solvers in Finite Element Methods for elliptic problems. However, the use of standard smoothers like (damped) Jacobi or Gauss-Seidel leads to convergence rates which deteriorate for increa-sing values of $p$ \cite{gah}. It has been noted in \cite{dona} that very small eigenvalues associated with high-frequency eigenvectors cause this behaviour. This has lead to the development of non-classical smoothers, such as smoothers based on mass smoothing \cite{hof,tak1,tak2} or overlapping multiplicative Schwarz methods \cite{gaspar}, showing convergence rates independent of both $h$ and $p$. 

\noindent $p$-Multigrid methods can be adopted as an alternative solution strategy.  In contrast to $h$-multigrid methods, a hierarchy is constructed where each level represents a different approximation order. Throughout this paper, the coarse grid correction is obtained at level $p=1$. Here, B-spline functions coincide with linear Lagrange basis functions, thereby enabling the use of well known solution techniques for standard Lagrangian Finite Elements. Furthermore, the stencil of the basis functions and bandwidth of the matrix is significantly smaller at level $p=1$, reducing both assembly and factorization costs of the multigrid method.     

\noindent $p$-Multigrid methods have mostly been used for solving linear systems arising within the Discontinuous Galerkin method \cite{fid, luo, luo2, paulien}, where a hierarchy was constructed until level $p=0$. However, some research has been performed for continuous Galerkin methods \cite{hel} as well, where the coarse grid correction was obtained at level $p=1$. 

\noindent Recently, the authors applied a $p$-multigrid method, using a Gauss-Seidel smoother, in the context of IgA \cite{roel}. As with $h$-multigrid methods, a depen-dence of the convergence rate on $p$ was reported. In this paper, a $p$-multigrid method is presented that makes use of an Incomplete LU Factorization based on a dual threshold strategy (ILUT) \cite{saad} for smoothing.  The spectral properties of the resulting $p$-multigrid method are analyzed adopting both smoothers. Numerical results are presented for Poisson's equation on a quarter annulus, an L-shaped (multipatch) geometry and the unit cube. Furthermore, the convection-diffusion-reaction (CDR) equation is considered on the unit square. The use of ILUT as a smoother improves the performance of the $p$-multigrid method significantly and leads to convergence rates which are independent of $h$ and $p$. 

\noindent Compared to standard $h$-multigrid methods, both the coarsening strategy and smoother are adjusted in the proposed $p$-multigrid method. Therefore, a comparison study is performed in terms of convergence rates and CPU times between $p$-multigrid and $h$-multigrid methods using both smoothers. Furthermore, the $p$-multigrid method with ILUT as a smoother is compared to an $h$-multigrid method adopting a smoother based on stable splittings of spline spaces \cite{tak1}. Finally, to show the versatility of the proposed $p$-multigrid method, it is applied to solve linear systems of equations resulting from THB-spline discretizations \cite{juttler}. 

\noindent This paper is organised as follows. In Section \ref{mp} the model problem, the basics of IgA and the spatial discretization are considered. Section \ref{pm} presents the $p$-multigrid method in detail, together with the proposed ILUT smoother. A spectral analysis is performed with both smoothers and coarsening strategies and is discussed in Section \ref{num}. In Section \ref{numex}, numerical results for the considered benchmarks are presented. Finally, conclusions are drawn in Section \ref{con}.

\section{Model problem and IgA discretization}
\label{mp}

To assess the quality of the $p$-multigrid method, the convection-diffusion-reaction (CDR) equation is considered as a model problem:
\vskip-.6cm
\begin{eqnarray} \label{pois}
-\nabla \cdot (\mathbf{D} \nabla u) + \mathbf{v} \cdot \nabla u + R u = f, \hspace{0.2cm} \text{on } \Omega,
\end{eqnarray}

\noindent where $\mathbf{D}$ denotes the diffusion tensor, $\mathbf{v}$ a divergence-free velocity field and $R$ a source term. Here, $\Omega \subset \mathbb{R}^2$ is a connected, Lipschitz domain, $f \in L^2(\Omega)$ and $u = 0$ on the boundary $\partial \Omega$. Let $\mathcal{V} = H^1_0(\Omega)$ denote the space of functions in the Sobolev space $H^1(\Omega)$ that vanish on $\partial \Omega$.  The variational form of  \eqref{pois} is then obtained by multiplication with an arbitrary test function $v \in \mathcal{V}$ and application of integration by parts: \\
\\
Find $u \in \mathcal{V}$ such that
\vskip-.6cm
\begin{eqnarray} \label{weak}
a(u, v) = (f,v) \hspace{0.5cm} \forall v \in \mathcal{V},
\end{eqnarray}
\noindent where
\vskip-.6cm
\begin{eqnarray}
a(u,v) = \int_{\Omega} (\mathbf{D}\nabla u) \cdot \nabla v + (\mathbf{v} \cdot \nabla u) v + R u v\ \text{d}\Omega 
\end{eqnarray}
and
\vskip-.6cm
\begin{eqnarray}
(f,v) = \int_{\Omega} f v \ \text{d}\Omega. 
\end{eqnarray}

\noindent The physical domain $\Omega$ is then parameterized by a geometry function 
\vskip-.6cm
\begin{eqnarray}
\mathbf{F}: \Omega_0 \rightarrow \Omega, \hspace{1.2cm} \mathbf{F}(\boldsymbol{\xi})=\mathbf{x},
\end{eqnarray}

\noindent The geometry function $\mathbf{F}$ describes an invertible mapping connecting the parameter domain $\Omega_0 = (0,1)^2$ with the physical domain $\Omega$. In case $\Omega$ cannot be described by a single geometry function, the physical domain is divided into a collection of $K$ non-overlapping subdomains $\Omega^{(k)}$ such that

\vskip-.6cm
\begin{eqnarray}
\overline{\Omega} = \bigcup_{k=1}^K \ \overline{\Omega}^{(k)}.
\end{eqnarray}

\noindent A family of geometry functions $\mathbf{F}^{(k)}$ is then defined to parameterize each subdomain $\Omega^{(k)}$ separately: 

\vskip-.6cm
\begin{eqnarray}
\mathbf{F}^{(k)}: \Omega_0 \rightarrow \Omega^{(k)}, \hspace{1.2cm} \mathbf{F}^{(k)}(\boldsymbol{\xi})=\mathbf{x}.
\end{eqnarray} 

\noindent In this case, we refer to $\Omega$ as a multipatch geometry consisting of $K$ patches.

\subsection*{B-spline basis functions}
Throughout this paper, the tensor product of univariate B-spline basis functions of order $p$ is used for spatial discretization, unless stated otherwise. Univariate B-spline basis functions are defined on the parameter domain $\hat{\Omega} = (0,1)$ and are uniquely determined by their underlying knot vector 
\vskip-.6cm
\begin{eqnarray}
\Xi = \{ \xi_1, \xi_2 , \ldots , \xi_{N+p}, \xi_{N+p+1} \},
\end{eqnarray}

\noindent consisting of a sequence of non-decreasing knots $\xi_i \in \hat{\Omega}$. Here, $N$ denotes the number of univariate basis functions of order $p$ defined by this knot vector. 

\noindent B-spline basis functions are defined recursively by the Cox-de Boor formula \cite{cox}. The resulting B-spline basis functions $\phi_{i,p}$ are non-zero on the interval $[\xi_{i}, \xi_{i+p+1})$, implying a compact support that increases with $p$. Furthermore, at every knot $\xi_i$ the basis functions are $C^{p-m_i}$-continuous, where $m_i$ denotes the mutiplicity of knot $\xi_i$. Finally, the basis functions possess the partition of unity property:
\vskip-.6cm
\begin{eqnarray}
\sum_{i=1}^{N} \phi_{i,p}(\xi) = 1 \quad \forall \xi \in [\xi_1,\xi_{n+p+1}].
\end{eqnarray}

\noindent Throughout this paper, B-spline basis functions are considered based on an open uniform knot vector with knot span size $h$, implying that the first and last knots are repeated $p+1$ times. As a consequence, the basis functions considered are $C^{p-1}$ continuous and interpolatory only at the two end points. 

\noindent For the two-dimensional case, the tensor product of univariate B-spline basis functions $\phi_{i_x,p}(\xi)$ and $\phi_{i_y,q}(\eta)$ of order $p$ and $q$, respectively, with maximum continuity is adopted for the spatial discretization:
\vskip-.6cm
\begin{eqnarray}
\Phi_{\vec{i},\vec{p}}(\boldsymbol{\xi}) := \phi_{i_x,p}(\xi) \phi_{i_y,q}(\eta), \hspace{1cm} \vec{i} = (i_x,i_y),\ \vec{p} = (p,q).
\end{eqnarray}

\noindent Here, $\vec{i}$ and $\vec{p}$ are multi indices, with $i_x = 1, \ldots , n_x$ and $i_y = 1, \ldots , n_y$ denoting the univariate basis functions in the $x$ and $y$-dimension, respectively. Furthermore, $i = i_x n_x + (i_y-1) n_y$ assigns a unique index to each pair of univariate basis functions, where $i = 1, \ldots N_{\rm dof}$. Here, $N_{\rm dof}$ denotes the number of degrees of freedom, or equivalently, the number of tensor-product basis functions and depends on both $h$ and $p$. In this paper, all univariate B-spline basis functions are assumed to be of the same order (i.e. $p=q$). The spline space $\mathcal{V}_{h,p}$ can then be written, using the inverse of the geometry mapping $\mathbf{F}^{-1}$ as pull-back operator, as follows: 
\vskip-.6cm
\begin{eqnarray}
\mathcal{V}_{h,p} := \text{span}\left  \{ \Phi_{\vec{i},\vec{p}} \circ \mathbf{F}^{-1} \right \}_{i=1, \ldots,N_{\rm dof}}.
\end{eqnarray}

\noindent The Galerkin formulation of  \eqref{weak} becomes: Find $u_{h,p} \in \mathcal{V}_{h,p}$ such that
\vskip-.6cm
\begin{eqnarray}
a(u_{h,p},v_{h,p}) = (f,v_{h,p}) \hspace{0.7cm} \forall v_{h,p} \in \mathcal{V}_{h,p}.
\end{eqnarray}

\noindent The discretized problem can be written as a linear system
\vskip-.6cm
\begin{eqnarray} \label{los}
\mathbf{A}_{h,p} \mathbf{u}_{h,p} = \mathbf{f}_{h,p},
\end{eqnarray}

\noindent where $\mathbf{A}_{h,p}$ denotes the system matrix resulting from this discretization with B-spline basis functions of order $p$ and mesh width $h$. For a more detailed description of the spatial discretization in Isogeometric Analysis, we refer to \cite{hugh}. Throughout this paper four benchmarks are considered, to investigate the influence of the geometric factor, the considered coefficients in the CDR-equation, the number of patches and the dimension on the proposed $p$-multigrid solver. \\
\\
\noindent \textbf{Benchmark }$\mathbf{1.}$ Let $\Omega$ be the quarter annulus with an inner and outer radius of $1$ and $2$, respectively. The coefficients are chosen as follows:
\vskip-.6cm
\begin{eqnarray}
\mathbf{D} = \begin{bmatrix*}[r]  1 & 0           \\ 0           & 1 \\    \end{bmatrix*}, \hspace{0.5cm} \mathbf{v} = \begin{bmatrix*}[r]  0           \\ 0 \\  \end{bmatrix*}, \hspace{0.5cm} R = 0.  
\end{eqnarray}

Furthermore, homogeneous Dirichlet boundary conditions are applied and the right-hand side is chosen such that the exact solution $u$ is given by:
\vskip-.6cm
\begin{eqnarray}
u(x,y) = -(x^2 + y^2 - 1)(x^2 + y^2 -4)xy^2. \nonumber
\end{eqnarray}       

\noindent \textbf{Benchmark }$\mathbf{2.}$ Here, the unit square is adopted as domain, i.e. $\Omega = [0,1]^2$, and the coefficients are chosen as follows:
\vskip-.6cm
\begin{eqnarray}
\mathbf{D} = \begin{bmatrix*}[r]  1.2 & -0.7           \\ -0.4           & 0.9 \\    \end{bmatrix*}, \hspace{0.5cm} \mathbf{v} = \begin{bmatrix*}[r]  0.4           \\ -0.2 \\  \end{bmatrix*}, \hspace{0.5cm} R = 0.3.  
\end{eqnarray}

Homogeneous Dirichlet boundary conditions are applied and the right-hand side is chosen such that the exact solution $u$ is given by:
\vskip-.6cm
\begin{eqnarray}
u(x,y) = \text{sin}(\pi x)\text{sin}(\pi y). \nonumber
\end{eqnarray}       

\noindent \textbf{Benchmark }$\mathbf{3.}$ Let $\Omega =  \{[-1,1] \times [-1,1]\}  \backslash \{[0,1] \times [0,1] \}$ be an L-shaped domain. A multipatch geometry is created, by splitting the single patch in each direction uniformly. The coefficients are chosen as follows:
\begin{eqnarray}
\mathbf{D} = \begin{bmatrix*}[r]  1 & 0 \\ 0 & 1 \\ \end{bmatrix*}, \hspace{0.5cm} \mathbf{v} = \begin{bmatrix*}[r] 0  \\ 0 \\  \end{bmatrix*}, \hspace{0.5cm} R = 0.  
\end{eqnarray}
\newpage
The exact solution is given by: 
\vskip-.6cm
\begin{eqnarray}
u(x,y) = \begin{cases} \sqrt[3]{x^2+y^2}\text{sin}\left ( \frac{2\text{atan}2(y,x)-\pi}{3}\right ) \hspace{0.48cm} \text{if} \ y>0 \\ \sqrt[3]{x^2+y^2}\text{sin}\left ( \frac{2\text{atan}2(y,x) +3\pi}{3} \right ) \quad \text{if} \ y<0 \end{cases}, \nonumber
\end{eqnarray}       

and the right-hand side is chosen accordingly. Inhomogeneous Dirichlet boundary conditions are prescribed for this benchmark. \\
\\
\noindent \textbf{Benchmark }$\mathbf{4.}$ Here, the unit cube is adopted as domain, i.e. $\Omega = [0,1]^3$, and the coefficients are chosen as follows:
\vskip-.6cm
\begin{eqnarray}
\mathbf{D} = \begin{bmatrix*}[r]  1 & 0 & 0 \\ 0  & 1  & 0 \\ 0 & 0 & 1 \end{bmatrix*}, \hspace{0.5cm} \mathbf{v} = \begin{bmatrix*}[r]  0 \\ 0 \\ 0 \\  \end{bmatrix*}, \hspace{0.5cm} R = 0.  
\end{eqnarray}

Homogeneous Dirichlet boundary conditions are applied and the right-hand side is chosen such that the exact solution $u$ is given by:
\vskip-.6cm
\begin{eqnarray}
u(x,y) = \text{sin}(\pi x)\text{sin}(\pi y)\text{sin}(\pi z). \nonumber
\end{eqnarray}     

\section{p-Multigrid method}
\label{pm}
Multigrid methods \cite{brandt,hackbush} aim to solve linear systems of equations by defining a hierarchy of discretizations. At each level of the multigrid hierarchy
a smoother is applied, whereas on the coarsest level a correction is determined by means of a direct solver. Starting from $\mathcal{V}_{h,1}$, a sequence of spaces $\mathcal{V}_{h,1}, \ldots , \mathcal{V}_{h,p}$ is obtained by applying $k$-refinement to solve Equation \eqref{los}. Note that, since basis functions with maximal continuity are considered, the spaces are not nested.   
\noindent A single step of the two-grid correction scheme for the $p$-multigrid method consists of the following steps \cite{roel}:
\begin{enumerate}
\item Starting from an initial guess $\mathbf{u}_{h,p}^{(0,0)}$, apply a f\/ixed number $\nu_1$ of pre-smoothing steps:
\vskip-.6cm
\begin{eqnarray} \label{smooth}
\mathbf{u}_{h,p}^{(0,m)} = \mathbf{u}_{h,p}^{(0,m-1)} + \mathcal{S}_{h,p} \left (\mathbf{f}_{h,p} - \mathbf{A}_{h,p} \mathbf{u}_{h,p}^{(0,m-1)} \right ),  \hspace{0.25cm} m=1,\ldots, \nu_1, 
\end{eqnarray}
where $\mathcal{S}_{h,p}$ is a smoothing operator applied to the high-order problem. 
\item Determine the residual  at level $p$ and project it onto the space $\mathcal{V}_{h,p-1}$ using the restriction operator $\mathcal{I}_{p}^{p-1}$:
\vskip-.6cm
\begin{eqnarray}
\mathbf{r}_{h,p-1} = \mathcal{I}_{p}^{p-1} \left (  \mathbf{f}_{h,p} - \mathbf{A}_{h,p} \mathbf{u}_{h,p}^{(0,\nu_1)} \right ).
\end{eqnarray} 
\item Solve the residual equation at level $p-1$ to determine the coarse grid error: 
\vskip-.6cm
\begin{eqnarray} \label{eq:res}
\mathbf{A}_{h,p-1} \mathbf{e}_{h,p-1} = \mathbf{r}_{h,p-1}. 
\end{eqnarray}
\item Project the error $\mathbf{e}_{h,p-1}$ onto the space $\mathcal{V}_{h,p}$ using the prolongation operator $\mathcal{I}_{p-1}^p$ and update $\mathbf{u}_{h,p}^{(0,\nu_1)}$:
\vskip-.6cm
\begin{eqnarray}
\mathbf{u}_{h,p}^{(0,\nu_1)} := \mathbf{u}_{h,p}^{(0,\nu_1)} + \mathcal{I}_{p-1}^p \left (\mathbf{e}_{h,p-1} \right ).
\end{eqnarray}
\item Apply $\nu_2$ postsmoothing steps of~\eqref{smooth} to obtain $\mathbf{u}_{h,p}^{(0,\nu_1 + \nu_2)} =: \mathbf{u}_{h,p}^{(1,0)}$.

\end{enumerate} 

\noindent In the literature, steps (2)-(4) are referred to as 'coarse grid correction'. Recursive application of this scheme on Equation \eqref{eq:res} until level $p=1$ is reached, results in a V-cycle. In contrast to $h$-multigrid methods, the coarsest problem in $p$-multigrid can still be large for small values of $h$. However, since we restrict to level $p=1$, the coarse grid problem corresponds to a standard low-order Lagrange FEM discretization of the problem at hand. Therefore, we use a standard $h$-multigrid method to solve the coarse grid problem in our $p$-multigrid scheme, which is known to be optimal (in particular $h$-independent) in this case. As a smoother, Gauss-Seidel is applied within the $h$-multigrid method, as it's both cheap and effective for low degree problems.   
Applying two V-cycles using canonical prolongation, weighted restriction and a single smoothing step turned out to be sufficient and has therefore been adopted throughout this paper as coarse grid solver. 

\begin{figure}[h!]
\begin{center}
\begin{tikzpicture}[xscale = 0.67, yscale = 0.62]
\draw (-14,4.5) node {\tiny $\textcolor{blue}{p=3}$};
\draw (-2,4.5) node {\tiny $\textcolor{red}{h = 2^{-5}}$};
\draw (-14, 2) node {\tiny $\textcolor{blue}{p=3}$};
\draw (-2, 2) node {\tiny $\textcolor{red}{h = 2^{-4}}$};  
\draw (-14,-0.5) node {\tiny $\textcolor{blue}{p=3}$};
\draw (-2, -0.5) node {\tiny $\textcolor{red}{h = 2^{-3}}$};
\draw (0, 2) node {\Bigg \} \tiny $h$-multigrid};
\draw (-15, 2) node { \Bigg\{ };
\draw (-16, 2) node {\tiny IgA};

\draw [->,thick,blue] (-10.85,4.25) --  node[sloped, anchor=center, below]{}(-9.62,2.2);
\draw [->,thick,blue] (-9.45,1.85) --  node[sloped, anchor=center, below]{}(-8.12,-0.3);
\draw [->,thick,blue] (-7.95,-0.3)  --  node[sloped, anchor=center, below]{}(-6.8,1.85);
\draw [->,thick,blue] (-6.65,2.2)  --  node[sloped, anchor=center, below]{}(-5.55,4.3);

\draw[dashed] (-12, 4.5) -- (-4, 4.5);
\draw[dashed] (-12, 2) -- (-4, 2);
\draw[dashed] (-12, -0.5) -- (-4, -0.5);

\node[scale=1.35] at (-11,4.5) {\pgfuseplotmark{triangle*}};
\node[scale=1.35] at (-9.5,2) {\pgfuseplotmark{triangle*}};
\node at (-8,-0.5) {\pgfuseplotmark{square*}};
\node[scale=1.35] at (-6.75,2) {\pgfuseplotmark{triangle*}};
\node[scale=1.35] at (-5.45,4.5) {\pgfuseplotmark{triangle*}};

\end{tikzpicture}
\caption{Illustration of the considered $p$-multigrid (top) and $h$-multigrid (bottom) method. At $p=1$, Gauss-Seidel is always adopted as a smoother ($\bullet$), whereas at the high order level Gauss-Seidel or ILUT can be applied ($\blacktriangle$). At the coarsest level, a direct solver is applied to solve the residual equation ($\blacksquare$).}
\label{fig:cycle}
\end{center}
\end{figure}
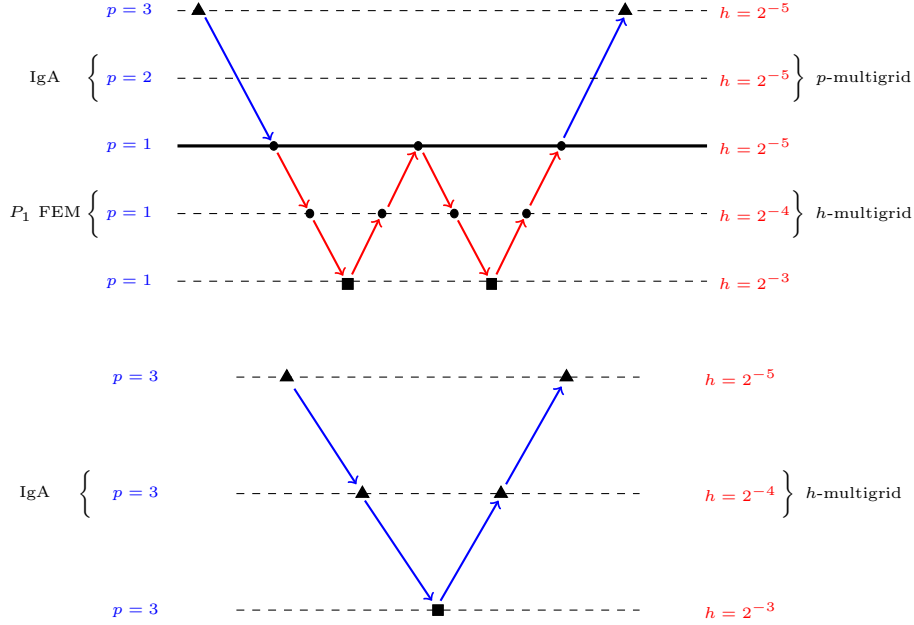

\noindent Note that, for the $p$-multigrid method, the residual can be projected directly to level $p=1$. It was shown in \cite{tielen_enumath} that the preformance of the $p$-multigrid method is hardly affected, while the set-up costs decrease significantly. In \ref{direct}, numerical results are presented for the first benchmark confirming this observation. Throughout this paper, a direct projection to level $p=1$ is adopted for the $p$-multigrid method, see Figure \ref{fig:cycle}. Results are compared to an $h$-multigrid method, which is shown as well in Figure \ref{fig:cycle}.

\begin{figure}[h!]
\begin{center}
\begin{tikzpicture}[xscale = 0.32, yscale = 0.36]
\draw (-14,4.5) node {\tiny $\textcolor{blue}{p=3}$};
\draw (12,4.5) node {\tiny $\textcolor{red}{h = 2^{-5}}$};
\draw (-14, 2) node {\tiny $\textcolor{blue}{p=2}$};
\draw (12, 2) node {\tiny $\textcolor{red}{h = 2^{-5}}$};  
\draw (-14,-0.5) node {\tiny $\textcolor{blue}{p=1}$};
\draw (12, -0.5) node {\tiny $\textcolor{red}{h = 2^{-5}}$};
\draw (-14,-3.0) node {\tiny $\textcolor{blue}{p=1}$};
\draw (12, -3.0) node {\tiny $\textcolor{red}{h = 2^{-4}}$};
\draw (-14,-5.5) node {\tiny $\textcolor{blue}{p=1}$};
\draw (12, -5.5) node {\tiny $\textcolor{red}{h = 2^{-3}}$};
\draw (16, 2) node {\Bigg \} \tiny $p$-multigrid};
\draw (16, -3) node {\Bigg \} \tiny $h$-multigrid};
\draw (-15.5, 2) node { \Bigg\{ };
\draw (-15.5, -3) node {\Bigg \{};
\draw (-17.5, 2) node {\tiny IgA};
\draw (-17.5, -3) node {\tiny $P_1$ FEM};

\draw [->,thick,blue] (-10.85,4.25) --  node[sloped, anchor=center, below]{}(-8.12,-0.3);
\draw [->,thick,red] (-7.85,-0.75) --  node[sloped, anchor=center, below]{}(-6.62,-2.8);
\draw [->,thick,red] (-6.35,-3.25) -- node[sloped, anchor=center, below]{}(-5.12,-5.3);
\draw [->,thick,red] (-4.75,-5.25) -- node[sloped, anchor=center, below]{} (-3.62,-3.2);
\draw [->,thick,red] (-3.25,-2.75) --  node[sloped, anchor=center, below]{}(-2.12,-0.7);

\draw [->,thick,red] (-1.80,-0.75) --  node[sloped, anchor=center, below]{}(-0.62,-2.8);
\draw [->,thick,red] (-0.35,-3.25) -- node[sloped, anchor=center, below]{}(0.88,-5.3);
\draw [->,thick,red] (1.25,-5.25) -- node[sloped, anchor=center, below]{} (2.38,-3.2);
\draw [->,thick,red] (2.65,-2.7) --  node[sloped, anchor=center, below]{}(3.80,-0.7);
\draw [->,thick,blue] (4.05,-0.2)  --  node[sloped, anchor=center, below]{}(6.55,4.25);

\draw[dashed] (-12, 4.5) -- (10, 4.5);
\draw[dashed] (-12, 2) -- (10, 2);
\draw[very thick] (-12, -0.5) -- (10, -0.5);
\draw[dashed] (-12, -3.0) -- (10, -3.0);
\draw[dashed] (-12, -5.5) -- (10, -5.5);

\node[scale=1.35] at (-11.125,4.5) {\pgfuseplotmark{triangle*}};
\filldraw (-8,-0.5) circle (4.5pt) node[]{};
\draw[black, fill] (-5.15,-5.8) rectangle (-4.70,-5.4);
\filldraw (-6.5,-3.0) circle (4.5pt) node[]{};
\filldraw (-3.5,-3.0) circle (4.5pt) node[]{};
\filldraw (-2,-0.5) circle (4.5pt) node[]{};
\draw[black, fill] (0.85,-5.8) rectangle (1.25,-5.4);
\filldraw (-0.5,-3.0) circle (4.5pt) node[]{};
\filldraw (2.5,-3.0) circle (4.5pt) node[]{};
\filldraw (3.95,-0.5) circle (4.5pt) node[]{};
\node[scale=1.35] at (6.60,4.5) {\pgfuseplotmark{triangle*}};
\end{tikzpicture}
\end{center}
\end{figure}

\subsection*{Prolongation and restriction}
To transfer both coarse grid corrections and residuals between different levels of the multigrid hierarchy, prolongation and restriction operators are defined.
The prolongation and restriction operator adopted in this paper are based on an $L_2$ projection and have been used extensively in the literature \cite{scot, brig, oct}.  
At level $p=1$, the coarse grid correction is prolongated to level $p$ by projection onto the space $\mathcal{V}_{h,p}$. The prolongation operator $\mathcal{I}_{1}^p: \mathcal{V}_{h,1} \rightarrow \mathcal{V}_{h,p}$ is given by
\vskip-.6cm
\begin{eqnarray} \label{prolongation}
\mathcal{I}_{1}^p(\mathbf{v}_{p}) = (\mathbf{M}_p)^{-1} \mathbf{P}_{1}^{p} \ \mathbf{v}_{1}, 
\end{eqnarray}
where the mass matrix $\mathbf{M}_p$ and transfer matrix $\mathbf{P}_{1}^{p}$ are defined, respectively, as follows:
\vskip-.6cm
\begin{eqnarray} \label{PM}
 ( \mathbf{M}_{p})_{(i,j)} := \int_{\Omega} \phi_{i,p} \phi_{j,p}  \hspace{0.1cm} \text{d} \Omega,  \hspace{1.0cm} (\mathbf{P}^{p}_{1})_{(i,j)} :=  \int_{\Omega} \phi_{i,p} \phi_{j,1}  \hspace{0.1cm} \text{d} \Omega
\end{eqnarray}

\noindent The residuals are restricted from level $p$ to $1$ by projection onto the space $\mathcal{V}_{h,1}$. The restriction operator  $\mathcal{I}_{p}^{1}: \mathcal{V}_{h,p} \rightarrow \mathcal{V}_{h,1}$ is defined by
\vskip-.6cm
\begin{eqnarray} \label{restriction}
\mathcal{I}_{p}^{1}(\mathbf{v}_{p}) = (\mathbf{M}_{1})^{-1} \mathbf{P}_{p}^{{1}} \ \mathbf{v}_{p}.
\end{eqnarray}
To prevent the explicit solution of a linear system of equations for each projection step, the consistent mass matrix $\mathbf{M}$ in both transfer operators is replaced by its lumped counterpart $\mathbf{M}^L$ by applying row-sum lumping:  
\vskip-.6cm
\begin{eqnarray}
\mathbf{M}_{(i,i)}^L = \sum_{j=1}^{N_{\rm dof}} \mathbf{M}_{(i,j)} .
\end{eqnarray}

Numerical experiments, presented in \ref{lumped}, show that lumping the mass matrix hardly influences the convergence behaviour of the resulting $p$-multigrid method. Neither does it affect the overall accuracy obtained with the $p$-multigrid method. Alternatively, one could invert the mass matrix efficiently by exploiting the tensor product structure, see \cite{gao}.

\noindent Note that this choice of prolongation and restriction operators yields a non-symmetric coarse grid correction and, hence, a non-symmetric multigrid solver. As a consequence, the multigrid solver can only be applied as a preconditioner for a Krylov method suited for non-symmetric matrices, like BiCGSTAB. \\
\\
\noindent \textbf{Remark 1:} Choosing the prolongation and restriction operator tranpose to each other would restore the symmetry of the multigrid method. However, numerical experiments, not presented in this paper, show that this leads to a less robust $p$-multigrid method. Therefore, the prolongation and restriction operator are adopted as defined in Equation \eqref{prolongation} and \eqref{restriction}, respectively.    

\subsection*{Smoother}
Within multigrid methods, a basic iterative method is typically used as a smoother. However, in IgA the performance of classical smoothers such as (damped) Jacobi and Gauss-Seidel decreases significantly for higher values of $p$. Therefore, an Incomplete LU Factorization is adopted with a dual threshold strategy (ILUT) \cite{saad} to approximate the operator $\mathbf{A}_{h,p}$:
\vskip-.6cm
\begin{eqnarray}
\mathbf{A}_{h,p} \approx \mathbf{L}_{h,p} \mathbf{U}_{h,p}.
\end{eqnarray}

\noindent The ILUT factorization is determined completely by a tolerance $\tau$  and fillfactor $m$. Two dropping rules are applied during factorization:
\vskip-.6cm
\begin{enumerate}
\item All elements smaller (in absolute value) than the dropping tolerance are dropped. The dropping tolerance is obtained by multiplying the tolerance $\tau$ with the average magnitude of all elements in the current row. 
\item Apart from the diagonal element, only the $M$ largest elements are kept in each row. Here, $M$ is determined by multiplying the fillfactor $m$ with the average number of non-zeros in each row of the original operator $\mathbf{A}_{h,p}$.
\end{enumerate}   

\noindent The ILUT factorization considered in this paper is closely related to an ILU(0) factorization. This factorization has been applied in the context of IgA as a preconditioner, showing good convergence behaviour \cite{collier}.

\begin{figure}[h!]
\begin{center}
\includegraphics[scale=0.31]{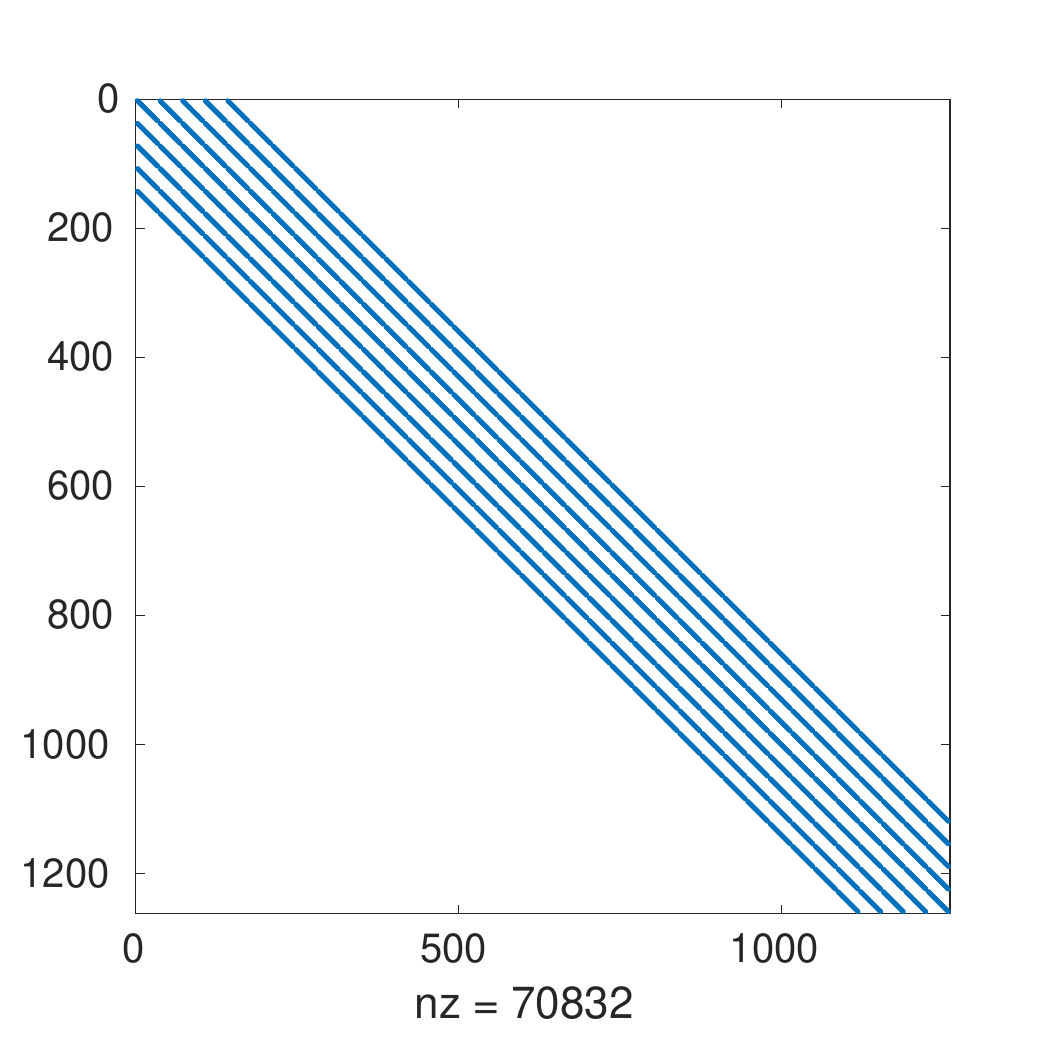}
\includegraphics[scale=0.31]{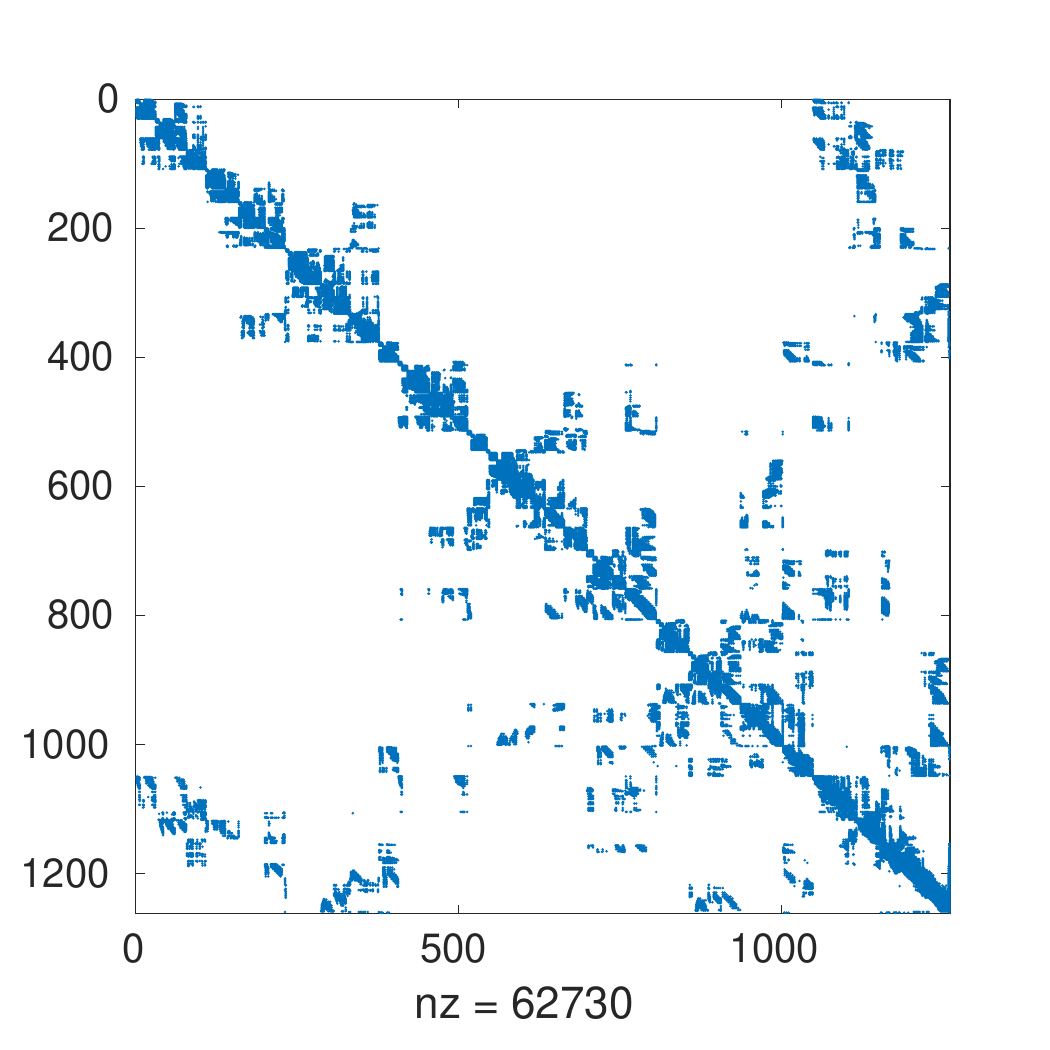}
\end{center} 
\caption{Sparsity pattern of $\mathbf{A}_{h,3}$ (left) and $\mathbf{L}_{h,3} + \mathbf{U}_{h,3}$ (right) for $h=2^{-5}$.}
\label{fig:2}
\end{figure}

\noindent An efficient implementation of ILUT is available in the Eigen library \cite{eigen} based on \cite{saad2}. Once the factorization is obtained, a single smoothing step is applied as follows:
\vskip-.6cm
\begin{eqnarray}
 \mathbf{e}_{h,p}^{(n)} &=& (\mathbf{L}_{h,p} \mathbf{U}_{h,p})^{-1} (\mathbf{f}_{h,p} - \mathbf{A}_{h,p} \mathbf{u}_{h,p}^{(n)}),  \\
                          &=&  \mathbf{U}_{h,p}^{-1}  \mathbf{L}_{h,p}^{-1} (\mathbf{f}_{h,p} - \mathbf{A}_{h,p} \mathbf{u}_{h,p}^{(n)}), \label{inverse} \\
 \mathbf{u}_{h,p}^{(n+1)} &=& \mathbf{u}_{h,p}^{(n)} + \mathbf{e}_{h,p}^{(n)}, 
\end{eqnarray}

\noindent where the two matrix inversions in Equation \eqref{inverse} amount to forward and backward substitution. Throughout this paper, the fillfactor $m = 1$ is used (unless stated otherwise) and the dropping tolerance equals $\tau = 10^{-12}$. Hence, the number of non-zero elements of $\mathbf{L}_{h,p} + \mathbf{U}_{h,p}$ is similar to the number of non-zero elements of $\mathbf{A}_{h,p}$. 
\noindent Figure \ref{fig:2} shows the sparsity pattern of the stiffness matrix $\mathbf{A}_{h,3}$ and $\mathbf{L}_{h,3} + \mathbf{U}_{h,3}$ for the first benchmark and $h=2^{-5}$. Since a fill-reducing permutation is performed during the ILUT factorization, sparsity patterns differ significantly. However, the number of non-zero entries is comparable. 

\subsection*{Coarse grid operator}
At each level of the multigrid hierarchy, the operator $\mathbf{A}_{h,p}$ is needed to apply smoothing or compute the residual. The operators at the coarser levels can be obtained by rediscretizing the bilinear form in \eqref{weak} with low-order basis functions. Alternatively, a Galerkin projection can be adopted: 
\vskip-.6cm
\begin{eqnarray}
\mathbf{A}_{h,p-1}^{G} = \mathcal{I}_{p}^{p-1} \ \mathbf{A}_{h,p} \ \mathcal{I}_{p-1}^{p}. 
\end{eqnarray}

\noindent However, since the condition number when using the Galerkin projection is significantly higher compared to the condition number of the rediscretized operator, the coarse grid operators in this paper are obtained by using the rediscretization approach.  

\subsection*{Computational costs}

To investigate the costs of the proposed $p$-multigrid method, the assembly costs, factorization costs and the costs of a single V-cycle are analyzed. Assuming a direct projection to level $p=1$, both $\mathbf{A}_{h,p}$ and $\mathbf{A}_{h,1}$ have to be assembled. Furthermore, the (variationally lumped) mass matrices $\mathbf{M^L}_{1}$, $\mathbf{M^L}_{p}$ and transfer matrix $\mathbf{P}_{p}^{1}$ have to be assembled.

\noindent Assuming an element-based assembly loop with standard Gauss-quadrature, assembling the stiffness matrix or transfer matrix at level $k$ costs $\mathcal{O}(N_{\rm dof} k^{3d})$ f\/lops. More efficient assembly techniques like weighted quadrature \cite{sangalli}, sum factorizations \cite{antolin} and low tensor rank approximations \cite{mantzaflaris} exist, but have not yet been explored. However, since assembly costs might dominate the overall computational costs, assembly costs will be presented separately in Section \ref{numex}. Assembling the (variationally lumped) mass matrices $\mathbf{M^L}_{1}$ and $\mathbf{M^L}_{p}$ costs $\mathcal{O}(N_{\rm dof})$ flops.  

\noindent At the high order level an ILUT factorization of $\mathbf{A}_{h,p}$ needs to be determined, costing $\mathcal{O}(N_{\rm dof} p^{2d}$) flops \cite{collier} in case $m=1$ and $\tau = 10^{-12}$. Alternatively to ILUT, Gauss-Seidel can be applied as a smoother, without any set-up costs. 

\noindent At the high order level of the V-cycle both pre- and postsmoothing is applied. Given the ILUT factorization, applying a single smoothing step costs $\mathcal{O}(N_{\rm dof} p^{d})$ flops. Applying Gauss-Seidel as a smoother at level $p=1$, costs $\mathcal{O}(N_{\rm dof} 1^{d})$ f\/lops \cite{collier}. For both prolongation and restriction, a sparse matrix-vector multiplication has to be performed, costing $\mathcal{O}(N_{\rm dof}p^{d})$ flops for each application. 

\noindent Finally, the residual equation \eqref{eq:res} is solved by applying a single V-cycle of an $h$-multigrid method, which uses Gauss-Seidel as a smoother. Prolongation and restriction operators of the $h$-multigrid method are based on canonical interpolation and weighted restriction, respectively. Table \ref{tab:cc} provides an overview of the computational costs of the proposed $p$-multigrid method. 

\begin{table}[h]
\centering
\caption{Total computational costs with $p$-multigrid for general values of the approximation order $p$ and dimension $d$.}
\begin{tabular}{c|c|c|c|c} \hline
 \multicolumn{5}{c}{Setup costs} \\ \hline
 Level $k$ & Assembly  $\mathbf{A}_{h,k}$    & Assembly $\mathbf{M^L}_{k}$         & Assembly $\mathbf{P}_{k}^{k-1}$   & ILU factorization \\
 $1$       &$\mathcal{O}(N_{\rm dof} 1^{3d}$) &  $\mathcal{O}(N_{\rm dof}$)& 								  & $\mathcal{O}(N_{\rm dof} 1^{2d}$) \\
 $p$       &$\mathcal{O}(N_{\rm dof} p^{3d}$) &  $\mathcal{O}(N_{\rm dof}$)& $\mathcal{O}(N_{\rm dof} p^{3d}$) & $\mathcal{O}(N_{\rm dof} p^{2d}$) \\ \hline
 \multicolumn{5}{c}{Costs V-cycle} \\ \hline
Level $k$ & Presmoothing    		        & Restriction      				& Prolongation  				& Postsmoothing \\
 $1$   & - 			  				    & 			-	                & $\mathcal{O}(N_{\rm dof}p^d)$	&  - \\
 $p$   & $\mathcal{O}(N_{\rm dof} p^{d}$)	&  $\mathcal{O}(N_{\rm dof}p^d)$&			-   	   			& $\mathcal{O}(N_{\rm dof} p^{d})$ \\ \hline
\end{tabular}
\label{tab:cc}
\end{table}

\noindent The memory requirements of the proposed $p$-multigrid method is strongly related to the number of nonzero entries of each operator. For the stiffness matrix in $d$ dimensions, the number of nonzero entries at level $k$ equals $\mathcal{O}(N_{\rm dof} k^{d}$). Table \ref{tab:nnz} shows the number of nonzero entries for all operators in the $p$-multigrid method for each level.

\begin{table}[h]
\centering
\caption{Memory requirements with $p$-multigrid for general values of the approximation order $p$ and dimension $d$.}
\begin{tabular}{c|c|c|c|c} \hline
 \multicolumn{5}{c}{Number of nonzero entries} \\ \hline
 Level $k$ &  $\mathbf{A}_{h,k}$             &   $\mathbf{M}_{k}$                      &   $\mathbf{P}_{k}^{k-1}$         &  ILU factorization \\
 $1$       &$\mathcal{O}(N_{\rm dof} 1^{d}$) &  $\mathcal{O}(N_{\rm dof}$)& $\mathcal{O}(N_{\rm dof} 1^{d}$) & $\mathcal{O}(N_{\rm dof} 1^{d}$) \\
 $p$       &$\mathcal{O}(N_{\rm dof} p^{d}$) &  $\mathcal{O}(N_{\rm dof}$)& $\mathcal{O}(N_{\rm dof} p^{d}$) & $\mathcal{O}(N_{\rm dof} p^{d}$) \\ \hline
\end{tabular}
\label{tab:nnz}
\end{table}

\noindent Note that, compared to $h$-multigrid methods, the $p$-multigrid method consists of one extra level. Since $k$-refinement is applied, the dimensions of the matrix remain of $\mathcal{O}(N_{\rm dof}$). However, at level $p=1$, coarsening in $h$ is applied which leads to a reduction of the number of degrees of freedom with a factor of $2^d$ from one level to the other, as with $h$-multigrid. Furthermore, the number of nonzero entries significantly reduces due to the smaller support of the piecewise linear B-spline basis functions. A more detailed comparison between $h$-multigrid and $p$-multigrid methods, also in terms of CPU times, can be found in Section \ref{num} and \ref{numex}.  
 
\newpage 

\section{Spectral Analysis}
\label{num}

In this section, the performance of the proposed $p$-multigrid method is analyzed and compared with $h$-multigrid methods in different ways. First, a spectral analysis is performed to investigate the interplay between the coarse grid correction and the smoother. In particular, we compare both smoothers (Gauss-Seidel and ILUT) and coarsening strategies (in $h$ or $p$). Then, the spectral radius of the iteration matrix is determined to obtain the asymptotic convergence factors of the $p$-multigrid and $h$-multigrid methods. Throughout this section, the first two benchmarks presented in Section \ref{mp} are considered for the analysis.     

\subsection*{Reduction factors}
\label{spa}

\begin{figure}[h!]
\begin{center}
\includegraphics[scale=0.31]{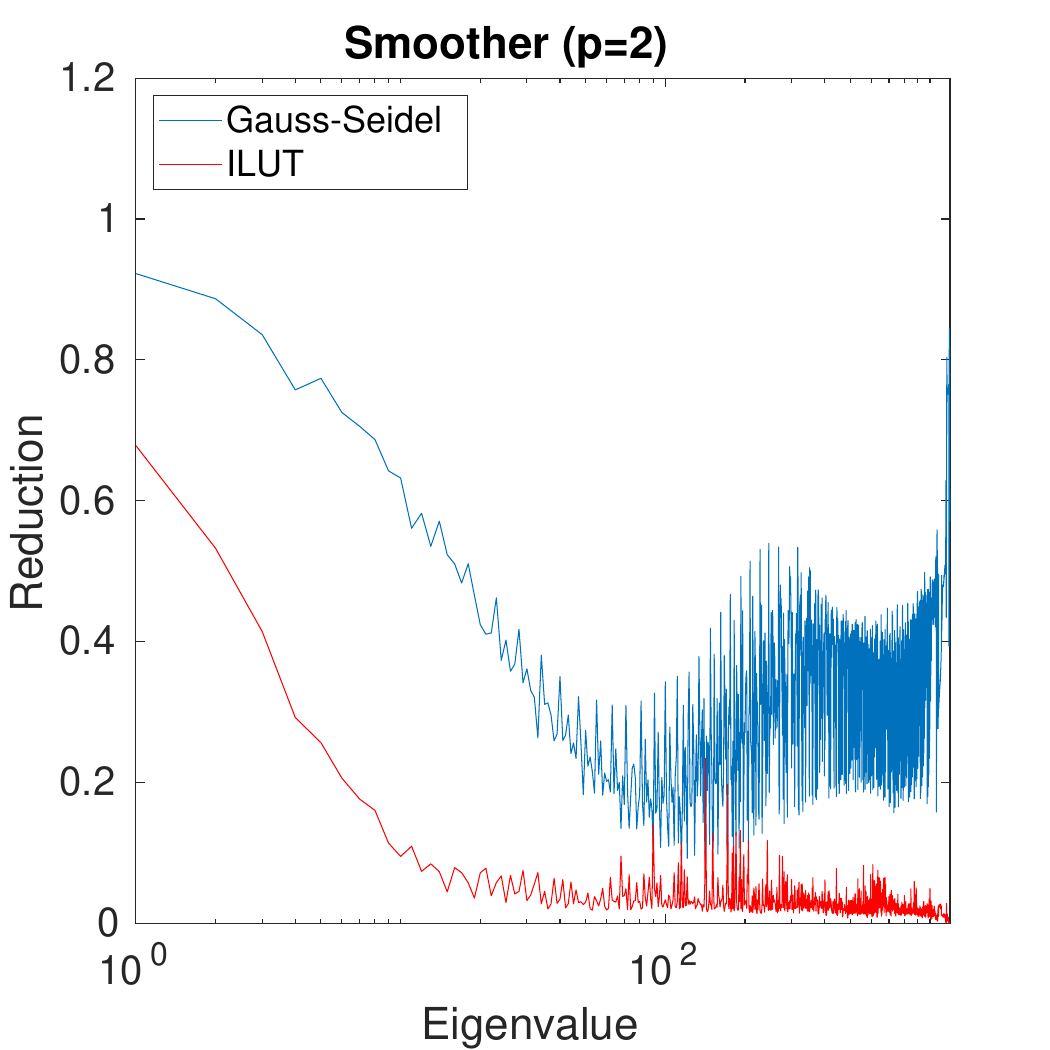}
\includegraphics[scale=0.31]{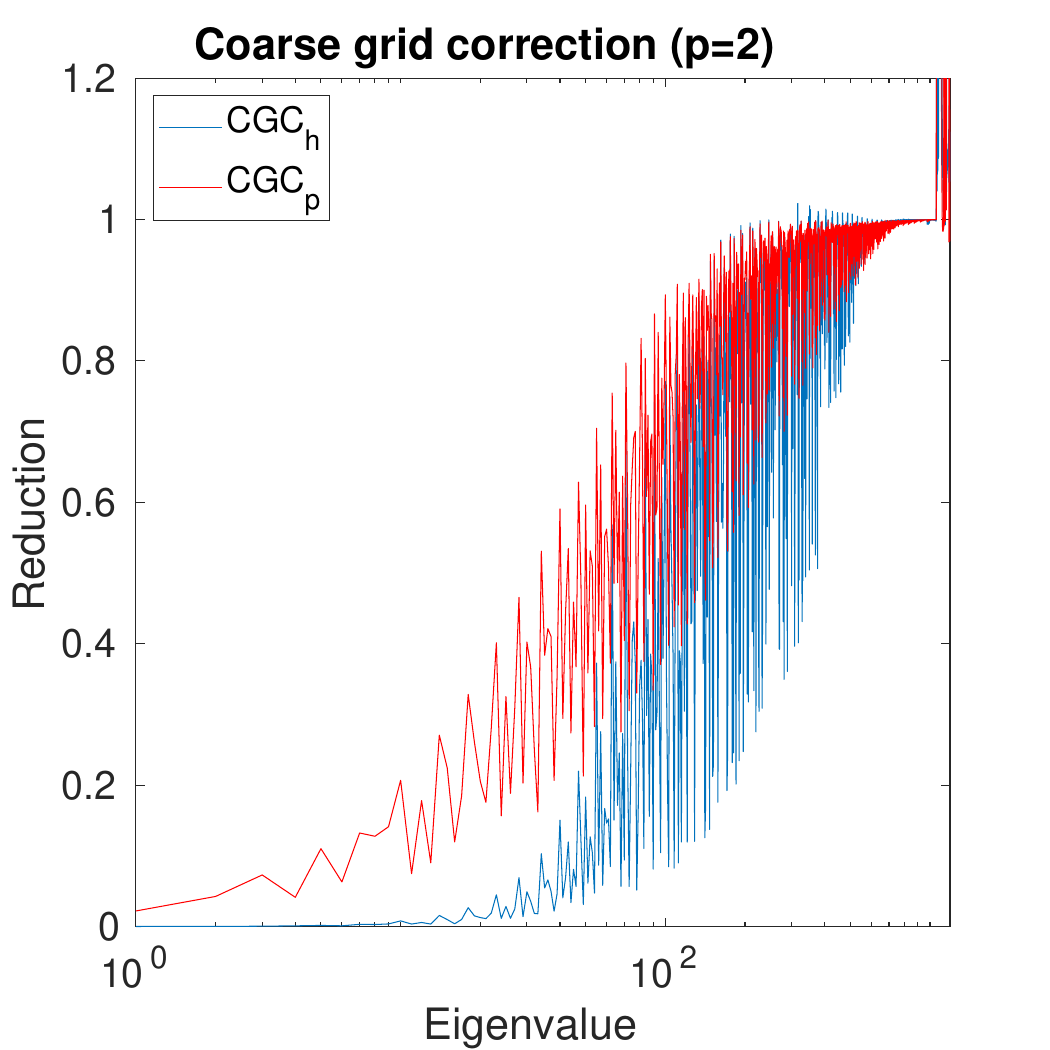}
\includegraphics[scale=0.31]{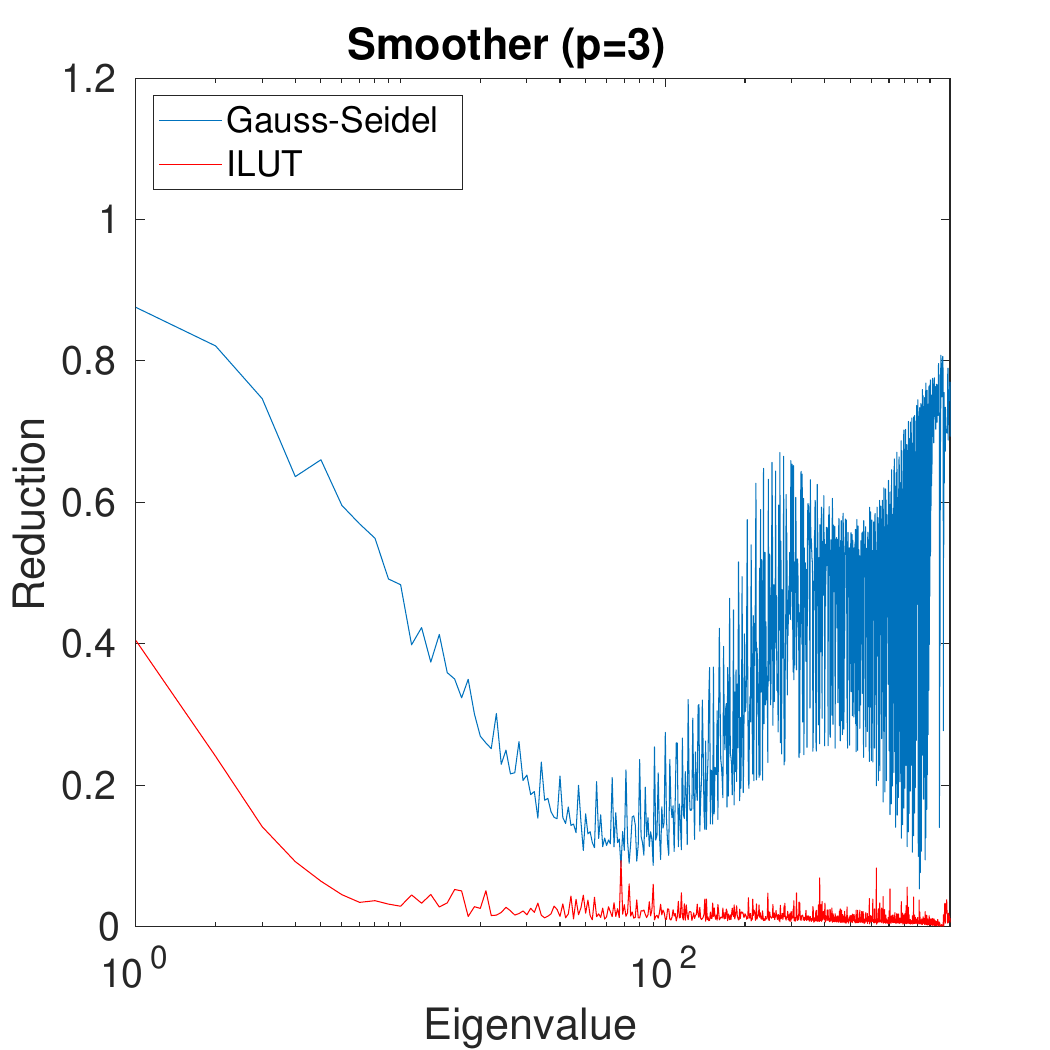}
\includegraphics[scale=0.31]{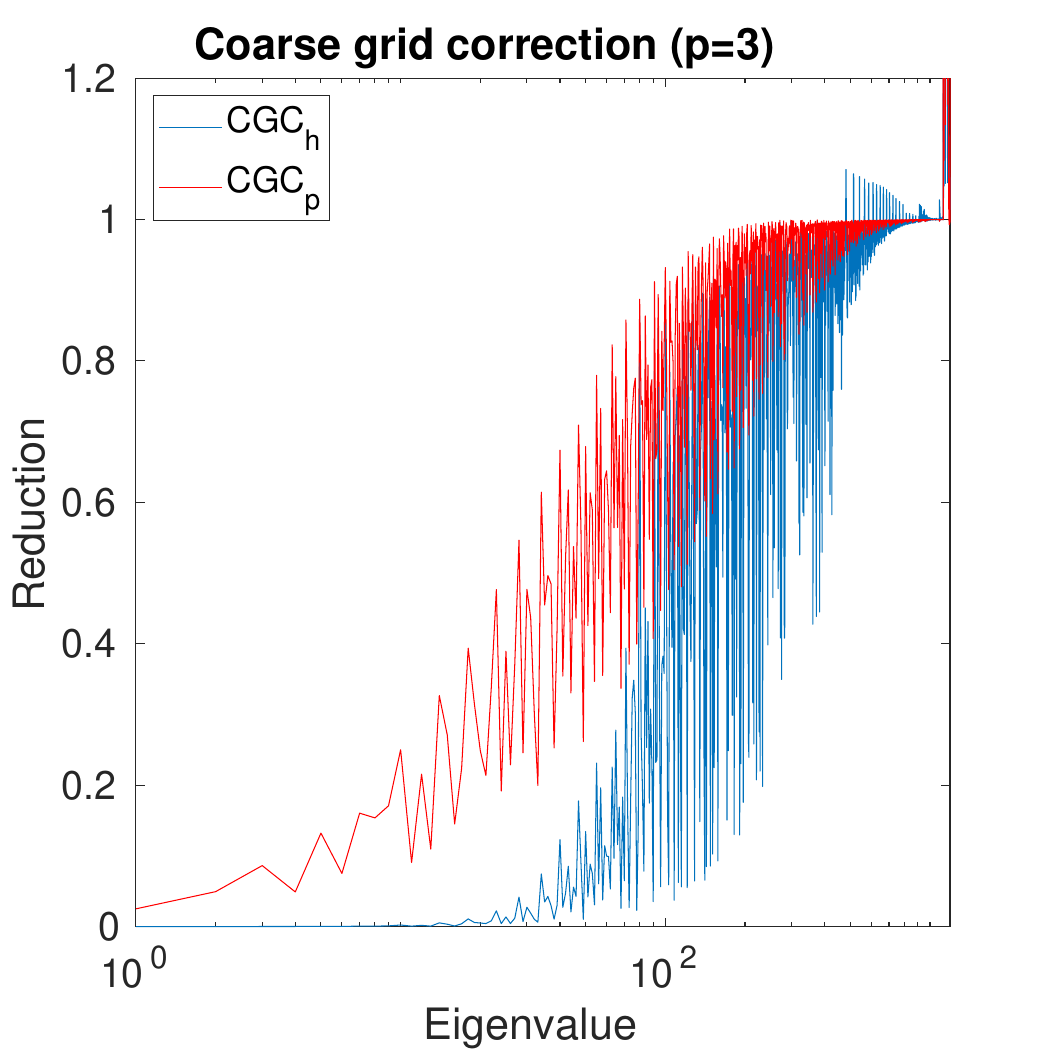}
\includegraphics[scale=0.31]{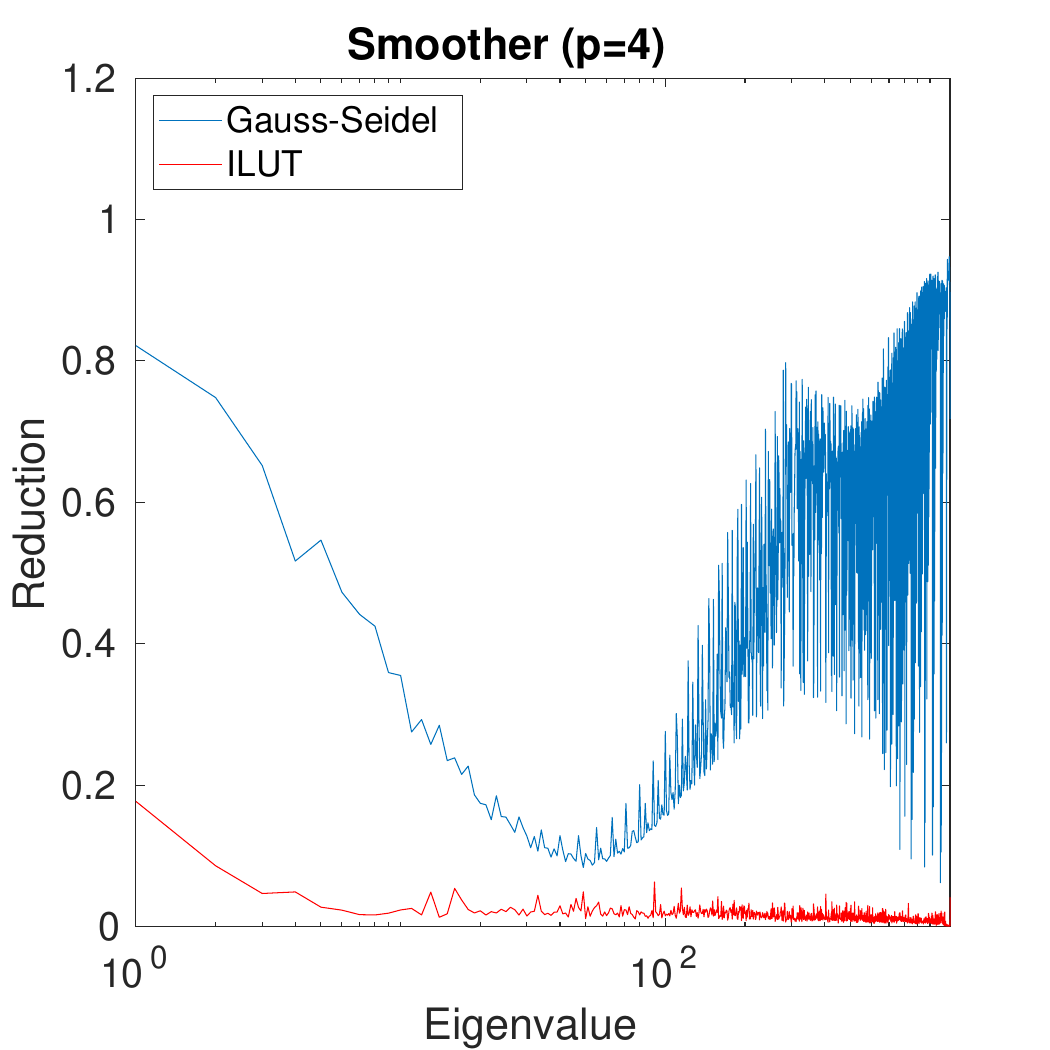}
\includegraphics[scale=0.31]{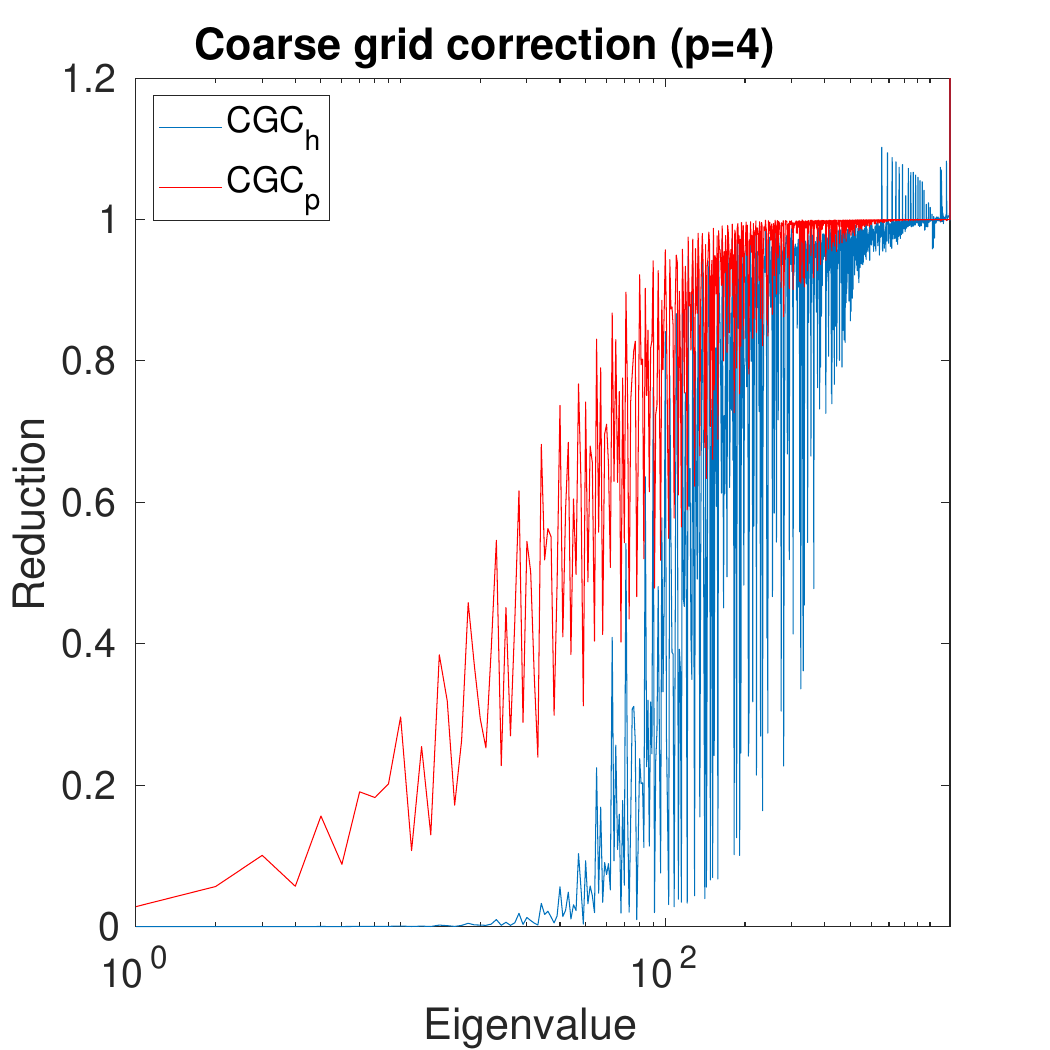}
\end{center} 
\caption{Error reduction in $(\mathbf{v}_j)$ for the first benchmark with $p=2,3,4$ and $h=2^{-5}$ obtained with different smoothers (left) and coarsening strategies (right).}
\label{fig:3}
\end{figure}

\begin{figure}[h!]
\begin{center}
\includegraphics[scale=0.31]{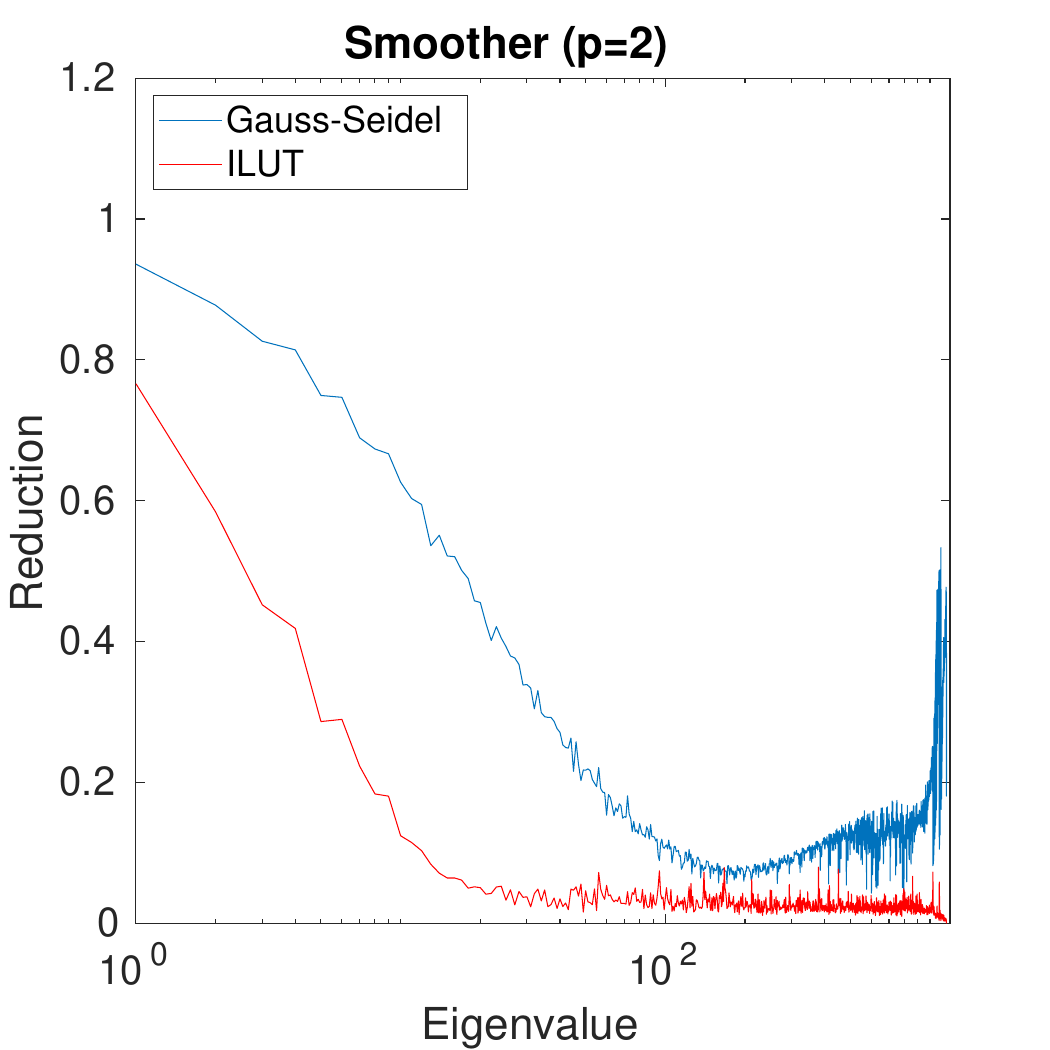}
\includegraphics[scale=0.31]{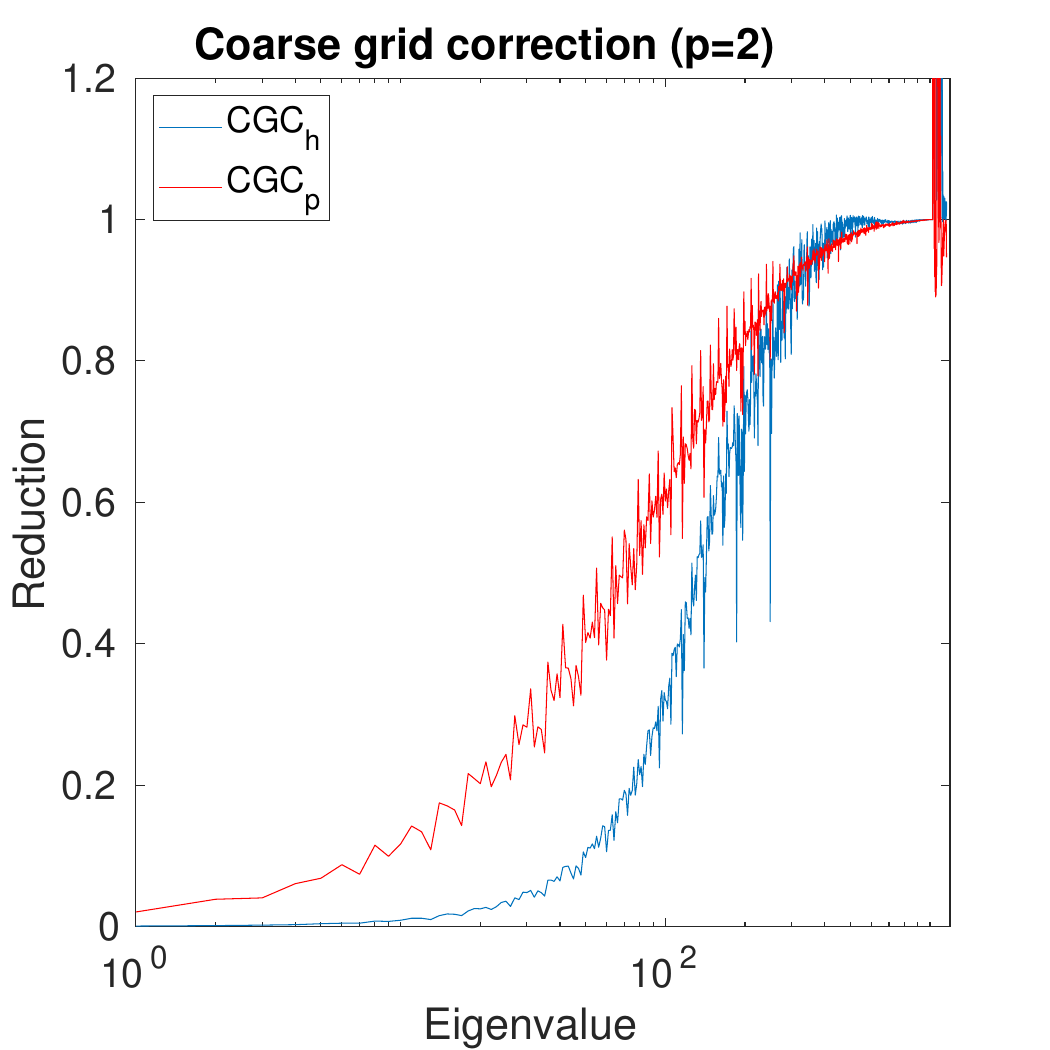}
\includegraphics[scale=0.31]{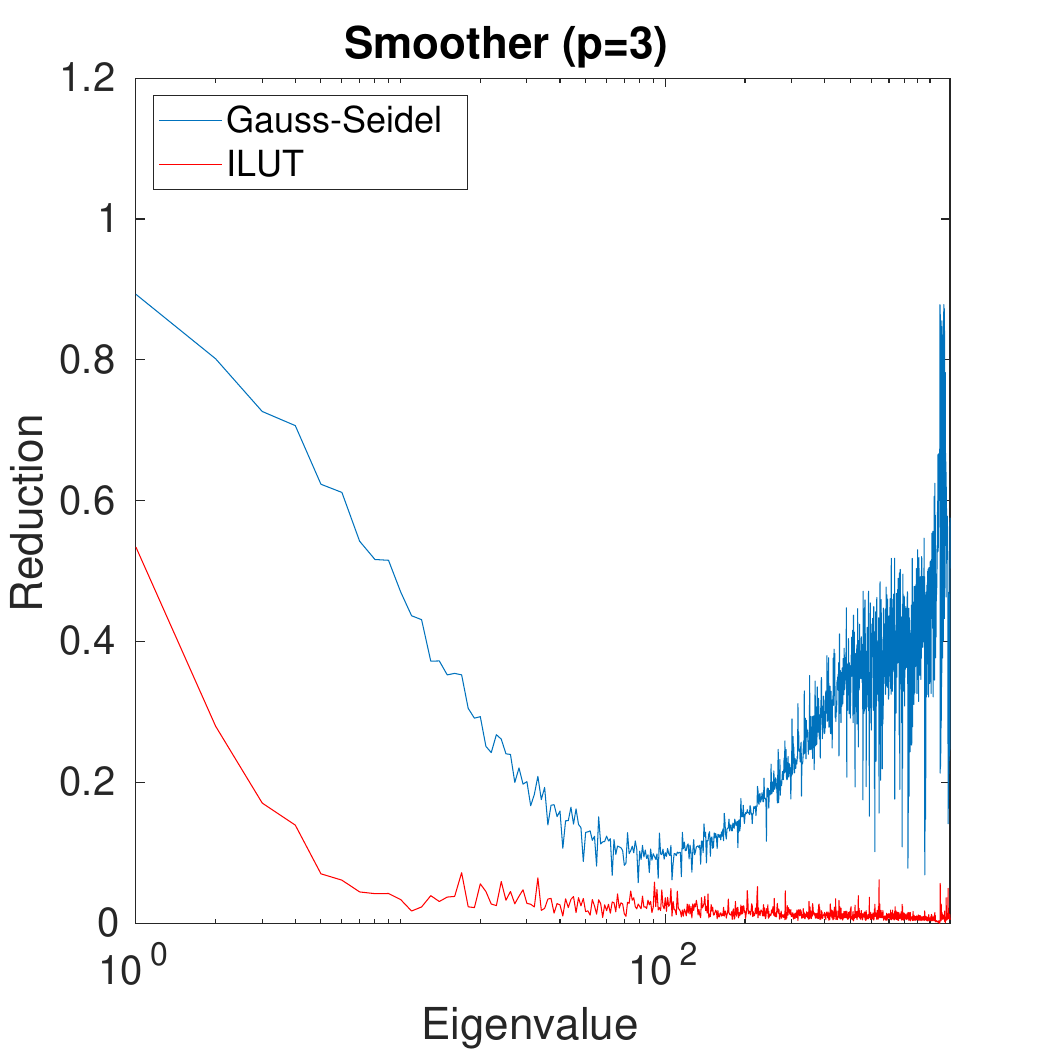}
\includegraphics[scale=0.31]{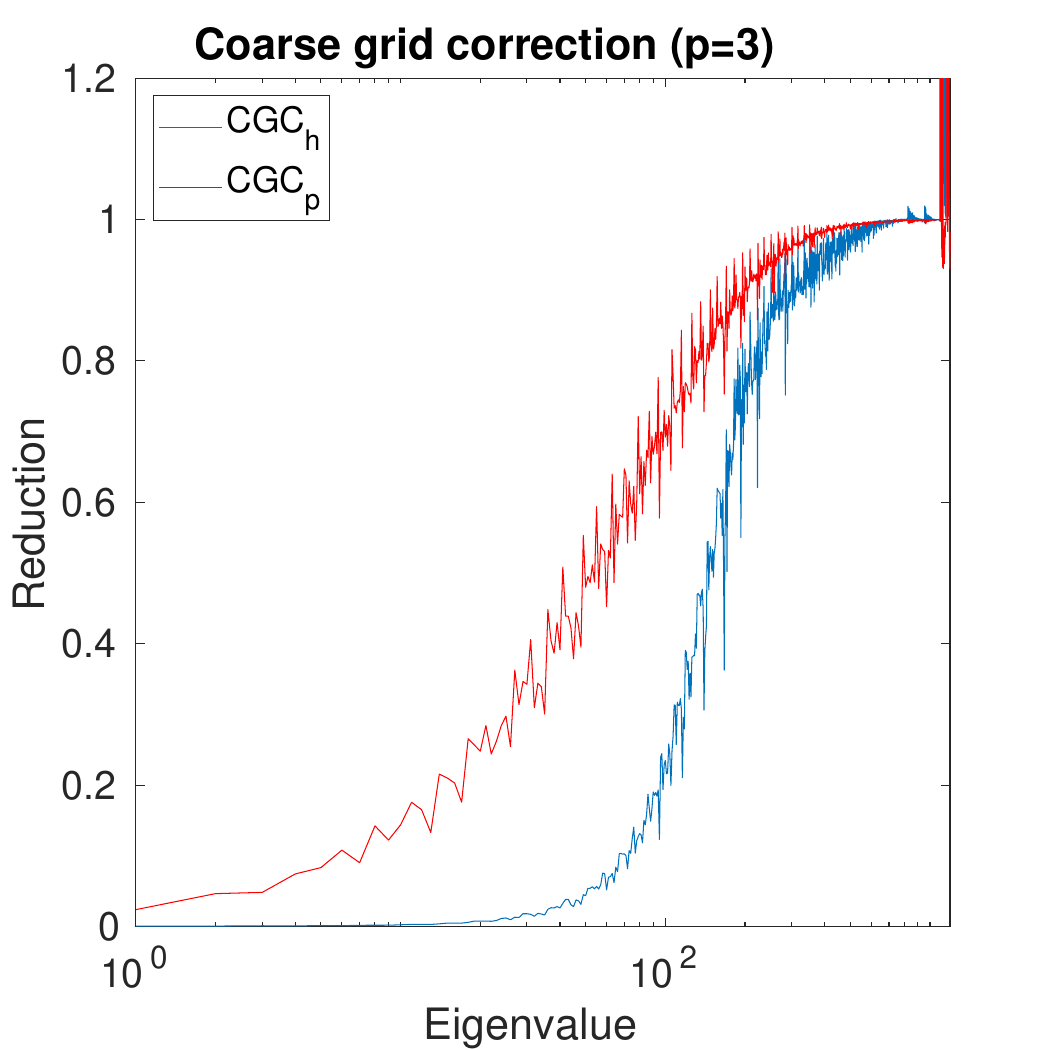}
\includegraphics[scale=0.31]{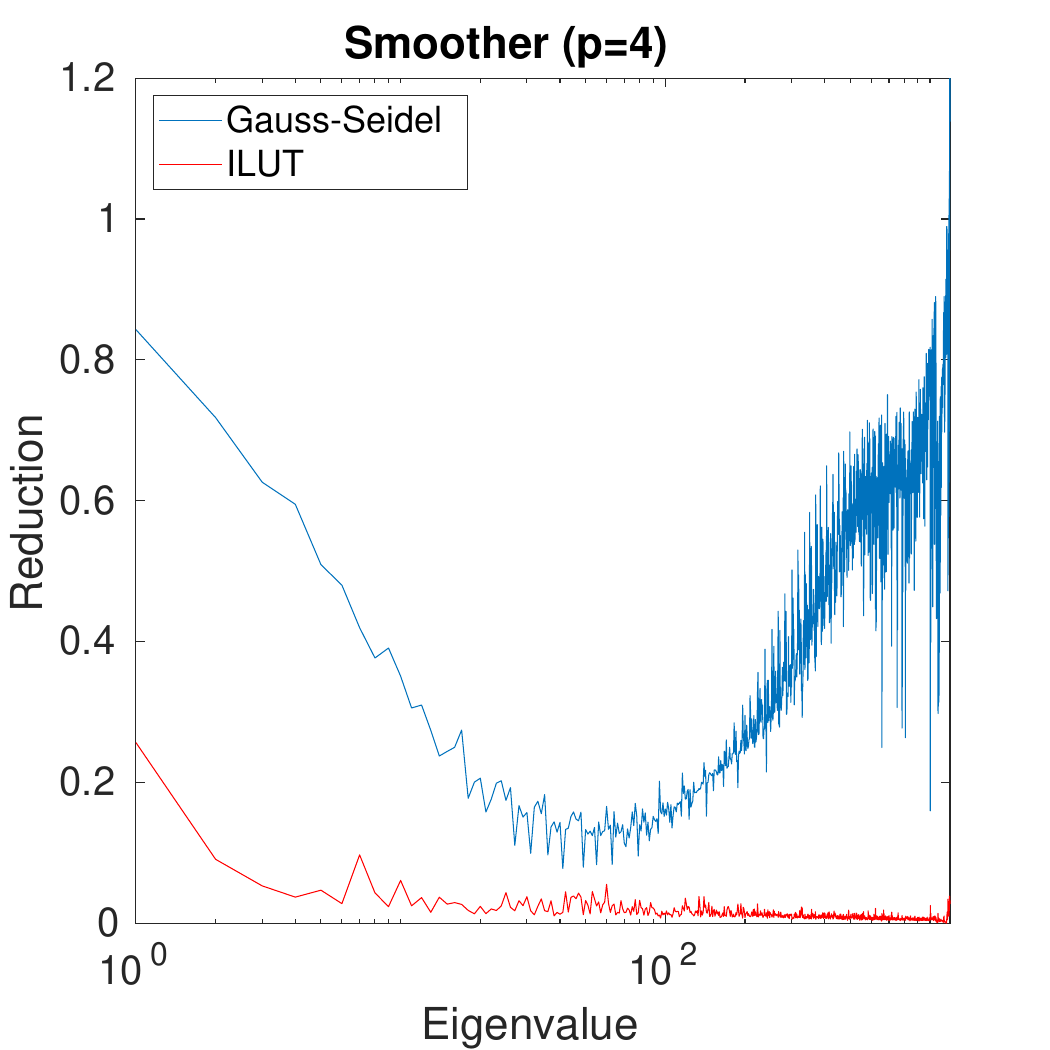}
\includegraphics[scale=0.31]{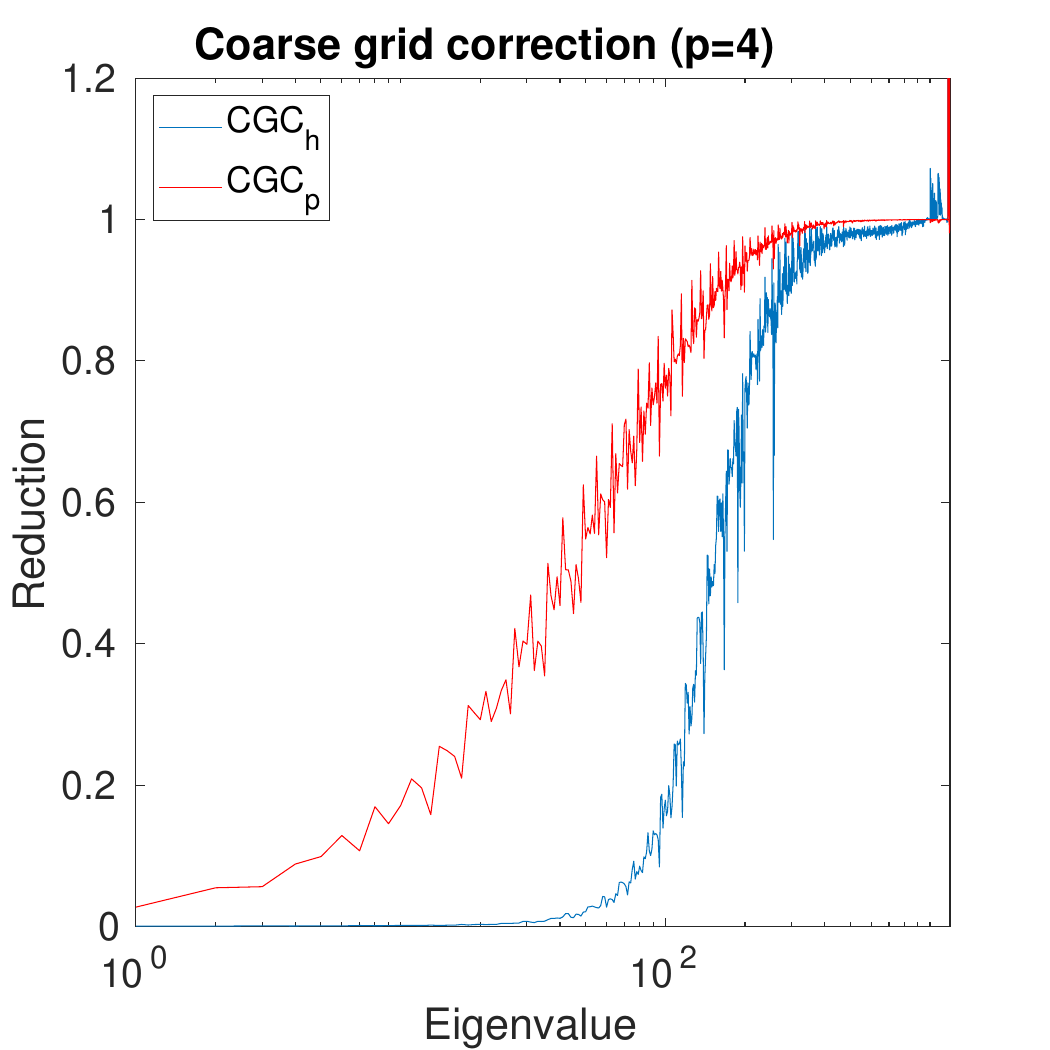}
\end{center} 
\caption{Error reduction in $(\mathbf{v}_j)$ for the second benchmark with $p=2,3,4$ and $h=2^{-5}$ obtained with different smoothers (left) and coarsening strategies (right).}
\label{fig:4}
\end{figure}

To investigate the effect of a single smoothing step or coarse grid correction, a spectral analysis \cite{spectra} is carried out for different values of $p$. For this analysis, we consider $ -\Delta u = 0$ with homogeneous Dirichlet boundary conditions and, hence, $u = 0$ as its exact solution. Let us define the error reduction factors as follows:
\vskip-.6cm
\begin{eqnarray}
 r^{\mathcal{S}}(\mathbf{u}_{h,p}^0) = \frac{|| \mathcal{S}_{h,p}(\mathbf{u}_{h,p}^0)||_2 }{||\mathbf{u}_{h,p}^0||_2}, \hspace{1.5cm} r^{\mathcal{CGC}}(\mathbf{u}_{h,p}^0) = \frac{|| \mathcal{CGC}(\mathbf{u}_{h,p}^0)||_2 }{||\mathbf{u}_{h,p}^0||_2},
 \end{eqnarray}

\noindent where $\mathcal{S}_{h,p}( \cdot)$ denotes a single smoothing step and $\mathcal{CGC}(\cdot)$ an exact coarse grid correction. We denote a coarse grid correction obtained by coarsening in $p$ and $h$ by $\mathcal{CGC}_p$ and $\mathcal{CGC}_h$, respectively. For $\mathcal{CGC}_p$, a direct projection to $p=1$ is considered. As an initial guess, the generalized eigenvectors $\mathbf{v}_i$ are chosen which satisfy
\vskip-.6cm
\begin{eqnarray}
\mathbf{A}_{h,p} \mathbf{v}_i = \lambda_i \mathbf{M}_{h,p} \mathbf{v}_i, \ i=1, \ldots ,N_{\rm dof}.
\end{eqnarray}
\noindent Here, $\mathbf{M}_{h,p}$ is the consistent mass matrix as defined in \eqref{PM}. The error reduction factors for the first benchmark for both smoothers and coarsening strategies are shown in Figure \ref{fig:3} for $h=2^{-5}$ and different values of $p$.  The reduction factors obtained with both smoothers are shown in the left column, while the plots in the right column show the reduction factors for both coarsening strategies. 
\noindent In general, the coarse grid corrections reduce the coefficients of the eigenvector expansion corresponding to the low-frequency modes, while the smoother reduces the coefficients associated with high-frequency modes. However, for increasing values of $p$, the reduction factors of the Gauss-Seidel smoother increase for the high-frequency modes, implying that the smoother becomes less effective. On the other hand, the use of ILUT as a smoother leads to decreasing reduction factors for all modes when the value of $p$ is increased. The coarse grid correction obtained by coarsening in $h$ (e.g. $\mathcal{CGC}_h$) is more effective compared to a correction obtained by coarsening $p$. Note that, for higher values of $p$, both types of coarse grid correction remain effective in reducing the coefficients of the eigenvector expansion corresponding to the low-frequency modes. Figure \ref{fig:4} shows the error reduction factors obtained for the second benchmark, showing similar, but less oscillatory, behaviour. These results indicate that the use of ILUT as a smoother (with $\nu_1 = \nu_2 = 1$) could significantly improve the convergence properties of the $p$-multigrid and $h$-multigrid method compared to the use of Gauss-Seidel as a smoother. \\

\subsection*{Iteration Matrix}
\label{im}

\begin{figure}[h!]
\begin{center}
\includegraphics[scale=0.31]{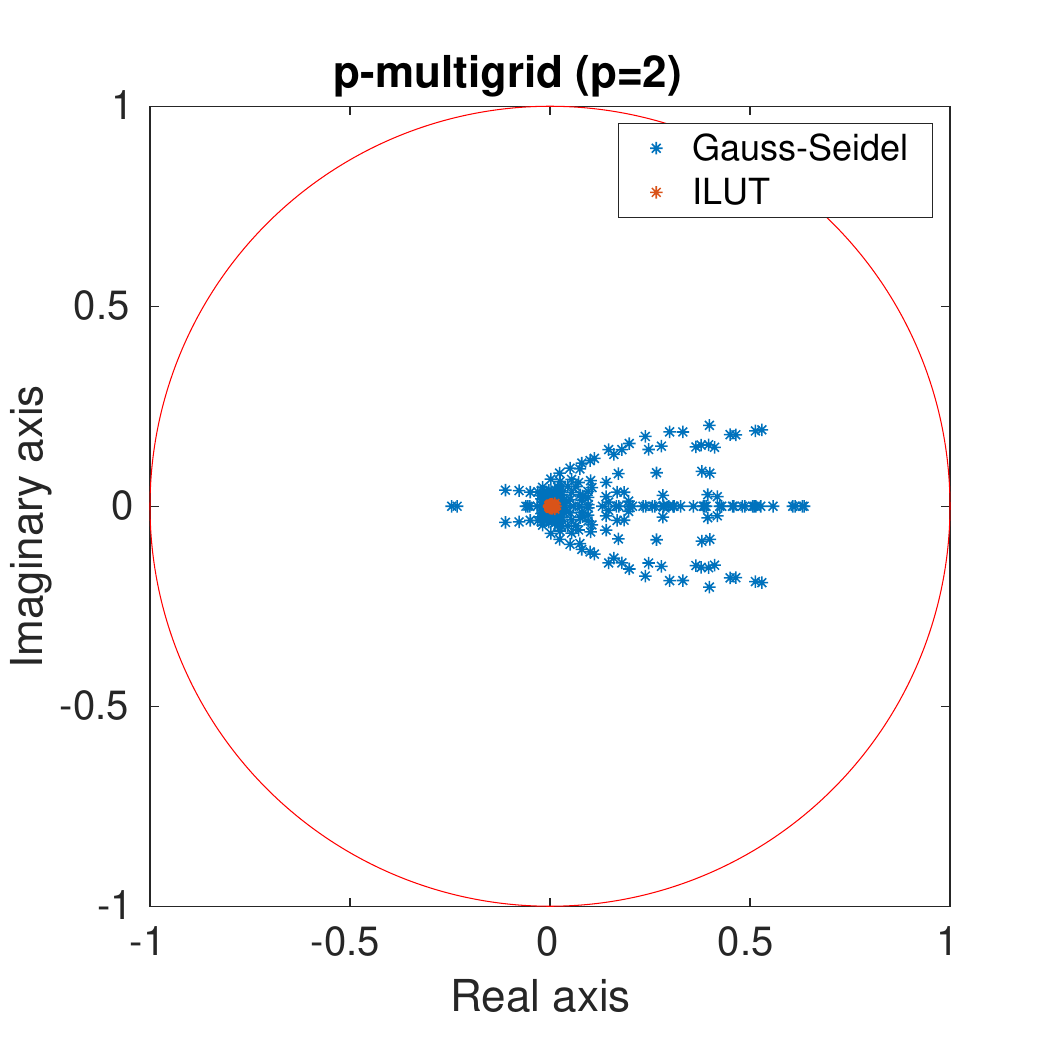}
\includegraphics[scale=0.31]{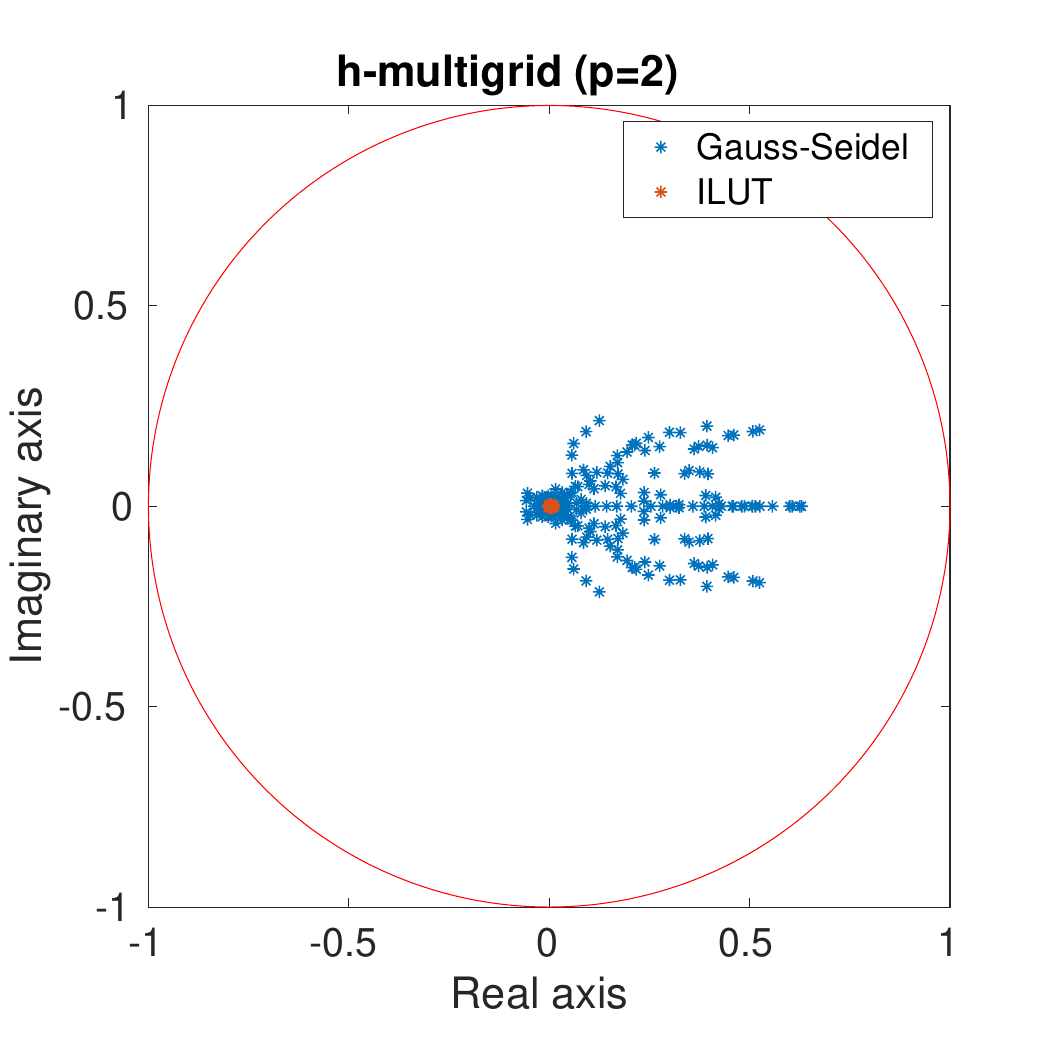}
\includegraphics[scale=0.31]{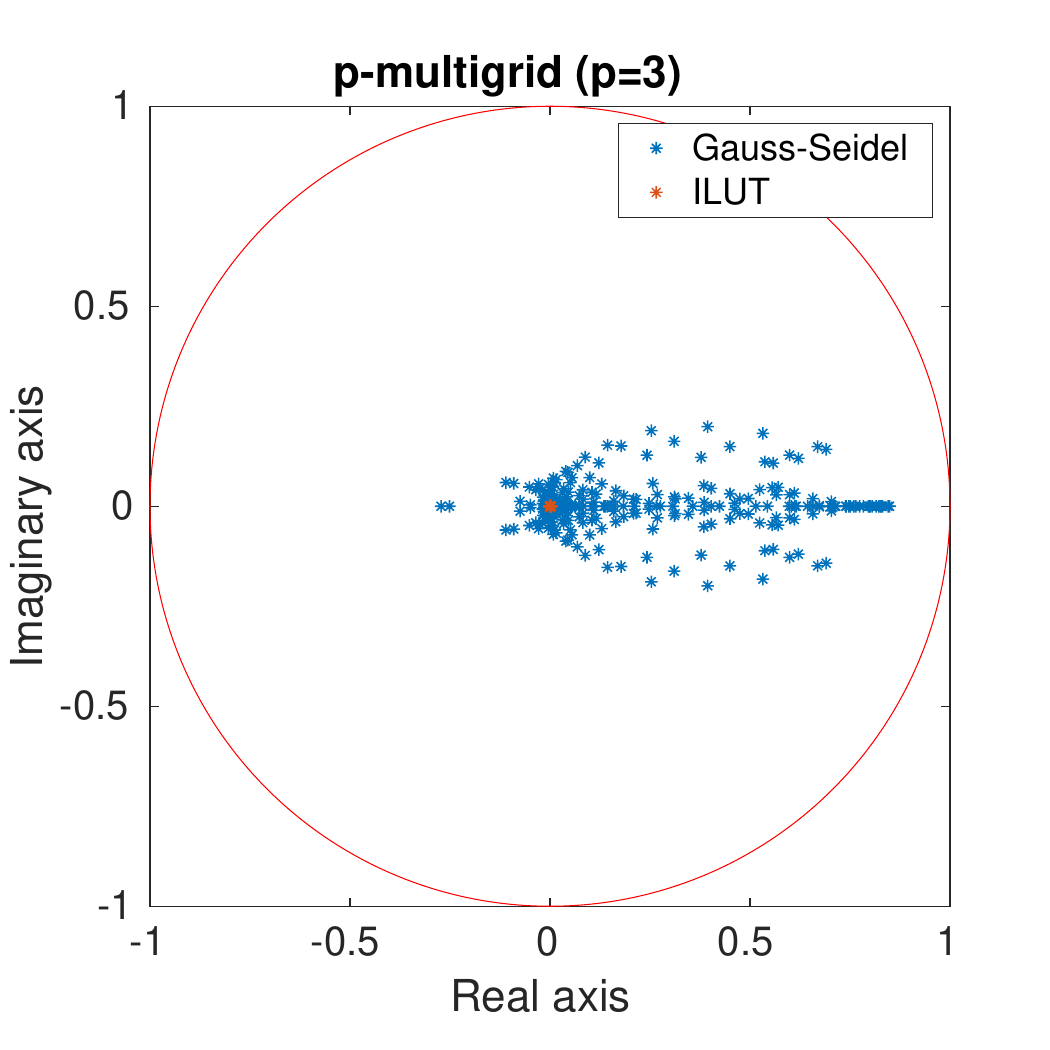}
\includegraphics[scale=0.31]{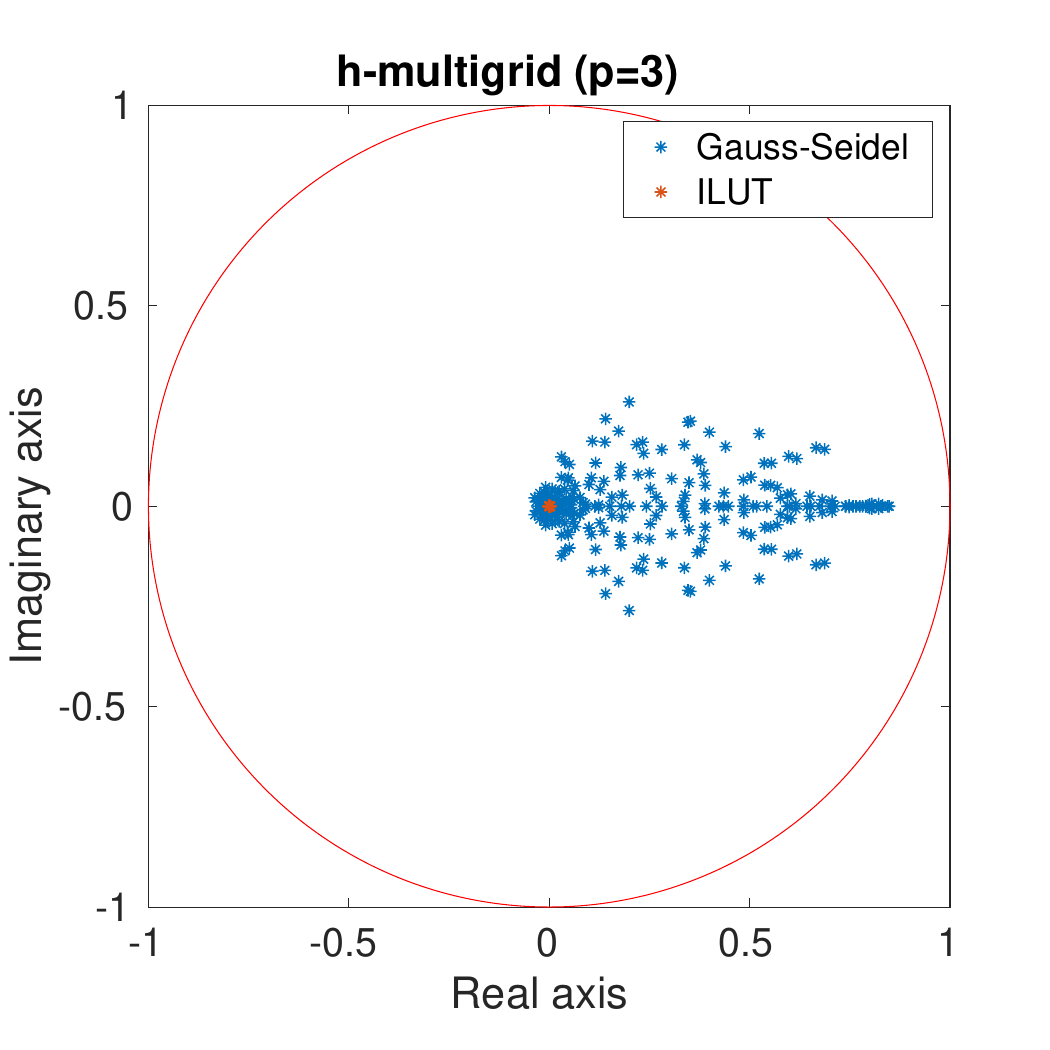}
\includegraphics[scale=0.31]{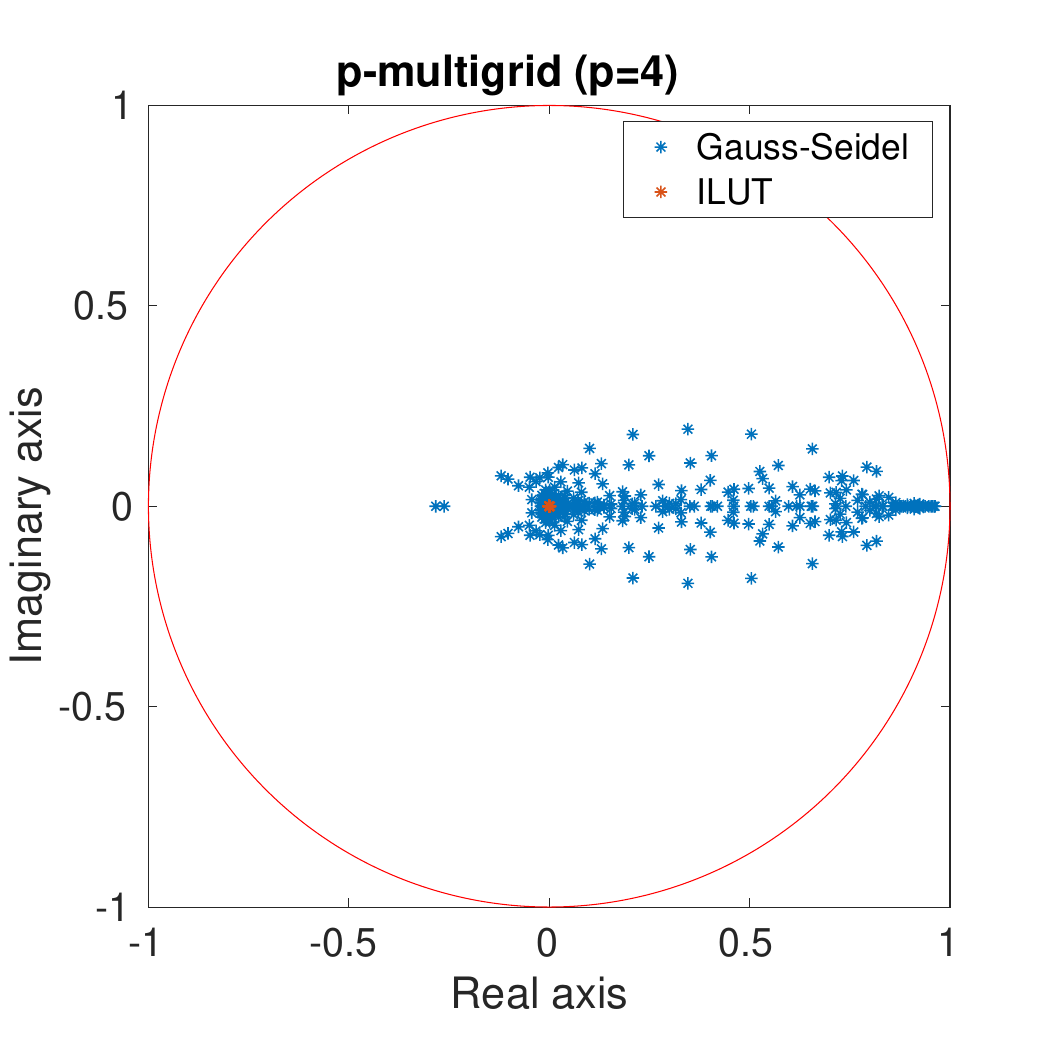}
\includegraphics[scale=0.31]{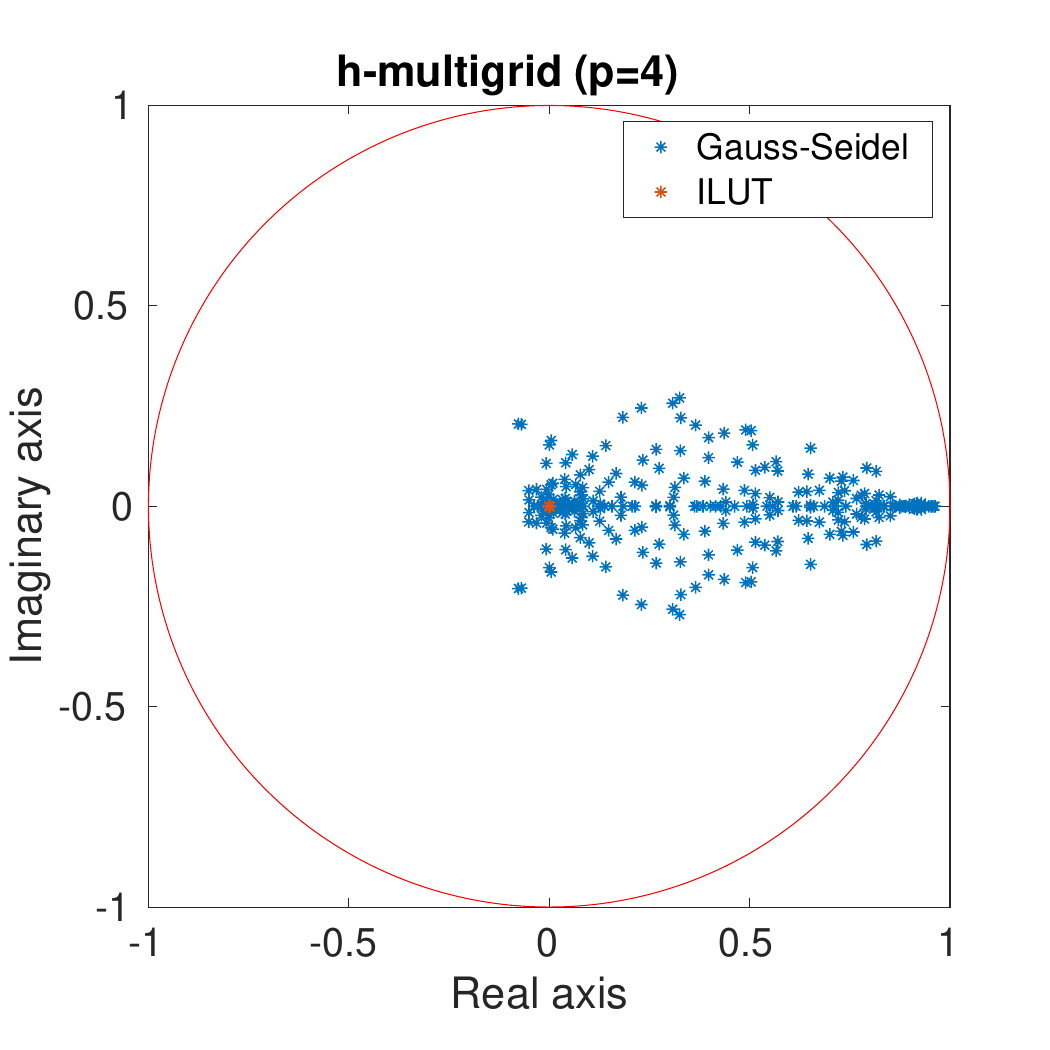}
\end{center}
\caption{Spectra of the iteration matrix for the first benchmark obtained with $p$-multigrid (left) and $h$-multigrid (right).}
\label{fig:5}
\end{figure}

For any multigrid method, the asymptotic convergence rate is determined by the spectral radius of the iteration matrix. To obtain this matrix explicitly, consider, again, $- \Delta u = 0$ with homogeneous Dirichlet boundary conditions. By applying a single iteration of the $p$-multigrid or $h$-multigrid method using the unit vector $\mathbf{e}_{h,p}^{i}$ as initial guess, one obtains the $i^{\rm th}$ column of the iteration matrix \cite{oosterlee}.  Figure \ref{fig:5} shows the spectra for the first benchmark for $h=2^{-5}$ and different values of $p$ obtained with both multigrid methods using Gauss-Seidel and ILUT as a smoother. For reference, the unit circle has been added in all plots. The spectral radius of the iteration matrix, defined as the maximum eigenvalue in absolute value, is then given by the eigenvalue that is the furthest away from the origin. Clearly, the spectral radius significantly increases for higher values of $p$ when adopting Gauss-Seidel as a smoother. The use of ILUT as a smoother results in spectra clustered around the origin, implying fast convergence of the resulting $p$-multigrid or $h$-multigrid method. 

\begin{figure}[h!]
\begin{center}
\includegraphics[scale=0.31]{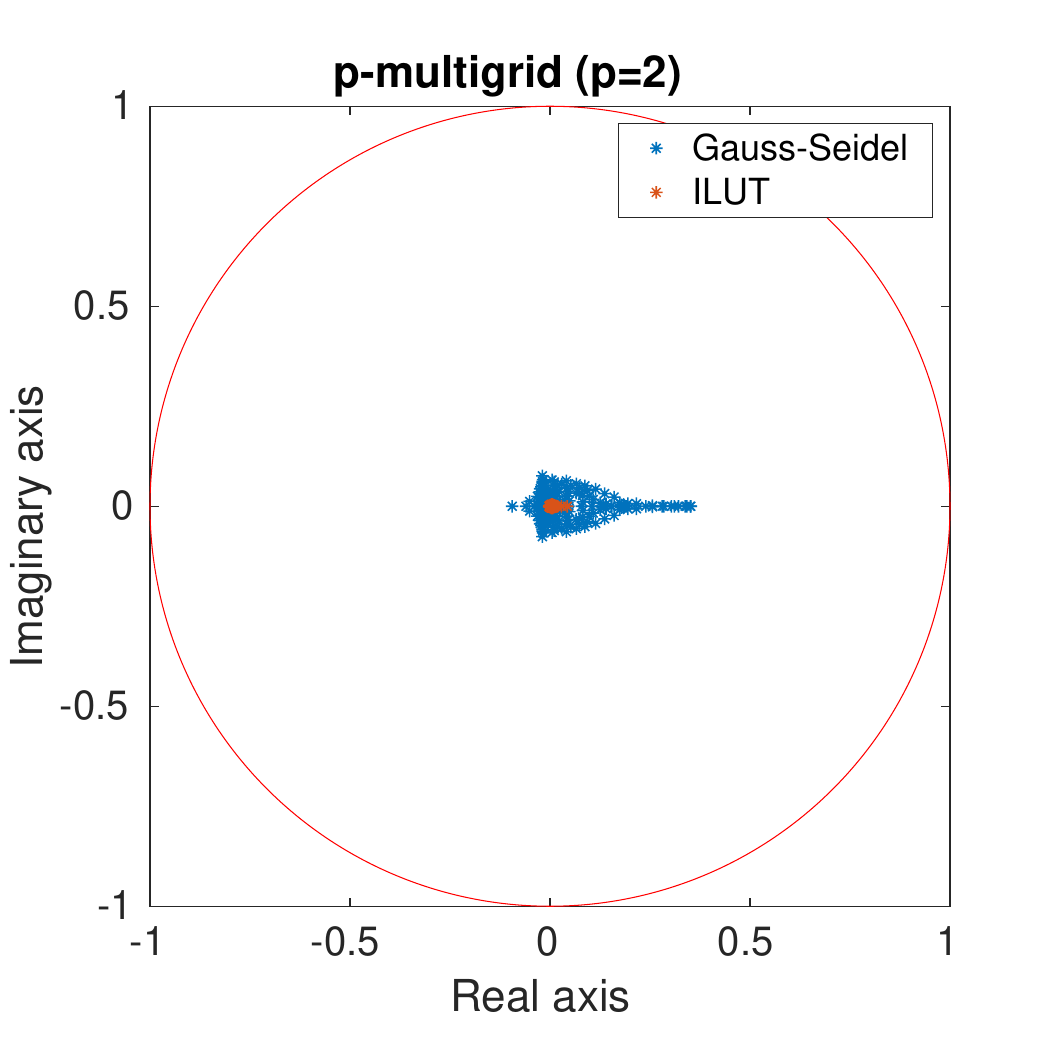}
\includegraphics[scale=0.31]{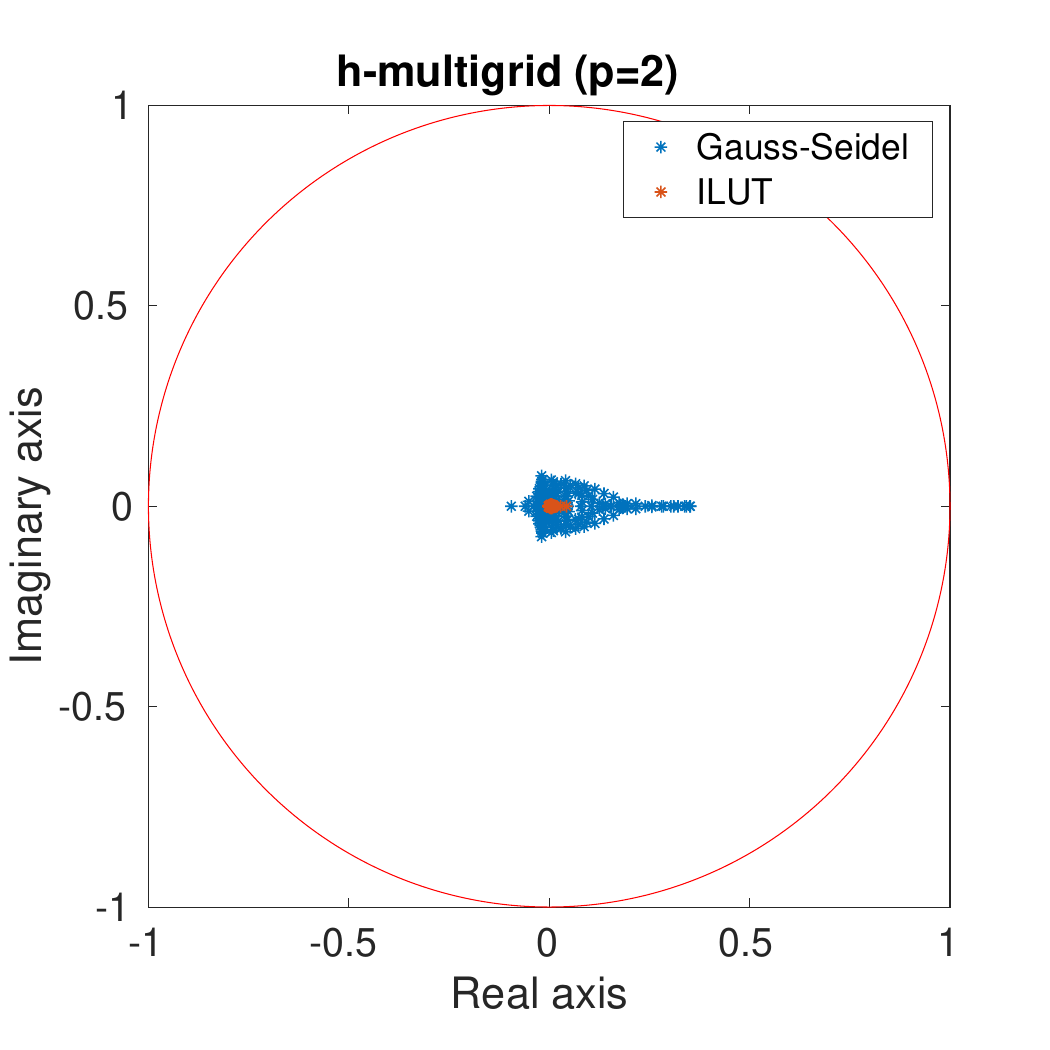}
\includegraphics[scale=0.31]{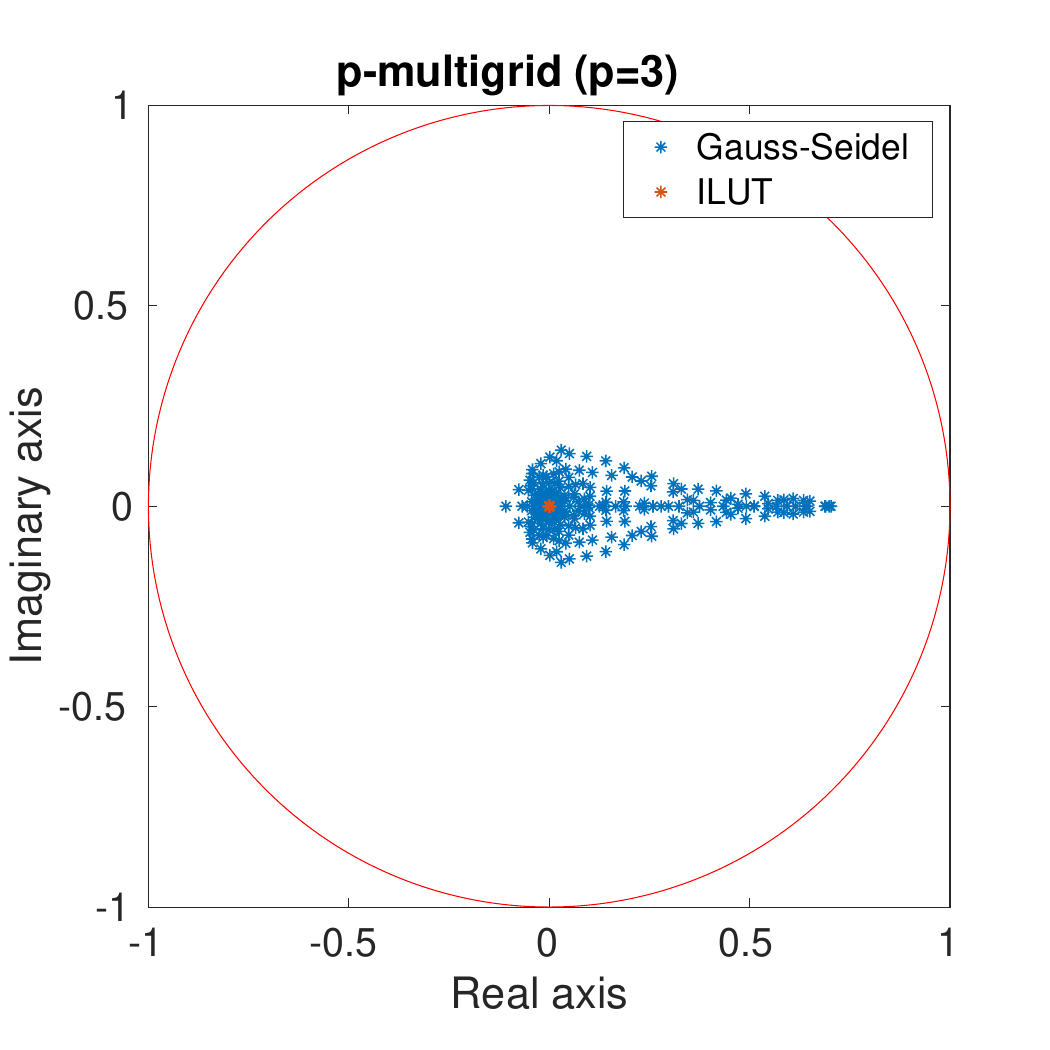}
\includegraphics[scale=0.31]{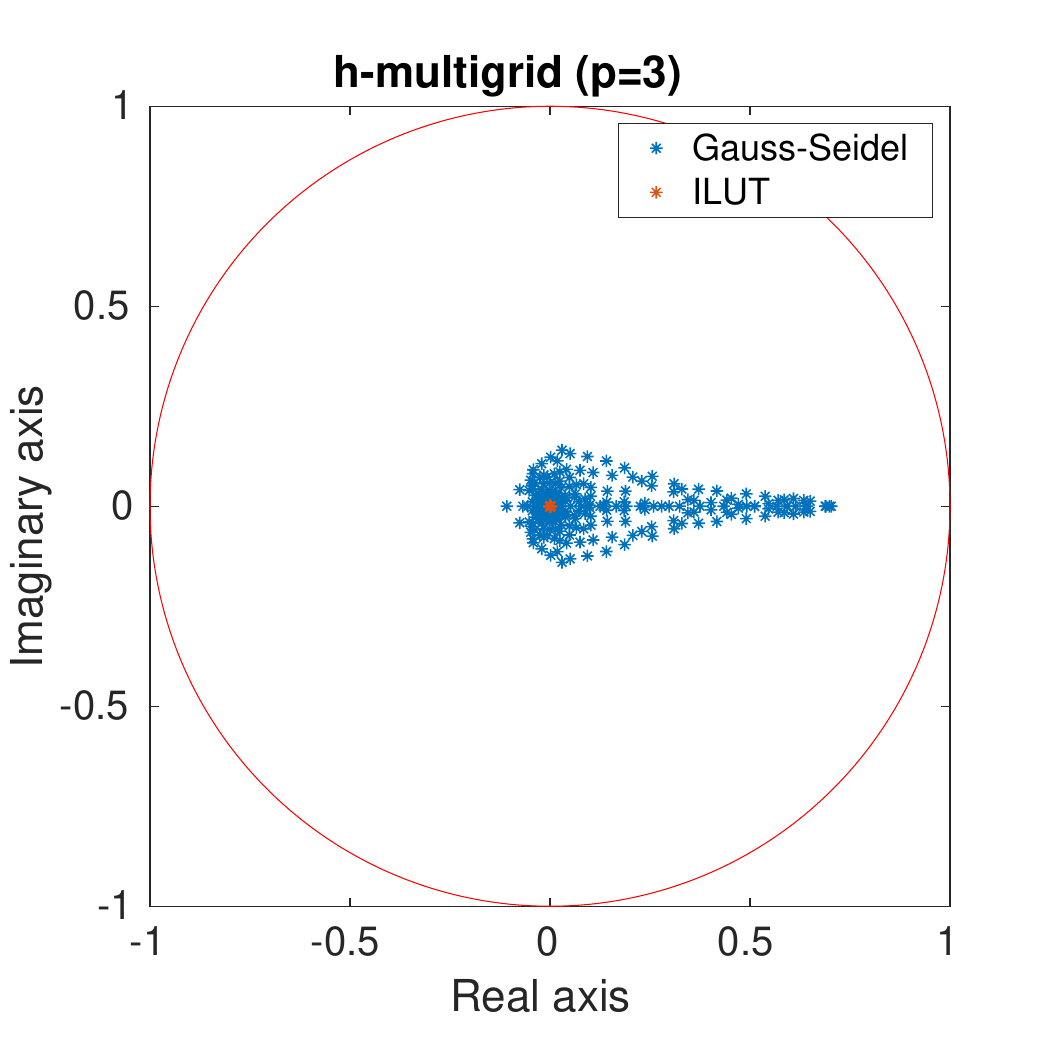}
\includegraphics[scale=0.31]{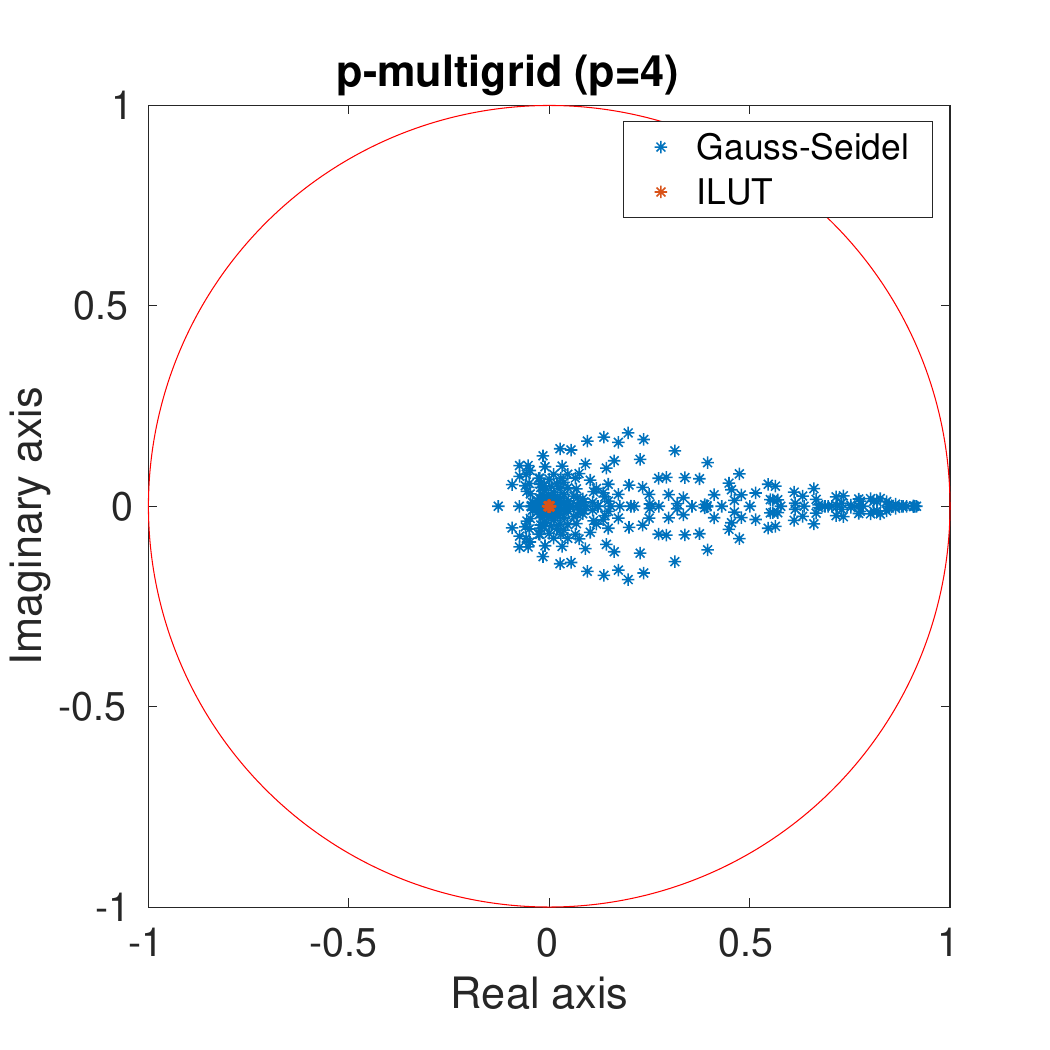}
\includegraphics[scale=0.31]{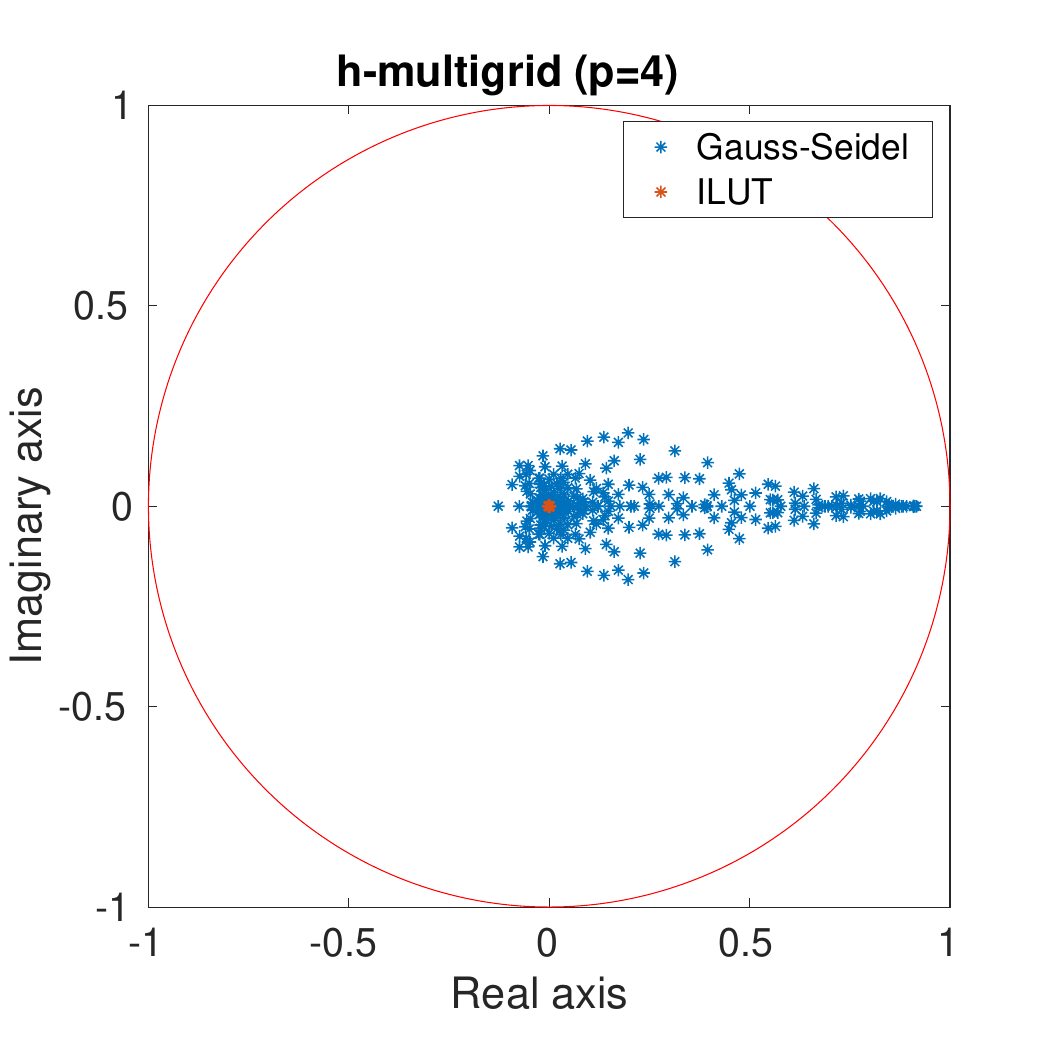}
\end{center}
\caption{Spectra of the iteration matrix for the second benchmark obtained with $p$-multigrid (left) and $h$-multigrid (right).}
\label{fig:6}
\end{figure}

\noindent The spectra of the iteration matrices for the second benchmark are presented in Figure \ref{fig:6}. Although the eigenvalues are more clustered with Gauss-Seidel compared to the first benchmark, the same behaviour can be observed.  

\noindent The spectral radii for both benchmarks, where $\nu_1 = \nu_2=1$, are presented in Table \ref{tab:2a} until Table \ref{tab:3b}. For Gauss-Seidel, the spectral radius of the iteration matrix is independent of the mesh width $h$ and coarsening strategy, but depends strongly on the approximation order $p$.  
\noindent The use of ILUT leads to a spectral radius which is significantly lower for all values of $h$ and $p$. Although ILUT exhibits a small $h$-dependence, the spectral radius remains low for all values of $h$ and both coarsening strategies. As a consequence, the $p$-multigrid and $h$-multigrid method are expected to show both $h$- and $p$-independence convergence behaviour when ILUT is adopted as a smoother.

\begin{table}[h!]
\centering
\caption{Spectral radius $\rho$ for the first benchmark for $p$-multigrid using Gauss-Seidel and ILUT for different values of $h$ and $p$.}
\begin{tabular}{l|ll|l|ll|l|ll}

$p=2$        &  GS     &   ILUT   &     $p=3$    &  GS    &  ILUT  & $p=4$       &   GS    & ILUT      \\ \hline
$h = 2^{-4}$ & $0.635$ &  $0.014$ &$h = 2^{-4}$  & $0.849$& $0.003$& $h = 2^{-4}$& $0.963$ & $0.003$   \\
$h = 2^{-5}$ & $0.631$ &  $0.039$ &$h = 2^{-5}$  & $0.845$& $0.019$& $h = 2^{-5}$& $0.960$ & $0.029$   \\
$h = 2^{-6}$ & $0.647$ &  $0.058$ &$h = 2^{-6}$  & $0.844$& $0.017$& $h = 2^{-6}$& $0.960$ & $0.023$   \\
\end{tabular}
\label{tab:2a}
\end{table}

\begin{table}[h!]
\centering
\caption{Spectral radius $\rho$ for the first benchmark for $h$-multigrid using Gauss-Seidel and ILUT for different values of $h$ and $p$.}
\begin{tabular}{l|ll|l|ll|l|ll}
$p=2$        &  GS     &   ILUT   &     $p=3$    &  GS    &  ILUT  & $p=4$       &   GS    & ILUT      \\ \hline
$h = 2^{-4}$ & $0.630$ &  $0.012$ &$h = 2^{-4}$  & $0.848$& $0.004$& $h = 2^{-4}$& $0.963$ & $0.003$   \\
$h = 2^{-5}$ & $0.627$ &  $0.039$ &$h = 2^{-5}$  & $0.845$& $0.018$& $h = 2^{-5}$& $0.960$ & $0.029$   \\
$h = 2^{-6}$ & $0.646$ &  $0.059$ &$h = 2^{-6}$  & $0.844$& $0.014$& $h = 2^{-6}$& $0.960$ & $0.023$   \\
\end{tabular}
\label{tab:2b}
\end{table}

\begin{table}[h!]
\centering
\caption{Spectral radius $\rho$ for the second benchmark for $p$-multigrid using Gauss-Seidel and ILUT for different values of $h$ and $p$.}
\begin{tabular}{l|ll|l|ll|l|ll}
$p=2$    	 &  GS      &   ILUT   & $p=3$       &   GS    &  ILUT  & $p=4$       &   GS    & ILUT      \\ \hline
$h = 2^{-4}$ & $0.352$  &  $0.043$ & $h = 2^{-4}$& $0.703$ & $0.002$& $h = 2^{-4}$&	$0.916$ & $0.003$	\\
$h = 2^{-5}$ & $0.351$  &  $0.037$ & $h = 2^{-5}$& $0.699$ & $0.011$& $h = 2^{-5}$&	$0.913$ & $0.020$   \\
$h = 2^{-6}$ & $0.352$  &  $0.042$ & $h = 2^{-6}$& $0.699$ & $0.017$& $h = 2^{-6}$& $0.914$	& $0.016$   \\
\end{tabular}
\label{tab:3a}
\end{table}

\begin{table}[h!]
\centering
\caption{Spectral radius $\rho$ for the second benchmark for $h$-multigrid using Gauss-Seidel and ILUT for different values of $h$ and $p$.}
\begin{tabular}{l|ll|l|ll|l|ll}

$p=2$    	 &  GS      &   ILUT   & $p=3$       &   GS    &  ILUT  & $p=4$       &   GS    & ILUT      \\ \hline
$h = 2^{-4}$ & $0.367$  &  $0.043$ & $h = 2^{-4}$& $0.698$ & $0.002$& $h = 2^{-4}$&	$0.916$ & $0.003$	\\
$h = 2^{-5}$ & $0.367$  &  $0.036$ & $h = 2^{-5}$& $0.696$ & $0.008$& $h = 2^{-5}$&	$0.913$ & $0.020$   \\
$h = 2^{-6}$ & $0.359$  &  $0.042$ & $h = 2^{-6}$& $0.698$ & $0.006$& $h = 2^{-6}$& $0.913$	& $0.016$   \\
\end{tabular}
\label{tab:3b}
\end{table}

\newpage

\section{Numerical Results}
\label{numex}

In the previous Section, a spectral analysis showed that the use of ILUT as a smoother within a $p$-multigrid or $h$-multigrid method significantly improves the asymptotic convergence rate compared to the use of Gauss-Seidel as a smoother. In this Section, $p$-multigrid and $h$-multigrid are both applied as a stand-alone solver and as a preconditioner within a stabilized BiConjugate Gradient (BiCGSTAB) to verify this analysis. Results in terms of iteration numbers and CPU times are obtained using Gauss-Seidel and ILUT as a smoother. Furthermore, the proposed $p$-multigrid method is compared to an $h$-multigrid method using a non-standard smoother. Finally, the $p$-multigrid method is adopted for discretizations using THB-splines.

\noindent For all numerical experiments, the initial guess $\mathbf{u}_{h,p}^{(0)}$ is chosen randomly, where each entry is sampled from a uniform distribution on the interval $[-1,1]$ using the same seed. Furthermore, we choose $\nu_1 = \nu_2 = 1$ for consistency. Application of multiple smoothing steps, which is in particular common for Gauss-Seidel, decreases the number of iterations until convergence, but does not qualitatively or quantitatively change the $p$-dependence. Furthermore, since the smoother costs dominate when solving the linear systems, CPU times are only mildly affected. The stopping criterion is based on a relative reduction of the initial residual, where a tolerance of  $\epsilon = 10^{-8}$ is adopted. Boundary conditions are imposed by eliminating the degrees of freedom associated to the boundary. Since the use of a V-cycle or W-cycle led to the same number of iterations and the computational costs per cycle is lower for V-cycles, they are considered throughout the remainder of this paper.  
 
\subsection*{$p$-Multigrid as stand-alone solver}
\noindent Table \ref{tab:4} shows the number of V-cycles needed to achieve convergence using different smoothers for all benchmarks. For the first three benchmarks, the number of V-cycles needed with Gauss-Seidel is in general independent of the mesh width $h$, but strongly depends on the approximation order $p$. For some configurations, however, the use of Gauss-Seidel leads to a method that diverges, indicated with $(-)$. The $p$-multigrid method was said to be diverged in case the relative residual exceeded $10^{10}$ at the end of a V-cycle. In general, the use of ILUT as a smoother leads to a $p$-multigrid which converges for all configurations and exhibits both independence of $h$ and $p$. Only for the third benchmark, a logarithmic dependence in $h$ is visible for $p=2$. Furthermore, the number of iterations needed for convergence is significantly lower. 

\noindent For Poisson's equation on the unit cube, the $p$-dependence when Gauss-Seidel is adopted as smoother is stronger compared to the dependence for the twodimensional benchmarks. Furthermore, the number of iterations slightly decreases for smaller values of the meshwidth $h$. The number of iterations needed with ILUT as a smoother, is effectively independent of the approximation order $p$.

\begin{table}[h!]
\begin{subtable}{1\textwidth}
\centering
\begin{tabular}{c|cc|cc|cc|cc}
     & \multicolumn{2}{c|}{$p=2$}  & \multicolumn{2}{c|}{$p=3$}  & \multicolumn{2}{c|}{$p=4$} & \multicolumn{2}{c}{$p=5$}  \\ 
           & ILUT  & GS  & ILUT &GS   &  ILUT & GS   & ILUT & GS  \\ \hline 
$h=2^{-6}$ & $4$   &$30$ & $3$	&$62$ & $3$	  &$176$	 &  $3$ &  $491$    \\
$h=2^{-7}$ & $4$   &$29$ & $3$	&$61$ & $3$	  &$172$	 &  $3$ &  $499$    \\
$h=2^{-8}$ & $5$   &$30$ & $3$	&$60$ & $3$	  &$163$     &	$3$ &  $473$    \\
$h=2^{-9}$ & $5$   &$32$ & $3$	&$61$ & $3$	  &$163$	 &  $3$ &  $452$    \\
\end{tabular}
\end{subtable}
\begin{subtable}{1\textwidth}
\centering
\subcaption{Poisson's equation on quarter annulus}
\begin{tabular}{c|cc|cc|cc|cc}
     & \multicolumn{2}{c|}{$p=2$}  & \multicolumn{2}{c|}{$p=3$}  & \multicolumn{2}{c|}{$p=4$} & \multicolumn{2}{c}{$p=5$}  \\ 
           & ILUT  & GS  & ILUT &GS   &  ILUT & GS   & ILUT & GS  \\ \hline 
$h=2^{-6}$ & $5$   &$-$  & $3$	&$-$  & $3$	  &$-$	 &  $4$ &  $-$  \\
$h=2^{-7}$ & $5$   &$-$  & $3$	&$-$  & $4$	  &$-$	 &  $4$ &  $-$  \\
$h=2^{-8}$ & $5$   &$-$  & $3$	&$-$  & $3$	  &$-$	 &  $4$ &  $-$  \\
$h=2^{-9}$ & $5$   &$-$  & $4$	&$-$  & $3$	  &$-$	 &  $4$ &  $-$  \\
\end{tabular}
\end{subtable}
\begin{subtable}{1\textwidth}
\centering
\subcaption{CDR-equation on unit square}
\begin{tabular}{c|cc|cc|cc|cc}
     & \multicolumn{2}{c|}{$p=2$}  & \multicolumn{2}{c|}{$p=3$}  & \multicolumn{2}{c|}{$p=4$} & \multicolumn{2}{c}{$p=5$}  \\ 
           & ILUT  & GS  & ILUT &GS   &  ILUT & GS   & ILUT & GS  \\ \hline 
$h=2^{-6}$ & $6$   &$24$  & $6$	&$52$  & $5$  &$115$ &  $5$ &  $333$     \\
$h=2^{-7}$ & $7$   &$24$  & $6$	&$53$  & $5$  &$133$ &  $5$ &  $325$      \\
$h=2^{-8}$ & $8$   &$26$  & $6$	&$54$  & $6$  &$126$ &	$6$ &  $322$    \\
$h=2^{-9}$ & $9$   &$26$  & $6$ &$55$  & $6$  &$131$ &  $5$ &  $327$    \\
\end{tabular}
\subcaption{Poisson's equation on L-shaped domain}
\end{subtable}
\begin{subtable}{1\textwidth}
\centering
\begin{tabular}{c|cc|cc|cc|cc}
     & \multicolumn{2}{c|}{$p=2$}  & \multicolumn{2}{c|}{$p=3$}  & \multicolumn{2}{c|}{$p=4$} & \multicolumn{2}{c}{$p=5$}  \\ 
           & ILUT  & GS  & ILUT &GS   &  ILUT & GS   & ILUT & GS  \\ \hline 
$h=2^{-2}$ & $3$   &$65$  & $3$	&$405$  & $3$  &$3255$ &  $5$ &  $22787$     \\
$h=2^{-3}$ & $3$   &$59$  & $3$	&$339$  & $3$  &$2063$ &  $3$ &  $8128$      \\
$h=2^{-4}$ & $3$   &$57$  & $3$	&$281$  & $3$  &$1652$ &  $5$ &  $7152$    \\
$h=2^{-5}$ & $4$   &$55$  & $3$ &$273$  & $5$  &$1361$ &  $10$&  $6250$    \\
\end{tabular}
\subcaption{Poisson's equation on the unit cube}
\end{subtable}
\caption{Number of V-cycles needed to achieve convergence with $p$-multigrid.}
\label{tab:4} 
\end{table} 

\newpage

\subsection*{$h$-Multigrid as stand-alone solver}
\noindent Table \ref{tab:4h} shows the number of V-cycles needed to achieve convergence using a $h$-multigrid method. As expected from the spectral analysis, the number of V-cycles needed with Gauss-Seidel is in general independent of the mesh width $h$, but strongly depends on the approximation order $p$. The use of ILUT as a smoother leads to a $h$-multigrid which converges for all configurations and exhibits both independence of $h$ and $p$. Furthermore, the number of iterations needed for convergence is significantly lower. Compared to the use of $p$-multigrid as a method, the results are very similar. For benchmark $3$, however, the number of iterations needed with $h$-multigrid using ILUT as a smoother is slightly lower compared to the $p$-multigrid method. For some configurations, the $h$-multigrid does not converge when applied to the threedimensional benchmark, which is denoted by $(-)$.

\begin{table}[h!]
\begin{subtable}{1\textwidth}
\centering
\begin{tabular}{c|cc|cc|cc|cc}
     & \multicolumn{2}{c|}{$p=2$}  & \multicolumn{2}{c|}{$p=3$}  & \multicolumn{2}{c|}{$p=4$} & \multicolumn{2}{c}{$p=5$}  \\ 
           & ILUT  & GS  & ILUT &GS   &  ILUT & GS   & ILUT & GS  \\ \hline 
$h=2^{-6}$ & $4$   &$30$ & $3$	&$62$ & $3$	  &$175$	 &  $3$ &  $492$    \\
$h=2^{-7}$ & $4$   &$29$ & $3$	&$61$ & $3$	  &$172$	 &  $3$ &  $499$    \\
$h=2^{-8}$ & $5$   &$30$ & $3$	&$60$ & $3$	  &$163$     &	$3$ &  $473$    \\
$h=2^{-9}$ & $5$   &$32$ & $3$	&$61$ & $3$	  &$164$	 &  $3$ &  $452$    \\
\end{tabular}
\end{subtable}
\begin{subtable}{1\textwidth}
\centering
\subcaption{Poisson's equation on quarter annulus}
\begin{tabular}{c|cc|cc|cc|cc}
     & \multicolumn{2}{c|}{$p=2$}  & \multicolumn{2}{c|}{$p=3$}  & \multicolumn{2}{c|}{$p=4$} & \multicolumn{2}{c}{$p=5$}  \\ 
           & ILUT  & GS  & ILUT &GS   &  ILUT & GS   & ILUT & GS  \\ \hline 
$h=2^{-6}$ & $5$   &$-$  & $3$	&$-$  & $3$	  &$-$	 &  $4$ &  $-$  \\
$h=2^{-7}$ & $5$   &$-$  & $3$	&$-$  & $4$	  &$-$	 &  $4$ &  $-$  \\
$h=2^{-8}$ & $5$   &$-$  & $3$	&$-$  & $3$	  &$-$	 &  $4$ &  $-$  \\
$h=2^{-9}$ & $5$   &$-$  & $3$	&$-$  & $3$	  &$-$	 &  $4$ &  $-$  \\
\end{tabular}
\end{subtable}
\begin{subtable}{1\textwidth}
\centering
\subcaption{CDR-equation on unit square}
\begin{tabular}{c|cc|cc|cc|cc}
     & \multicolumn{2}{c|}{$p=2$}  & \multicolumn{2}{c|}{$p=3$}  & \multicolumn{2}{c|}{$p=4$} & \multicolumn{2}{c}{$p=5$}  \\ 
           & ILUT  & GS  & ILUT &GS   &  ILUT & GS   & ILUT & GS  \\ \hline 
$h=2^{-6}$ & $6$   &$25$  & $3$	&$53$  & $3$  &$113$ &  $3$ &  $332$     \\
$h=2^{-7}$ & $7$   &$25$  & $3$	&$53$  & $3$  &$133$ &  $3$ &  $325$      \\
$h=2^{-8}$ & $8$   &$26$  & $3$	&$54$  & $3$  &$127$ &	$3$ &  $322$    \\
$h=2^{-9}$ & $9$   &$26$  & $3$ &$55$  & $3$  &$131$ &  $3$ &  $327$    \\
\end{tabular}
\subcaption{Poisson's equation on L-shaped domain}
\end{subtable}
\end{table}
\begin{table}
\ContinuedFloat
\begin{subtable}{1\textwidth}
\centering
\begin{tabular}{c|cc|cc|cc|cc}
     & \multicolumn{2}{c|}{$p=2$}  & \multicolumn{2}{c|}{$p=3$}  & \multicolumn{2}{c|}{$p=4$} & \multicolumn{2}{c}{$p=5$}  \\ 
           & ILUT  & GS  & ILUT &GS   &  ILUT & GS   & ILUT & GS  \\ \hline 
$h=2^{-2}$ & $3$   &$65$  & $3$	&$405$  & $3$  &$2269$ &  $5$ &  $-$     \\
$h=2^{-3}$ & $3$   &$59$  & $3$	&$339$  & $3$  &$-$    &  $3$ &  $-$      \\
$h=2^{-4}$ & $3$   &$57$  & $3$	&$281$  & $3$  &$1652$ &  $5$ &  $-$    \\
$h=2^{-5}$ & $4$   &$55$  & $3$ &$273$  & $5$  &$1361$ &  $10$&  $6248$    \\
\end{tabular}
\subcaption{Poisson's equation on the unit cube}
\end{subtable}
\caption{Number of V-cycles needed to achieve convergence with $h$-multigrid.}
\label{tab:4h} 
\end{table} 

\subsection*{$p$-Multigrid as a preconditioner}
\noindent As an alternative, the $p$-multigrid method can be applied as a preconditioner within a BiCGSTAB method. In the preconditioning phase of each iteration, a single V-cycle is applied. The number of iterations needed to achieve convergence can be found in Table \ref{tab:5}. When applying Gauss-Seidel as a smoother, the number of iterations needed with BiCGSTAB is significantly lower compared to the number of $p$-multigrid V-cycles and even restores stability for higher values of $p$ (see Table \ref{tab:4}). However, a dependence of the iteration numbers on $p$ is still present. When adopting ILUT as a smoother, the number of iterations needed for convergence slightly decreases compared to the number of $p$-multigrid V-cycles for all configurations and benchmarks. Furthermore, the number of iterations is independent of both $h$ and $p$.

\begin{table}[h!]
\begin{subtable}{1\textwidth}
\centering
\begin{tabular}{c|cc|cc|cc|cc}
     & \multicolumn{2}{c|}{$p=2$}  & \multicolumn{2}{c|}{$p=3$}  & \multicolumn{2}{c|}{$p=4$} & \multicolumn{2}{c}{$p=5$}  \\ 
           & ILUT  & GS  & ILUT   &GS   &  ILUT & GS  & ILUT & GS     \\ \hline 
$h=2^{-6}$ & $2$   &$13$ &  $2$	  &$18$ & 	$2$	&$41$ &	$2$  & $78$\\
$h=2^{-7}$ & $2$   &$12$ &  $2$	  &$20$ & 	$2$	&$41$ &	$2$  & $92$\\
$h=2^{-8}$ & $3$   &$13$ &  $2$	  &$19$ &   $2$	&$43$ &	$2$  & $95$    \\
$h=2^{-9}$ & $3$   &$13$ &  $2$	  &$21$ &   $2$	&$41$ & $2$  & $95$    \\
\end{tabular}
\subcaption{Poisson's equation on quarter annulus}
\end{subtable}
\begin{subtable}{1\textwidth}
\centering
\begin{tabular}{c|cc|cc|cc|cc}
     & \multicolumn{2}{c|}{$p=2$}  & \multicolumn{2}{c|}{$p=3$}  & \multicolumn{2}{c|}{$p=4$} & \multicolumn{2}{c}{$p=5$}  \\ 
           & ILUT  & GS  & ILUT   &GS   &  ILUT & GS  & ILUT & GS     \\ \hline 
$h=2^{-6}$ & $2$   &$7$ & $2$	  &$13$ &	$2$	&$29$ & $2$ &	$65$ \\
$h=2^{-7}$ & $2$   &$8$ & $2$	  &$13$ &	$2$	&$29$ & $2$ &	$70$\\
$h=2^{-8}$ & $2$   &$7$ & $2$	  &$12$ &   $2$	&$29$ & $2$ &   $64$  \\
$h=2^{-9}$ & $2$   &$7$ & $2$	  &$14$ &   $2$	&$28$ & $2$ &   $72$  \\
\end{tabular}
\subcaption{CDR-equation on unit square}
\end{subtable}
\end{table}
\begin{table}
\ContinuedFloat
\begin{subtable}{1\textwidth}
\centering
\begin{tabular}{c|cc|cc|cc|cc}
     & \multicolumn{2}{c|}{$p=2$}  & \multicolumn{2}{c|}{$p=3$}  & \multicolumn{2}{c|}{$p=4$} & \multicolumn{2}{c}{$p=5$}  \\ 
           & ILUT  & GS  & ILUT   &GS   &  ILUT & GS  & ILUT & GS     \\ \hline 
$h=2^{-6}$ & $3$   &$10$& $2$	  &$16$ & 	$2$	&$26$ &$2$&  $52$  \\
$h=2^{-7}$ & $3$   &$10$& $2$	  &$17$ & 	$2$	&$32$ &$2$&  $57$  \\
$h=2^{-8}$ & $3$   &$10$& $2$	  &$17$ &   $2$	&$33$ &$2$&  $66$ \\
$h=2^{-9}$ & $4$   &$11$& $2$	  &$18$ &   $2$	&$36$ &$2$&  $64$ \\
\end{tabular}
\subcaption{Poisson's equation on L-shaped domain}
\end{subtable}
\begin{subtable}{1\textwidth}
\centering
\begin{tabular}{c|cc|cc|cc|cc}
     & \multicolumn{2}{c|}{$p=2$}  & \multicolumn{2}{c|}{$p=3$}  & \multicolumn{2}{c|}{$p=4$} & \multicolumn{2}{c}{$p=5$}  \\ 
           & ILUT  & GS  & ILUT &GS   &  ILUT & GS   & ILUT & GS  \\ \hline 
$h=2^{-2}$ & $2$   &$14$  & $2$	&$30$  & $2$  &$94$  &  $3$ &  $276$     \\
$h=2^{-3}$ & $2$   &$16$  & $2$	&$40$  & $2$  &$105$ &  $2$ &  $229$      \\
$h=2^{-4}$ & $2$   &$19$  & $2$	&$44$  & $2$  &$119$ &  $3$ &  $285$    \\
$h=2^{-5}$ & $2$   &$19$  & $2$ &$49$  & $3$  &$136$ &  $3$ &  $310$    \\
\end{tabular}
\subcaption{Poisson's equation on the unit cube}
\end{subtable}
\caption{Number of iterations needed to achieve convergence with BiCGSTAB, using $p$-multigrid as preconditioner.}
\label{tab:5}
\end{table}

\subsection*{$h$-Multigrid as a preconditioner}
\noindent The number of iterations needed to achieve convergence with $h$-multigrid can be found in Table \ref{tab:5h}. Note that, since the $h$-multigrid method is symmetric, a Conjugate Gradient (CG) method can be applied as a Krylov solver. In general, a single iteration performed with a BiCGSTAB method is twice as expensive compared to a single CG iteration. Results with a CG method have been added between brackets in Table \ref{tab:5h}. 

\noindent When applying $h$-multigrid as a preconditioner for a BiCGSTAB method, the number of iterations needed to achieve convergence is very similar to the use of $p$-multigrid as a preconditioner. Note that, the use of CG as outer Krylov solver, approximately doubles the number of iterations when ILUT is applied as a smoother. For Gauss-Seidel, the use of CG yields similar iteration numbers as the use of BiCGSTAB as a Krylov solver. However, due to the lower costs per iteration, the overall costs are lower when CG is applied.
\newpage

\begin{table}[h!]
\begin{subtable}{1\textwidth}
\centering
\begin{tabular}{c|cc|cc|cc|cc}
     & \multicolumn{2}{c|}{$p=2$}  & \multicolumn{2}{c|}{$p=3$}  & \multicolumn{2}{c|}{$p=4$} & \multicolumn{2}{c}{$p=5$}  \\ 
           & ILUT  & GS  & ILUT   &GS   &  ILUT & GS  & ILUT & GS     \\ \hline 
$h=2^{-6}$ & $2(4)$   &$13(15)$ &  $2(3)$	  &$18(23)$ & 	$2(3)$	&$41(42)$ &	$2(4)$  & $78(80)$\\
$h=2^{-7}$ & $2(4)$   &$12(15)$ &  $2(3)$	  &$20(23)$ & 	$2(3)$	&$41(42)$ &	$2(4)$  & $92(80)$\\
$h=2^{-8}$ & $3(5)$   &$13(16)$ &  $2(3)$	  &$19(23)$ &   $2(3)$	&$43(41)$ &	$2(4)$  & $95(78)$    \\
$h=2^{-9}$ & $3(5)$   &$13(16)$ &  $2(3)$	  &$21(23)$ &   $2(3)$	&$41(41)$ & $2(4)$  & $95(79)$    \\
\end{tabular}
\subcaption{Poisson's equation on quarter annulus}
\end{subtable}
\begin{subtable}{1\textwidth}
\centering
\begin{tabular}{c|cc|cc|cc|cc}
     & \multicolumn{2}{c|}{$p=2$}  & \multicolumn{2}{c|}{$p=3$}  & \multicolumn{2}{c|}{$p=4$} & \multicolumn{2}{c}{$p=5$}  \\ 
           & ILUT  & GS  & ILUT   &GS   &  ILUT & GS  & ILUT & GS     \\ \hline 
$h=2^{-6}$ & $2(4)$   &$7(11)$ & $2(3)$	  &$13(18)$ &	$2(3)$	&$29(35)$ & $2(4)$ & $65(68)$ \\
$h=2^{-7}$ & $2(4)$   &$8(11)$ & $2(3)$	  &$13(18)$ &	$2(4)$	&$29(33)$ & $2(4)$ & $70(66)$\\
$h=2^{-8}$ & $2(4)$   &$7(11)$ & $2(3)$	  &$12(18)$ &   $2(4)$	&$29(34)$ & $2(4)$ & $64(64)$  \\
$h=2^{-9}$ & $2(4)$   &$7(11)$ & $2(3)$	  &$14(18)$ &   $2(4)$	&$28(34)$ & $3(4)$ & $72(67)$  \\
\end{tabular}
\subcaption{CDR-equation on unit square}
\end{subtable}
\begin{subtable}{1\textwidth}
\centering
\begin{tabular}{c|cc|cc|cc|cc}
     & \multicolumn{2}{c|}{$p=2$}  & \multicolumn{2}{c|}{$p=3$}  & \multicolumn{2}{c|}{$p=4$} & \multicolumn{2}{c}{$p=5$}  \\ 
           & ILUT  & GS  & ILUT   &GS   &  ILUT & GS  & ILUT & GS     \\ \hline 
$h=2^{-6}$ & $3(5)$   &$10(15)$& $2(3)$	  &$16(22)$ & 	$2(3)$	&$26(35)$ &$2(4)$&  $52(62)$  \\
$h=2^{-7}$ & $3(5)$   &$10(15)$& $2(3)$	  &$17(22)$ & 	$2(3)$	&$32(35)$ &$2(3)$&  $57(61)$  \\
$h=2^{-8}$ & $3(6)$   &$10(16)$& $2(3)$	  &$17(22)$ &   $2(3)$	&$33(35)$ &$2(4)$&  $66(61)$ \\
$h=2^{-9}$ & $4(7)$   &$11(16)$& $2(3)$	  &$18(22)$ &   $2(3)$	&$36(36)$ &$2(4)$&  $64(61)$ \\
\end{tabular}
\subcaption{Poisson's equation on L-shaped domain}
\end{subtable}
\begin{subtable}{1\textwidth}
\centering
\begin{tabular}{c|cc|cc|cc|cc}
     & \multicolumn{2}{c|}{$p=2$}  & \multicolumn{2}{c|}{$p=3$}  & \multicolumn{2}{c|}{$p=4$} & \multicolumn{2}{c}{$p=5$}  \\ 
           & ILUT  & GS  & ILUT &GS   &  ILUT & GS   & ILUT & GS  \\ \hline 
$h=2^{-2}$ & $2(3)$   &$14(24)$  & $2(3)$	&$30(58)$  & $2(3)$  &$94(158)$  &  $3(7)$ &  $276(459)$     \\
$h=2^{-3}$ & $2(3)$   &$16(25)$  & $2(3)$	&$40(63)$  & $2(3)$  &$105(168)$ &  $2(3)$ &  $229(490)$      \\
$h=2^{-4}$ & $2(3)$   &$19(25)$  & $2(3)$	&$44(65)$  & $2(4)$  &$119(188)$ &  $3(8)$ &  $285(543)$    \\
$h=2^{-5}$ & $2(3)$   &$19(25)$  & $2(3)$   &$49(67)$  & $3(5)$  &$136(197)$ &  $3(12)$ & $310(591)$    \\
\end{tabular}
\subcaption{Poisson's equation on the unit cube}
\end{subtable}
\caption{Number of iterations needed to achieve convergence with BiCGSTAB (CG), using $h$-multigrid as preconditioner.}
\label{tab:5h}
\end{table}

\newpage

\subsubsection*{CPU times}

\begin{figure}
\begin{tikzpicture}[scale=0.85]
\definecolor{assembly}{RGB}{130,187,206}
\definecolor{factorize}{RGB}{61,152,222}
\definecolor{solve}{RGB}{34,70,122}

\pgfplotstableread{
name
pMG-ILUT
pMG-ILUT
pMG-ILUT
pMG-ILUT
hMG-ILUT
hMG-ILUT
hMG-ILUT
hMG-ILUT
pMG-GS
pMG-GS
pMG-GS
pMG-GS
hMG-GS
hMG-GS
hMG-GS
hMG-GS
}\datatable

    \begin{axis}[ybar stacked,
                 legend style={ legend columns=3,
                                at={(xticklabel cs:0.5)},
                                anchor=north,
                                draw=none},  
                 xtick=data,
                 bar width=2mm,
                 ymin=0,
                 axis y line*=none,
                 axis x line*=none,
                 xticklabels from table={\datatable}{name},
                 x tick label style={rotate=90,anchor=east,font=\tiny,color=black},
                 tick label style={font=\footnotesize},
                 legend style={font=\footnotesize,yshift=-3ex},
                 label style={font=\footnotesize},
                 xlabel style={yshift=-5ex},
                 ylabel={CPU time in seconds},
                 title =\textbf{\large $h=2^{-6}$},
                 area legend]    
                \addplot [assembly,fill=assembly,x tick label style={xshift=-0.3cm}] table[x=Clusters,y=assembly] {table_h1.txt};
                \addlegendentry[]{assembly};
                \addplot [factorize,fill=factorize,x tick label style={xshift=-0.3cm}] table[x=Clusters,y=factorize] {table_h1.txt};
                \addlegendentry{smoother(setup)};
                \addplot [solve,fill=solve,x tick label style={xshift=-0.3cm}] table[x=Clusters,y=solve] {table_h1.txt};
                \addlegendentry{solve};

    \end{axis}
        \draw (1.05,-1.75)  node{\small $p=2$};
        \draw (2.65,-1.75)  node{\small $p=3$};
        \draw (4.20,-1.75) node{\small $p=4$};
        \draw (5.80,-1.75) node{\small $p=5$};
\end{tikzpicture}\vspace{0.5cm}
\begin{tikzpicture}[scale=0.85]
\definecolor{assembly}{RGB}{130,187,206}
\definecolor{factorize}{RGB}{61,152,222}
\definecolor{solve}{RGB}{34,70,122}

\pgfplotstableread{
name
pMG-ILUT
pMG-ILUT
pMG-ILUT
pMG-ILUT
hMG-ILUT
hMG-ILUT
hMG-ILUT
hMG-ILUT
pMG-GS
pMG-GS
pMG-GS
pMG-GS
hMG-GS
hMG-GS
hMG-GS
hMG-GS
}\datatable

    \begin{axis}[ybar stacked,
                 legend style={ legend columns=3,
                                at={(xticklabel cs:0.5)},
                                anchor=north,
                                draw=none},  
                 xtick=data,
                 bar width=2mm,
                 ymin=0,
                 axis y line*=none,
                 axis x line*=none,
                 xticklabels from table={\datatable}{name},
                 x tick label style={rotate=90,anchor=east,font=\tiny,color=black},
                 tick label style={font=\footnotesize},
                 legend style={font=\footnotesize,yshift=-3ex},
                 label style={font=\footnotesize},
                 xlabel style={yshift=-5ex},
                 ylabel={CPU time in seconds},
                title =\textbf{\large $h=2^{-7}$},    
                 area legend]    
                \addplot [assembly,fill=assembly,x tick label style={xshift=-0.3cm}] table[x=Clusters,y=assembly] {table_h2.txt};
                \addlegendentry[]{assembly};
                \addplot [factorize,fill=factorize,x tick label style={xshift=-0.3cm}] table[x=Clusters,y=factorize] {table_h2.txt};
                \addlegendentry{smoother(setup)};
                \addplot [solve,fill=solve,x tick label style={xshift=-0.3cm}] table[x=Clusters,y=solve] {table_h2.txt};
                \addlegendentry{solve};

    \end{axis}
        \draw (1.05,-1.75)  node{\small $p=2$};
        \draw (2.65,-1.75)  node{\small $p=3$};
        \draw (4.20,-1.75) node{\small $p=4$};
        \draw (5.80,-1.75) node{\small $p=5$};
\end{tikzpicture}

\begin{tikzpicture}[scale=0.85]
\definecolor{assembly}{RGB}{130,187,206}
\definecolor{factorize}{RGB}{61,152,222}
\definecolor{solve}{RGB}{34,70,122}

\pgfplotstableread{
name
pMG-ILUT
pMG-ILUT
pMG-ILUT
pMG-ILUT
hMG-ILUT
hMG-ILUT
hMG-ILUT
hMG-ILUT
pMG-GS
pMG-GS
pMG-GS
pMG-GS
hMG-GS
hMG-GS
hMG-GS
hMG-GS
}\datatable

    \begin{axis}[ybar stacked,
                 legend style={ legend columns=3,
                                at={(xticklabel cs:0.5)},
                                anchor=north,
                                draw=none},  
                 xtick=data,
                 bar width=2mm,
                 ymin=0,
                 axis y line*=none,
                 axis x line*=none,
                 xticklabels from table={\datatable}{name},
                 x tick label style={rotate=90,anchor=east,font=\tiny,color=black},
                 tick label style={font=\footnotesize},
                 legend style={font=\footnotesize,yshift=-3ex},
                 label style={font=\footnotesize},
                 xlabel style={yshift=-5ex},
                title =\textbf{\large $h=2^{-8}$},    
                 ylabel={CPU time in seconds},
                 area legend]    
                \addplot [assembly,fill=assembly,x tick label style={xshift=-0.3cm}] table[x=Clusters,y=assembly] {table_h3.txt};
                \addlegendentry[]{assembly};
                \addplot [factorize,fill=factorize,x tick label style={xshift=-0.3cm}] table[x=Clusters,y=factorize] {table_h3.txt};
                \addlegendentry{smoother(setup)};
                \addplot [solve,fill=solve,x tick label style={xshift=-0.3cm}] table[x=Clusters,y=solve] {table_h3.txt};
                \addlegendentry{solve};

    \end{axis}
        \draw (1.05,-1.75)  node{\small $p=2$};
        \draw (2.65,-1.75)  node{\small $p=3$};
        \draw (4.20,-1.75) node{\small $p=4$};
        \draw (5.80,-1.75) node{\small $p=5$};
\end{tikzpicture}
\begin{tikzpicture}[scale=0.85]
\definecolor{assembly}{RGB}{130,187,206}
\definecolor{factorize}{RGB}{61,152,222}
\definecolor{solve}{RGB}{34,70,122}

\pgfplotstableread{
name
pMG-ILUT
pMG-ILUT
pMG-ILUT
pMG-ILUT
hMG-ILUT
hMG-ILUT
hMG-ILUT
hMG-ILUT
pMG-GS
pMG-GS
pMG-GS
pMG-GS
hMG-GS
hMG-GS
hMG-GS
hMG-GS
}\datatable

    \begin{axis}[ybar stacked,
                 legend style={ legend columns=3,
                                at={(xticklabel cs:0.5)},
                                anchor=north,
                                draw=none},  
                 xtick=data,
                 bar width=2mm,
                 ymin=0,
                 axis y line*=none,
                 axis x line*=none,
                 xticklabels from table={\datatable}{name},
                 x tick label style={rotate=90,anchor=east,font=\tiny,color=black},
                 tick label style={font=\footnotesize},
                 legend style={font=\footnotesize,yshift=-3ex},
                 label style={font=\footnotesize},
                 xlabel style={yshift=-5ex},
                 ylabel={CPU time in seconds},
                title = {\large $h=2^{-9}$},    
                 area legend]    
                \addplot [assembly,fill=assembly,x tick label style={xshift=-0.3cm}] table[x=Clusters,y=assembly] {table_h4.txt};
                \addlegendentry[]{assembly};
                \addplot [factorize,fill=factorize,x tick label style={xshift=-0.3cm}] table[x=Clusters,y=factorize] {table_h4.txt};
                \addlegendentry{smoother(setup)};
                \addplot [solve,fill=solve,x tick label style={xshift=-0.3cm}] table[x=Clusters,y=solve] {table_h4.txt};
                \addlegendentry{solve};

    \end{axis}
        \draw (1.05,-1.75)  node{\small $p=2$};
        \draw (2.65,-1.75)  node{\small $p=3$};
        \draw (4.20,-1.75) node{\small $p=4$};
        \draw (5.80,-1.75) node{\small $p=5$};
\end{tikzpicture}
\caption{CPU timings for $p$-multigrid and $h$-multigrid adopting ILUT and Gauss-Seidel (GS) as a smoother for different values of $h$ for the first benchmark.}
\label{fig:CPU}
\end{figure}
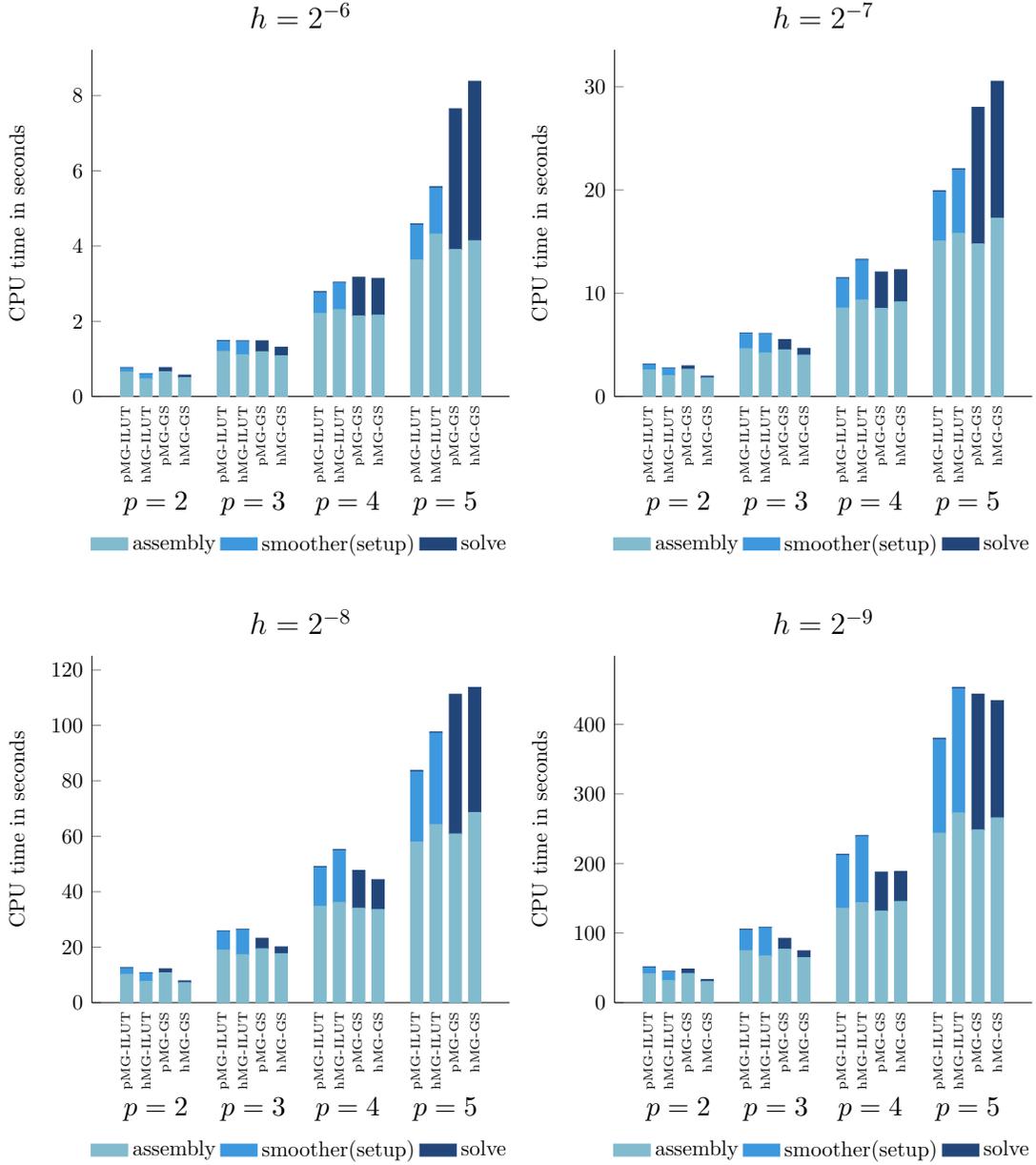

Besides iteration numbers, computational times have been determined when adopting $p$-multigrid and $h$-multigrid as a stand-alone solver. A serial implementation in the C++ library G+Smo \cite{gismo} is considered on a Intel(R) Core(TM) i7-8650U CPU (1.90GHz). Figure \ref{fig:CPU} illustrates the CPU times obtained for the $p$-multigrid and $h$-multigrid method for the first benchmark. Tables with the detailed CPU times can be found in \ref{cpu_pmg} and \ref{cpu_hmg}.

\noindent The assembly times denote the CPU time needed to assemble all operators, including the prolongation and restriction operators. Note that, for the $p$-multigrid method more operators have to be assembled. However, most of the operators in the $p$-multigrid method are assembled at level $p=1$, where the number of nonzero entries is significantly lower compared to the matrices resulting from high order discretizations. As a consequence the total assembly costs are lower with $p$-multigrid compared to $h$-multigrid for higher values of the approximation order $p$.

\noindent With respect to the setup costs of the smoother, similar observations can be made: For higher values of $p$, the ILUT factorization costs are significantly higher for the $h$-multigrid method. The time needed to solve linear systems is slightly lower for the $h$-multigrid methods, since the costs of a single V-cycle is lower compared to the $p$-multigrid method. When adopting Gauss-Seidel as a smoother, the time needed to solve the linear systems is significantly higher compared to the use of ILUT. However, since the factorizations costs are relatively high, the $p$-multigrid/$h$-multigrid methods using ILUT as a smoother are faster for only a limited amount of configurations.
\\
\\
\textbf{Remark 2:} For all numerical experiments, the 'coarse grid' operators of the multigrid hierarchy have been obtained by rediscretizing the bilinear form in Equation \eqref{weak}). Alternatively, all operators of the $h$-multigrid hierarchy could be obtained by applying the Galerkin projection. Furtermore, alternative (and more efficient) assembly strategies exist, as mentioned in Section \ref{pm}. Therefore, the assembly, smoother setup and solving costs are presented separately in this Section. 

\newpage 

\subsubsection*{Comparison with an alternative smoother}
Throughout this paper, the use of ILUT and Gauss-Seidel as a smoother has been investigated within a $p$-multigrid and $h$-multigrid method. However, alternative smoothers have been developed for $h$-multigrid methods in recent years. For example, a smoother based on stable splittings of spline spaces \cite{tak1}. In this section, we compare the $p$-multigrid (adopting ILUT as a smoother) with an $h$-multigrid method using this smoother. Note that this smoother has been extended to multipatch domains as well \cite{tak2}. For this comparison, we consider the CDR-equation on the unit square with:
\begin{eqnarray}
\mathbf{D} = \begin{bmatrix*}[r]  1 & 0 \\ 0 & 1 \\ \end{bmatrix*}, \hspace{0.5cm} \mathbf{v} = \begin{bmatrix*}[r] 0  \\ 0 \\  \end{bmatrix*}, \hspace{0.5cm} R = 1.  
\end{eqnarray}

\noindent Homogeneous Neumann boundary conditions are applied and the right-hand side is given by:
\vskip-.6cm
\begin{eqnarray}
f(x,y) = 2 \pi^2 \text{sin}(\pi (x+ \frac{1}{2}))\text{sin}(\pi (y+ \frac{1}{2})). \nonumber
\end{eqnarray}       
\noindent Table \ref{tab:6} shows the number of iterations needed to reach convergence with the $p$-multigrid method and the $h$-multigrid method based on stable splittings. Both methods show iteration numbers which are independent of both $h$ and $p$. With the $p$-multigrid method, the number of iterations needed to reach convergence is significantly lower compared to the considered $h$-multigrid method for all configurations. 
\begin{table}[h]
\centering
\begin{tabular}{c|cc|cc|cc|cc}
& \multicolumn{2}{c|}{$p=2$}  & \multicolumn{2}{c|}{$p=3$}  & \multicolumn{2}{c|}{$p=4$} & \multicolumn{2}{c}{$p=5$}  \\ 
           & $p$-MG   & $h$-MG  & $p$-MG   &$h$-MG  & $p$-MG & $h$-MG & $p$-MG & $h$-MG      \\ \hline 
$h=2^{-6}$ & $5$   	  &$48$ & $5$ &$48$ & $5$&$48$  &$5$     &	$48$  \\
$h=2^{-7}$ & $5$   	  &$49$ & $5$ &$50$ & $4$&$49$  &$5$     &	$49$  \\
$h=2^{-8}$ & $5$      &$49$ & $4$ &$50$ & $5$&$50$  &$4$     &  $49$ \\
$h=2^{-9}$ & $5$      &$49$ & $4$ &$50$ & $5$&$50$  &$4$     &  $50$  \\
\end{tabular}
\caption{Comparison of a $p$-multigrid using ILUT as a smoother with a $h$-multigrid method based on stable splittings of spline spaces \cite{tak1}.}
\label{tab:6}
\end{table}

\begin{figure}[h!]
\begin{tikzpicture}[yscale = 0.85, xscale=0.85]
\definecolor{assembly}{RGB}{130,187,206}
\definecolor{factorize}{RGB}{61,152,222}
\definecolor{solve}{RGB}{34,70,122}

\pgfplotstableread{
name
pMG-ILUT
hMG-SCMS
pMG-ILUT
hMG-SCMS
pMG-ILUT
hMG-SCMS
pMG-ILUT
hMG-SCMS
}\datatable

    \begin{axis}[ybar stacked,
                 legend style={ legend columns=3,
                                at={(xticklabel cs:0.5)},
                                anchor=north,
                                draw=none},  
                 xtick=data,
                 bar width=4mm,
                 ymin=0,
                 axis y line*=none,
                 axis x line*=none,
                 xticklabels from table={\datatable}{name},
                 x tick label style={rotate=90,anchor=east,font=\tiny,color=black},
                 tick label style={font=\footnotesize},
                 legend style={font=\footnotesize,yshift=-3ex},
                 label style={font=\footnotesize},
                 xlabel style={yshift=-5ex},
                 ylabel={CPU time in seconds},
                 title = {\large $h=2^{-7}$},    
                 area legend]    
                \addplot [assembly,fill=assembly,x tick label style={xshift=-0.3cm}] table[x=Clusters,y=assembly] {table_h5.txt};
                \addlegendentry[]{assembly};
                \addplot [factorize,fill=factorize,x tick label style={xshift=-0.3cm}] table[x=Clusters,y=factorize] {table_h5.txt};
                \addlegendentry{smoother(setup)};
                \addplot [solve,fill=solve,x tick label style={xshift=-0.3cm}] table[x=Clusters,y=solve] {table_h5.txt};
                \addlegendentry{solve};

    \end{axis}
        \draw (0.90,-1.75)  node{\small $p=2$};
        \draw (2.60,-1.75)  node{\small $p=3$};
        \draw (4.30,-1.75) node{\small $p=4$};
        \draw (6.00,-1.75) node{\small $p=5$};
\end{tikzpicture}
\begin{tikzpicture}[yscale = 0.85, xscale=0.85]
\definecolor{assembly}{RGB}{130,187,206}
\definecolor{factorize}{RGB}{61,152,222}
\definecolor{solve}{RGB}{34,70,122}

\pgfplotstableread{
name
pMG-ILUT
hMG-SCMS
pMG-ILUT
hMG-SCMS
pMG-ILUT
hMG-SCMS
pMG-ILUT
hMG-SCMS
}\datatable

    \begin{axis}[ybar stacked,
                 legend style={ legend columns=3,
                                at={(xticklabel cs:0.5)},
                                anchor=north,
                                draw=none},  
                 xtick=data,
                 bar width=4mm,
                 ymin=0,
                 axis y line*=none,
                 axis x line*=none,
                 xticklabels from table={\datatable}{name},
                 x tick label style={rotate=90,anchor=east,font=\tiny,color=black},
                 tick label style={font=\footnotesize},
                 legend style={font=\footnotesize,yshift=-3ex},
                 label style={font=\footnotesize},
                 xlabel style={yshift=-5ex},
                 ylabel={CPU time in seconds},
                 title = {\large $h=2^{-9}$},     
                 area legend]    
                \addplot [assembly,fill=assembly,x tick label style={xshift=-0.3cm}] table[x=Clusters,y=assembly] {table_h6.txt};
                \addlegendentry[]{assembly};
                \addplot [factorize,fill=factorize,x tick label style={xshift=-0.3cm}] table[x=Clusters,y=factorize] {table_h6.txt};
                \addlegendentry{smoother(setup)};
                \addplot [solve,fill=solve,x tick label style={xshift=-0.3cm}] table[x=Clusters,y=solve] {table_h6.txt};
                \addlegendentry{solve};

    \end{axis}
        \draw (0.90,-1.75)  node{\small $p=2$};
        \draw (2.60,-1.75)  node{\small $p=3$};
        \draw (4.30,-1.75) node{\small $p=4$};
        \draw (6.00,-1.75) node{\small $p=5$};
\end{tikzpicture}
\caption{CPU timings for $p$-multigrid and $h$-multigrid adopting ILUT and the smoother from \cite{tak1}, respectively, for $h=2^{-7}$ (left) and $h=2^{-9}$ (right).}
\label{fig:CPU_SCMS}
\end{figure}
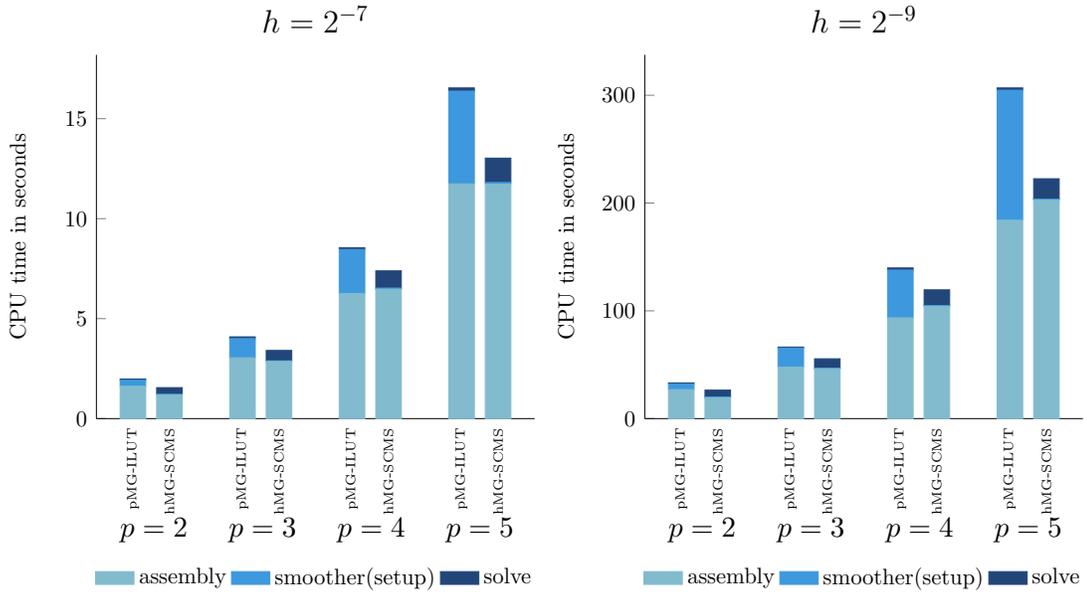

\noindent CPU times for assembly, setting up the smoother and solving the linear system are presented in Figure \ref{fig:CPU_SCMS}. Again, a serial implementation in the C++ library G+Smo \cite{gismo} is considered on a Intel(R) Core(TM) i7-8650U CPU (1.90GHz). Detailed CPU times can be found in \ref{cpu_scms}. The time needed to assembly the operators is comparable for the $p$-multigrid and $h$-multigrid method. However, setting up the ILUT smoother is significantly more expensive compared to the smoother from \cite{tak1}. On the other hand, the CPU needed to solve the problem is lower when adopting the $p$-multigrid method. The total solver costs are lower for all configurations when adopting the smoother based on stable splitting of subspaces. However, we would like to emphasize that the proposed $p$-multigrid method can easily be implemented and applied for a wide variety of problems (multipatch, variable coefficients) without the need of tuning a parameter or development of a specific smoother.  

\newpage

\section*{Truncated Hierarchical B-splines (THB-splines)}
Finally, to illustrate the versatility of the proposed $p$-multigrid method, we consider discretizations obtained with THB-splines \cite{juttler}. THB-splines are the result of a local refinement strategy, in which a subset of the basis functions on the fine level are truncated. As a result, not only linear indepen-dence and non-negativity are preserved (as with HB-splines \cite{kraft,vuong}), but also the partition of unity property.

\begin{figure}[h!]
\centering
\begin{tikzpicture}[scale=0.21]
\draw[] node at (28,-2){(a)};
\draw[] node at (48,-2){(b)};
\filldraw[fill=green!60!white, draw=black] (24,4) rectangle (32,12);
\filldraw[fill=orange!60!white, draw=black] (26,6) rectangle (30,10);
\draw[step=4cm,gray,very thin] (20,0) grid (36,16);
\draw[step=2cm,gray,very thin] (24,4) grid (32,12);
\draw[step=1cm,gray,very thin] (26,6) grid (30,10);
\filldraw[fill=green!60!white, draw=black] (44,4) rectangle (52,12);
\filldraw[fill=orange!60!white, draw=black] (46,6) rectangle (50,10);
\filldraw[fill=red!60!white, draw=black] (47,7) rectangle (49,9);

\draw[step=4cm,gray,very thin] (40,0) grid (56,16);
\draw[step=2cm,gray,very thin] (44,4) grid (52,12);
\draw[step=1cm,gray,very thin] (46,6) grid (50,10);
\draw[step=0.5cm,gray,very thin] (47,7) grid (49,9);
\end{tikzpicture}
\caption{Two hierarchical mesh adopted for THB-Spline basis with the second (green), third (orange) and fourth (red) refinement levels coloured.}
\label{fig:mesh}
\end{figure}
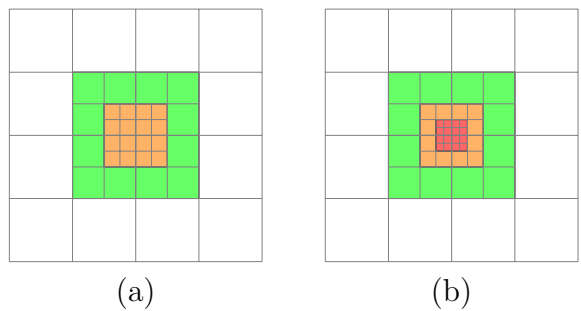

\noindent In the literature, the use of multigrid methods for THB-spline discretizations is limited and an ongoing topic of research \cite{thb,bpx,localmultigrid}. We consider Poisson's equation on the unit square, where the exact solution is the same as for the second benchmark. Starting from a tensor product B-spline basis with meshwidth $h$ and order $p$, two and three levels of refinement are added as shown in Figure \ref{fig:mesh}, leading to a THB-spline basis consisting of, respectively, three and four levels. 

\begin{figure}[h!]
\begin{center}
\includegraphics[scale=0.29]{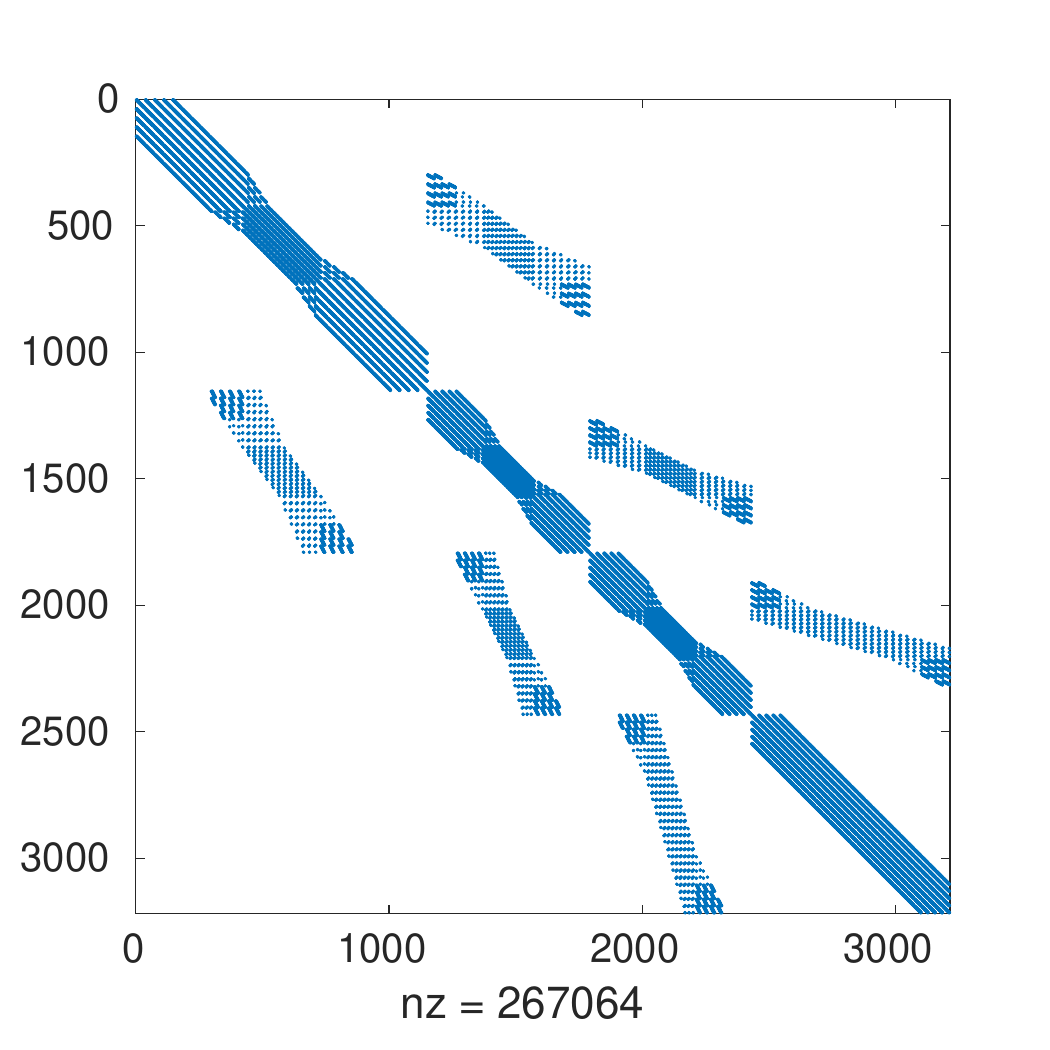}
\includegraphics[scale=0.29]{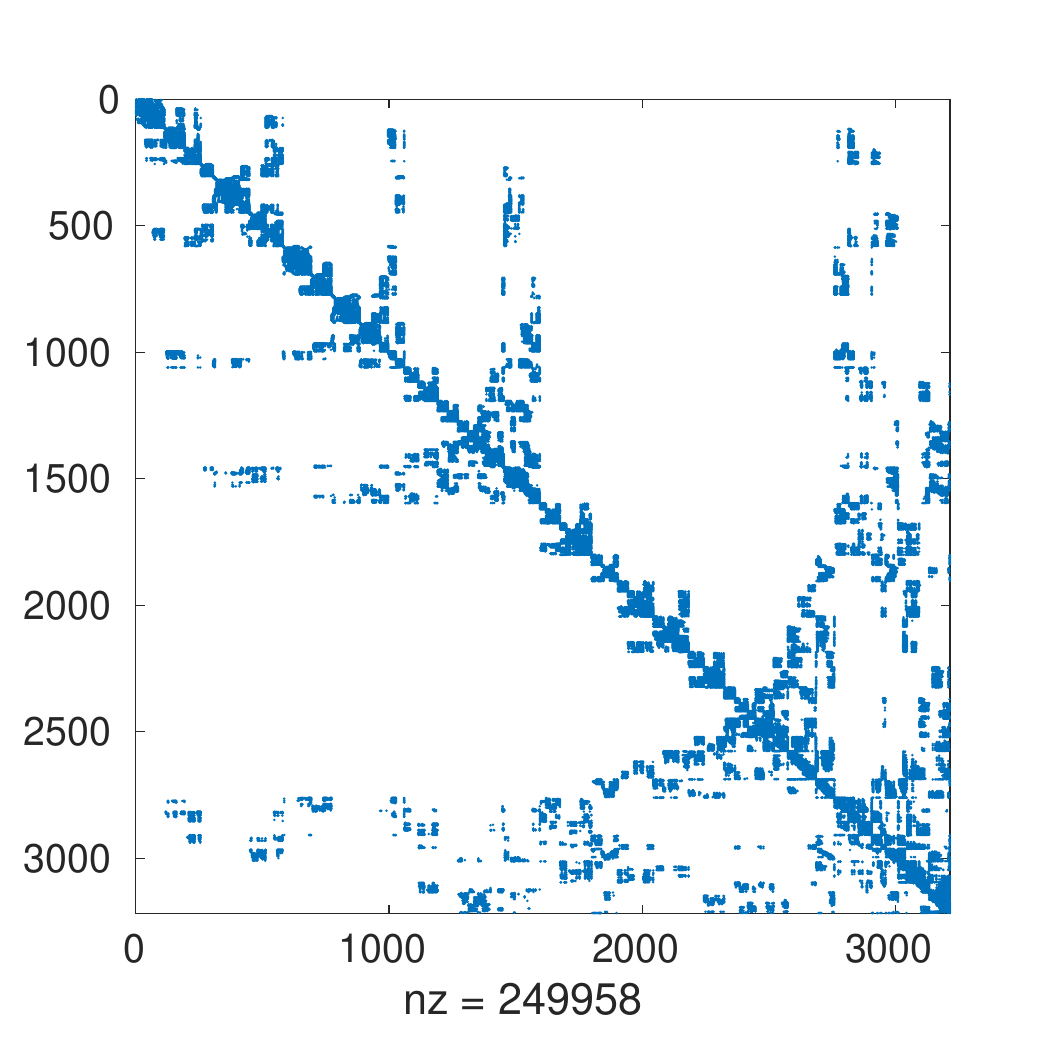}
\end{center} 
\caption{Sparsity pattern of the stiffness matrix $\mathbf{A}_{h,4}$ (left) and $\mathbf{L}_{h,4} + \mathbf{U}_{h,4}$ (right).}
\label{fig:8}
\end{figure}

\noindent Figure \ref{fig:8} shows the sparsity pattern of the stiffness matrix and the ILUT factorization for $p=4$ and $h=2^{-5}$ for configuration (b). Compared to the (standard) tensor-product B-spline basis the bandwith of the stiffness matrix significantly increases. Table \ref{tab:8} shows the results obtained with $p$-multigrid applied as a stand-alone solver. The number of iterations needed with $p$-multigrid (and ILUT as a smoother) depends only mildly on $p$. Furthermore, the number of iterations are significantly lower compared to the use of Gauss-Seidel as a smoother. 
 
\begin{figure}[h!]
\begin{center}
\includegraphics[scale=0.31]{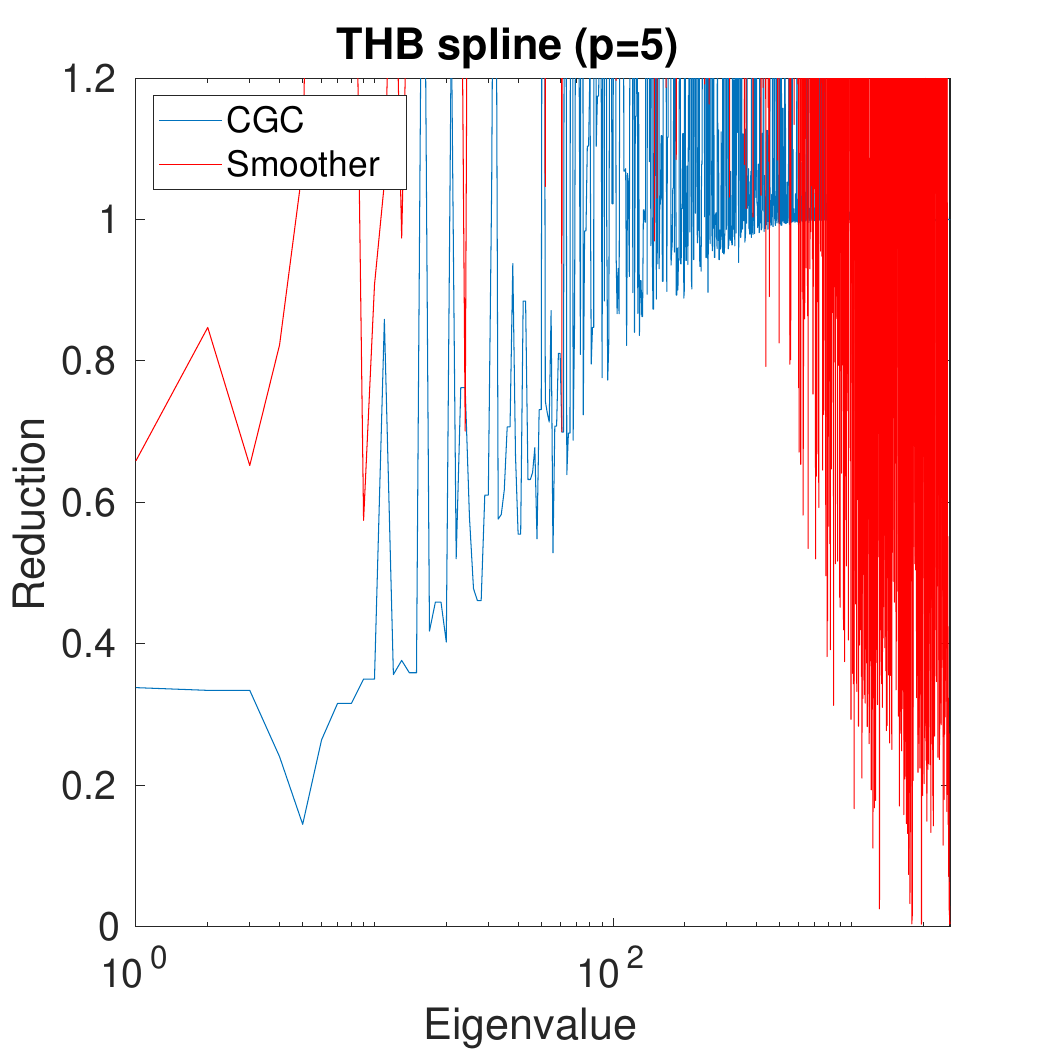}
\includegraphics[scale=0.31]{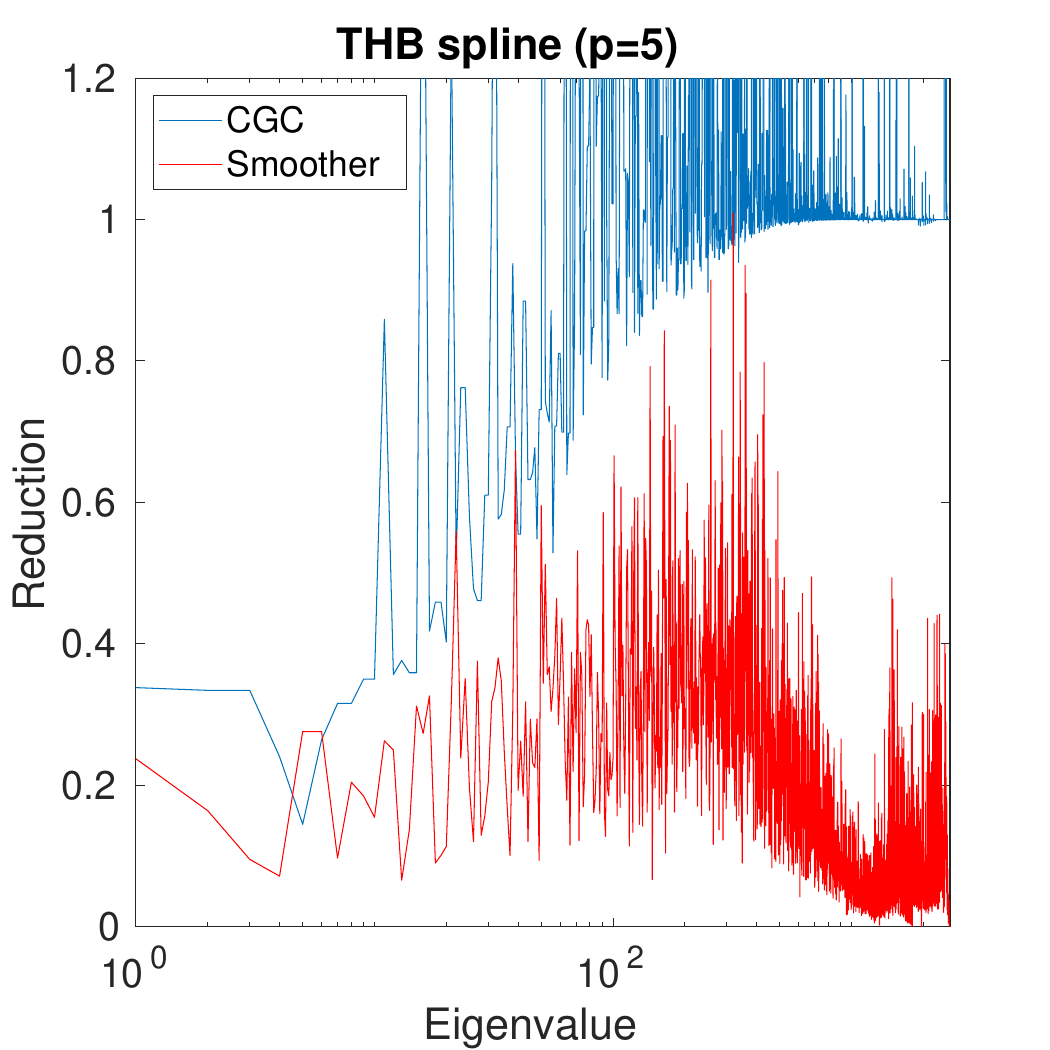}
\end{center}
\caption{Reduction factors obtained for fillfactor $1$ (left) and $2$ (right).}
\label{fig:7}
\end{figure}

\begin{table}[h!]
\begin{subtable}{\textwidth}
\centering
\begin{tabular}{c|cc|cc|cc|cc}
     & \multicolumn{2}{c|}{$p=2$}  & \multicolumn{2}{c|}{$p=3$}  & \multicolumn{2}{c|}{$p=4$} & \multicolumn{2}{c}{$p=5$}  \\ 
           & ILUT & GS  & ILUT& GS  & ILUT & GS   & ILUT & GS    \\ \hline 
$h=2^{-4}$ & $5$  &$16$ & $6$ &$45$ & $5$  &$178$ &$5$   &$713$ \\
$h=2^{-5}$ & $5$  &$17$ & $6$ &$40$ & $7$  &$182$ &$\mathbf{5}$   &$882$ \\
$h=2^{-6}$ & $5$  &$17$ & $5$ &$41$ & $7$  &$189$ &$\mathbf{11}$   &$936$ \\
\end{tabular}
\subcaption{THB-spline basis with tree levels of refinement.}
\end{subtable}
\begin{subtable}{\textwidth}
\centering
\begin{tabular}{c|cc|cc|cc|cc}
     & \multicolumn{2}{c|}{$p=2$}  & \multicolumn{2}{c|}{$p=3$}  & \multicolumn{2}{c|}{$p=4$} & \multicolumn{2}{c}{$p=5$}  \\ 
           & ILUT & GS  & ILUT& GS  & ILUT & GS   & ILUT & GS    \\ \hline 
$h=2^{-4}$ & $6$  &$17$ & $8$ &$47$ & $7$    &$177$   &$10$             &$1033$ \\
$h=2^{-5}$ & $6$  &$16$ & $7$ &$44$ & $8$    &$182$   &$\mathbf{7}$   &$923$ \\
$h=2^{-6}$ & $6$  &$17$ & $5$ &$43$ & $6$    &$201$   &$\mathbf{12}$  &$1009$ \\
\end{tabular}
\subcaption{THB-spline basis with four levels of refinement.}
\end{subtable}
\caption{Number of V-cycles needed for different THB-spline discretizations.}
\label{tab:8}
\end{table}
\noindent For the configurations denoted in bold, a fillfactor of $2$ was adopted, to prevent the $p$-multigrid from diverging. Figure \ref{fig:7} illustrates the reason for it in the case $p=4$ and $h=2^{-4}$ for configuration (a). A fillfactor of $1$ does not reduce the norm of the (generalized) eigenvectors, while a fillfactor of $2$ reduces the eigenvectors over the entire spectrum. In general, a higher fillfactor was necessary for only a limited amount of configurations.   

\noindent For all numerical experiments, smoothing is performed globally at each level of the multigrid hierarchy. In general, local smoothing is often adopted to ensure optimal order of the complexity. Results presented in this Section should be considered as a first step towards the use of $p$-multigrid methods for THB-spline discretizations. Future research should focus on more efficient applications of $p$-multigrid solvers for THB-spline discretizations. 

\section{Conclusions}
\label{con}
In this paper, we presented a $p$-multigrid method that uses ILUT factori-zation as a smoother and compared this with different smoothers and coarse-ning strategy (e.g. $h$-multigrid). In contrast to classical smoothers, (i.e. Gauss-Seidel), the reduction factors of the general eigenvectors associated with high-frequency modes do not increase when adopting ILUT as a smoother for higher values of $p$. This results in asymptotic convergence factors which are independent of both the mesh width $h$ and approximation order $p$ for both $p$-multigrid and $h$-multigrid methods adopting this smoother. Furthermore, we observed that, assuming an exact coarse grid correction, coarsening in $h$ leads to a more effective coarse grid correction compared to a correction obtained by coarsening in $p$.

\noindent Numerical results, obtained for Poisson's equation on a variety of domains and the CDR equation on the unit square have been presented when using $p$-multigrid and $h$-multigrid as stand-alone solver or as a preconditioner within a BiCGSTAB or CG method. For all configurations, the number of iterations needed when using ILUT as a smoother are significantly lower compared to the use of Gauss-Seidel, while the number of iterations needed with $p$-multigrid are very similar to those needed with an $h$-multigrid method. Hence, the smoother determines to a great extent the resulting convergence rate of the multigrid method. CPU times have been presented for the $p$-multigrid and $h$-multigrid method using both smoothers. For low values of $p$, the use of $h$-multigrid combined with Gauss-Seidel as a smoother lead to the lowest CPU times. For higher values of $p$, however, the use of $p$-multigrid adopting ILUT becomes more efficient, due to the lower assembly and factorizations costs. Note that this is the result of the smaller stencil of the B-spline functions at level $p=1$ compared to high order B-spline functions.   

\noindent The $p$-multigrid method using ILUT as a smoother has been compared as well to an $h$-multigrid method with a non-standard smoother \cite{tak1}. Result show that the total solving costs are lower when adopting $h$-multigrid with this smoother due to the lower setup costs of the smoother. Finally, the $p$-multigrid method has been succesfully applied to solve linear systems of equations arising from THB-spline discretizations. In general, a significantly lower number of iterations was needed compared to the use of Gauss-Seidel as a smoother. For a limited number of configurations, a higher fillfactor of $2$ (instead of $1$) was necessary to achieve convergence. 

\noindent Future research will focus on the application of $p$-multigrid methods for higher-order partial differential equations (i.e. biharmonic equation), where the use of basis functions with high continuity is necessary. Furthermore, local smoothing within the $p$-multigrid method should be considered to make it more efficient for THB-spline discretizations. Finally, the use of block ILUT as a smoother in case of a multipatch geometry will be investigated.

\section{Acknowledgements}
The authors would like to thank Prof. Kees Oosterlee from TU Delft for fruitful discussions with respect to $p$-multigrid methods.

\newpage

\appendix

\section{Direct or indirect projection}
\label{direct}
To investigate the effect of a direct projection to $p=1$, we consider the first benchmark. The number of V-cycles needed to achieve convergence with a direct projection and indirect projection has been determined for different values of $h$ and $p$. Table \ref{tab:direct_projection} and \ref{tab:indirect_projection} show the number of iterations needed to achieve convergence with a direct and indirect projection, respectively. For most configurations, the number of iterations is very similar. Only for higher values of $p$, the indirect project leads to diverging method when Gauss-Seidel is applied as a smoother. With a direct projection, all configurations lead to a converging multigrid method.

\begin{table}[h!]
\centering
\begin{tabular}{c|cc|cc|cc|cc}
     & \multicolumn{2}{c|}{$p=2$}  & \multicolumn{2}{c|}{$p=3$}  & \multicolumn{2}{c|}{$p=4$} & \multicolumn{2}{c}{$p=5$}  \\ 
           & ILUT  & GS  & ILUT &GS   &  ILUT & GS   & ILUT & GS  \\ \hline 
$h=2^{-6}$ & $4$   &$30$ & $3$	&$62$ & $3$	  &$176$	 &  $3$ &  $491$    \\
$h=2^{-7}$ & $4$   &$29$ & $3$	&$61$ & $3$	  &$172$	 &  $3$ &  $499$    \\
$h=2^{-8}$ & $5$   &$30$ & $3$	&$60$ & $3$	  &$163$     &	$3$ &  $473$    \\
$h=2^{-9}$ & $5$   &$32$ & $3$	&$61$ & $3$	  &$163$	 &  $3$ &  $452$    \\
\end{tabular}
\centering
\caption{Number of V-cycles needed to achieve convergence with $p$-multigrid for the first benchmark with a direct projection.}
\label{tab:direct_projection}
\end{table}

\begin{table}[h!]
\centering
\begin{tabular}{c|rr|rr|rr|rr}
     & \multicolumn{2}{c|}{$p=2$}  & \multicolumn{2}{c|}{$p=3$}  & \multicolumn{2}{c|}{$p=4$} & \multicolumn{2}{c}{$p=5$}  \\ 
           & ILUT  & GS  & ILUT    & GS  & ILUT   & GS  & ILUT & GS     \\ \hline 
$h=2^{-6}$ & $4$  &$30$ & $3$  &$62$ & $3$	&$-$   &$3$  & $-$\\
$h=2^{-7}$ & $4$  &$29$ & $3$  &$61$ & $3$	&$-$   &$3$  & $-$\\
$h=2^{-8}$ & $5$  &$30$ & $3$  &$61$ & $3$	&$-$   &$3$  & $-$ \\
$h=2^{-9}$ & $5$  &$32$ & $3$  &$63$ & $3$	&$-$   &$3$  & $-$\\
\end{tabular}
\centering
\caption{Number of V-cycles needed to achieve convergence with $p$-multigrid for the first benchmark with an indirect projection.}
\label{tab:indirect_projection}
\end{table}

\newpage

\section{Consistent vs. lumped projection}
\label{lumped}

\noindent In Section \ref{mp}, the prolongation and restriction operator to transfer residuals and corrections from level $p$ to $1$ and vice versa have been defined. Note that, the mass matrix in Equation \eqref{prolongation} and \eqref{restriction} can be lumped to reduce computational costs. To investigate the effect of lumping the mass matrix within the $L_2$ projection, the first benchmark is considered. 

\noindent Table \ref{tab:lumped} shows the number of V-cycles needed to achieve convergence using the lumped or consistent mass matrix in Equation \eqref{prolongation} and \eqref{restriction}. When ILUT is adopted as a smoother, the number of V-cycles needed to reach convergence is identical for almost all configurations. For Gauss-Seidel, the use of the consistent mass matrix leads to a slightly lower number of iterations. 
Considering the decrease of computational costs, however, the lumped mass matrix is adopted throughout the entire paper in the prolongation and restriction operator. 

\begin{table}[h!]
\centering
\begin{subtable}{\textwidth}
\centering
\begin{tabular}{c|cc|cc|cc|cc}
     & \multicolumn{2}{c|}{$p=2$}  & \multicolumn{2}{c|}{$p=3$}  & \multicolumn{2}{c|}{$p=4$} & \multicolumn{2}{c}{$p=5$}  \\ 
           & ILUT  & GS  & ILUT    & GS  & ILUT   & GS  & ILUT & GS     \\ \hline 
$h=2^{-6}$ & $4$  &$30$ & $3$  &$62$ & $3$	&$176$   &$3$  & $491$\\
$h=2^{-7}$ & $4$  &$29$ & $3$  &$61$ & $3$	&$172$   &$3$  & $499$\\
$h=2^{-8}$ & $5$  &$30$ & $3$  &$61$ & $3$	&$163$   &$3$  & $473$ \\
$h=2^{-9}$ & $5$  &$32$ & $3$  &$63$ & $3$	&$163$   &$3$  & $452$\\
\end{tabular}
\subcaption{Lumped mass matrix $\mathbf{M}_k^L$.}
\end{subtable}
\begin{subtable}{\textwidth}
\centering
\begin{tabular}{c|cc|cc|cc|cc}
     & \multicolumn{2}{c|}{$p=2$}  & \multicolumn{2}{c|}{$p=3$}  & \multicolumn{2}{c|}{$p=4$} & \multicolumn{2}{c}{$p=5$}  \\ 
           & ILUT  & GS  & ILUT    & GS  & ILUT   & GS  & ILUT & GS     \\ \hline 
$h=2^{-6}$ & $4$  &$29$ & $3$  &$56$ & $3$	&$170$   &$3$  & $474$\\
$h=2^{-7}$ & $4$  &$29$ & $3$  &$51$ & $3$	&$173$   &$3$  & $523$\\
$h=2^{-8}$ & $4$  &$29$ & $3$  &$53$ & $3$	&$164$   &$3$  & $445$ \\
$h=2^{-9}$ & $5$  &$31$ & $3$  &$51$ & $3$	&$163$   &$3$  & $440$\\
\end{tabular}
\subcaption{Consistent mass matrix $\mathbf{M}_k$.}
\end{subtable}
\caption{Number of V-cycles to reach convergence with $p$-multigrid adopting a lumped or consistent mass matrix in the prolongation and restriction operator.}
\label{tab:lumped}
\end{table}

\newpage

Although lumping the mass matrix in Equation \eqref{prolongation} and \eqref{restriction} slightly influences the number of iterations needed to achieve convergence with the $p$-multigrid method, the overall accuracy of the $p$-multigrid method is not affected. To illustrate this, the $L_2$ error is presented in Table \ref{tab:lumped_norm} for all configurations when adopting a lumped or consistent mass matrix in the prolongation and restriction operator. 

\begin{table}[h!]
\centering
\begin{subtable}{\textwidth}
\centering
\begin{tabular}{c|cc|cc}
     & \multicolumn{2}{c|}{$p=2$}  & \multicolumn{2}{c}{$p=3$}    \\ 
           & ILUT  & GS  & ILUT    & GS   \\ \hline 
$h=2^{-6}$ & $1.29\cdot10^{-5}$  &$1.28\cdot10^{-5}$ & $1.03\cdot10^{-6}$  &$4.78\cdot10^{-7}$ \\
$h=2^{-7}$ & $1.98\cdot10^{-6}$  &$1.70\cdot10^{-6}$ & $1.35\cdot10^{-6}$  &$6.62\cdot10^{-7}$ \\
$h=2^{-8}$ & $6.48\cdot10^{-7}$  &$5.49\cdot10^{-7}$ & $1.67\cdot10^{-6}$  &$4.87\cdot10^{-7}$  \\
$h=2^{-9}$ & $6.37\cdot10^{-7}$  &$3.68\cdot10^{-7}$ & $4.13\cdot10^{-6}$  &$3.24\cdot10^{-7}$ \\
\end{tabular}
\begin{tabular}{c|cc|cc}
           & \multicolumn{2}{c|}{$p=4$} & \multicolumn{2}{c}{$p=5$} \\
           & ILUT  & GS  & ILUT    & GS     \\ \hline
$h=2^{-6}$ & $4.45\cdot10^{-7}$	&$4.22\cdot10^{-6}$   &$1.83\cdot10^{-6}$  & $1.47\cdot10^{-6}$ \\
$h=2^{-7}$ & $2.69\cdot10^{-7}$	&$3.06\cdot10^{-6}$   &$1.44\cdot10^{-6}$  & $1.01\cdot10^{-5}$ \\
$h=2^{-8}$ & $4.58\cdot10^{-7}$	&$1.93\cdot10^{-6}$   &$2.42\cdot10^{-7}$  & $7.15\cdot10^{-6}$ \\
$h=2^{-9}$ & $2.03\cdot10^{-6}$	&$1.52\cdot10^{-6}$   &$1.83\cdot10^{-6}$  & $5.06\cdot10^{-6}$ 
\end{tabular}
\subcaption{Lumped mass matrix $\mathbf{M}_k^L$.}
\end{subtable}
\begin{subtable}{\textwidth}
\centering
\begin{tabular}{c|cc|cc}
     & \multicolumn{2}{c|}{$p=2$}  & \multicolumn{2}{c}{$p=3$}    \\ 
           & ILUT  & GS  & ILUT    & GS     \\ \hline 
$h=2^{-6}$ & $1.29\cdot10^{-5}$  &$1.28\cdot10^{-5}$ & $2.79\cdot10^{-7}$  &$6.69\cdot10^{-7}$ \\
$h=2^{-7}$ & $1.62\cdot10^{-6}$  &$1.64\cdot10^{-6}$ & $2.76\cdot10^{-7}$  &$7.14\cdot10^{-7}$ \\
$h=2^{-8}$ & $1.50\cdot10^{-6}$  &$4.75\cdot10^{-7}$ & $7.03\cdot10^{-7}$  &$4.58\cdot10^{-7}$ \\
$h=2^{-9}$ & $1.56\cdot10^{-7}$  &$3.32\cdot10^{-7}$ & $2.86\cdot10^{-6}$  &$3.27\cdot10^{-7}$ \\
\end{tabular}
\begin{tabular}{c|cc|cc}
             & \multicolumn{2}{c|}{$p=4$} & \multicolumn{2}{c}{$p=5$} \\
             & ILUT  & GS  & ILUT    & GS     \\ \hline
$h=2^{-6}$ & $3.32\cdot10^{-7}$	&$4.24\cdot10^{-6}$   &$1.85\cdot10^{-6}$  & $1.20\cdot10^{-6}$\\
$h=2^{-7}$ & $1.63\cdot10^{-7}$	&$3.16\cdot10^{-6}$   &$1.44\cdot10^{-6}$  & $7.71\cdot10^{-6}$\\
$h=2^{-8}$ & $2.10\cdot10^{-7}$	&$1.84\cdot10^{-6}$   &$2.36\cdot10^{-6}$  & $7.33\cdot10^{-6}$ \\
$h=2^{-9}$ & $1.52\cdot10^{-6}$	&$1.51\cdot10^{-6}$   &$1.73\cdot10^{-6}$  & $4.95\cdot10^{-6}$\\
\end{tabular}
\subcaption{Consistent mass matrix $\mathbf{M}_k$.}
\end{subtable}
\caption{$L_2$ error after convergence with $p$-multigrid adopting a lumped or consistent mass matrix for prolongation and restriction.}
\label{tab:lumped_norm}
\end{table}

\newpage


\section{CPU times $p$-multigrid}
\label{cpu_pmg}

Table \ref{tab:CPU_p} shows the CPU timings with $p$-multigrid as a stand-alone solver for the first benchmark. For each configuration, the assembly, factorization and solver costs are shown separately. 
\begin{table}[h!]
\begin{subtable}{1\textwidth}
\centering
\begin{tabular}{c|rr|rr|rr|rr}
     & \multicolumn{2}{c|}{$p=2$}  & \multicolumn{2}{c|}{$p=3$}  & \multicolumn{2}{c|}{$p=4$} & \multicolumn{2}{c}{$p=5$}  \\ 
           & ILUT  & GS  & ILUT   &GS   &  ILUT & GS  & ILUT & GS     \\ \hline 
$h=2^{-6}$ & $0.65$   &$0.66$  &  $1.20$   &$1.19$ &   $2.21$	&$2.14$  & $3.63$   & $3.91$\\
$h=2^{-7}$ & $2.58$   &$2.63$  &  $4.64$   &$4.52$ &   $8.58$	&$8.53$  & $15.06$  & $14.78$\\
$h=2^{-8}$ & $10.22$  &$10.77$ &  $18.93$  &$19.41$&   $34.77$	&$34.07$ & $57.95$  & $60.81$ \\
$h=2^{-9}$ & $41.52$  &$41.86$ &  $74.74$  &$76.70$&   $135.79$ &$131.66$& $243.64$ & $248.11$ \\
\end{tabular}
\subcaption{Assembly costs in seconds}
\end{subtable}
\begin{subtable}{1\textwidth}
\centering
\begin{tabular}{c|rr|rr|rr|rr}
     & \multicolumn{2}{c|}{$p=2$}  & \multicolumn{2}{c|}{$p=3$}  & \multicolumn{2}{c|}{$p=4$} & \multicolumn{2}{c}{$p=5$}  \\ 
           & ILUT  & GS  & ILUT   &GS   &  ILUT & GS  & ILUT & GS     \\ \hline 
$h=2^{-6}$ & $0.10$   &$-$ & $0.27$	  &$-$ & $0.55$	&$-$ & $0.93$   &$-$ \\
$h=2^{-7}$ & $0.50$   &$-$ & $1.43$	  &$-$ & $2.85$	&$-$ & $4.76$   &$-$\\
$h=2^{-8}$ & $2.13$   &$-$ & $6.70$	  &$-$ & $13.98$&$-$ & $25.32$  &$-$  \\
$h=2^{-9}$ & $8.66$  &$-$ & $29.78$   &$-$ & $75.89$&$-$ & $134.57$ &$-$  \\
\end{tabular}
\subcaption{Factorization costs in seconds}
\end{subtable}
\begin{subtable}{1\textwidth}
\centering
\begin{tabular}{c|rr|rr|rr|rr}
     & \multicolumn{2}{c|}{$p=2$}  & \multicolumn{2}{c|}{$p=3$}  & \multicolumn{2}{c|}{$p=4$} & \multicolumn{2}{c}{$p=5$}  \\ 
           & ILUT  & GS  & ILUT   &GS   &  ILUT & GS  & ILUT & GS     \\ \hline 
$h=2^{-6}$ & $0.02$   &$0.11$& $0.02$ &$0.29$ & $0.03$	&$1.03$  &$0.03$&  $3.74$  \\
$h=2^{-7}$ & $0.07$   &$0.35$& $0.07$ &$1.00$ & $0.09$	&$3.54$  &$0.11$&  $13.24$  \\
$h=2^{-8}$ & $0.32$   &$1.46$& $0.26$ &$3.82$ & $0.36$	&$13.62$ &$0.49$&  $50.39$ \\
$h=2^{-9}$ & $1.30$   &$6.48$& $1.07$ &$15.80$& $1.46$	&$56.15$ &$1.88$&  $195.28$ \\
\end{tabular}
\subcaption{Solver costs in seconds}
\end{subtable}
\begin{subtable}{1\textwidth}
\centering
\begin{tabular}{c|rr|rr|rr|rr}
     & \multicolumn{2}{c|}{$p=2$}  & \multicolumn{2}{c|}{$p=3$}  & \multicolumn{2}{c|}{$p=4$} & \multicolumn{2}{c}{$p=5$}  \\ 
           & ILUT  & GS  & ILUT   &GS   &  ILUT & GS  & ILUT & GS     \\ \hline 
$h=2^{-6}$ & $0.77$   &$0.77$&  $1.49$  &$1.48$ & $2.79$	&$3.17$   &$4.59$&   $7.65$  \\
$h=2^{-7}$ & $3.15$   &$2.98$&  $6.14$  &$5.52$ & $11.52$	&$12.07$  &$19.93$&  $28.02$  \\
$h=2^{-8}$ & $12.67$  &$12.23$& $25.89$ &$23.23$& $49.11$	&$47.69$  &$75.97$&  $111.20$ \\
$h=2^{-9}$ & $51.48$  &$48.34$& $105.59$&$92.50$& $213.14$	&$187.81$ &$380.09$& $443.39$ \\
\end{tabular}
\subcaption{Total costs in seconds}
\end{subtable}
\caption{CPU timings for the first benchmark using $p$-multigrid.}
\label{tab:CPU_p}
\end{table}
\section{CPU times $h$-multigrid}
\label{cpu_hmg}
Table \ref{tab:CPU_h} shows the CPU timings with $h$-multigrid as a stand-alone solver for the first benchmark. For each configuration, the assembly, factorization and solver costs are shown separately. 

\begin{table}[h!]
\begin{subtable}{1\textwidth}
\centering
\begin{tabular}{c|rr|rr|rr|rr}
     & \multicolumn{2}{c|}{$p=2$}  & \multicolumn{2}{c|}{$p=3$}  & \multicolumn{2}{c|}{$p=4$} & \multicolumn{2}{c}{$p=5$}  \\ 
           & ILUT  & GS  & ILUT   &GS   &  ILUT & GS  & ILUT & GS     \\ \hline 
$h=2^{-6}$ & $0.47$   &$0.50$ &  $1.11$	  &$1.08$ &	$2.31$	&$2.16$  &	$4.32$   & $4.14$\\
$h=2^{-7}$ & $2.03$   &$1.80$ &  $4.21$	  &$4.00$ &	$9.35$	&$9.16$  &	$15.82$  & $17.28$\\
$h=2^{-8}$ & $7.68$   &$7.17$ &  $17.33$  &$17.64$& $36.17$	&$33.53$ &	$64.25$  & $68.57$    \\
$h=2^{-9}$ & $31.77$  &$30.14$&  $67.10$  &$64.69$& $143.86$&$145.38$&  $272.68$ & $265.49$    \\
\end{tabular}
\subcaption{Assembly costs in seconds}
\end{subtable}
\begin{subtable}{1\textwidth}
\centering
\begin{tabular}{c|rr|rr|rr|rr}
     & \multicolumn{2}{c|}{$p=2$}  & \multicolumn{2}{c|}{$p=3$}  & \multicolumn{2}{c|}{$p=4$} & \multicolumn{2}{c}{$p=5$}  \\ 
           & ILUT  & GS  & ILUT   &GS   &  ILUT & GS  & ILUT & GS     \\ \hline 
$h=2^{-6}$ & $0.13$   &$-$ & $0.35$	  &$-$ & $0.71$	&$-$ & $1.22$   &$-$ \\
$h=2^{-7}$ & $0.70$   &$-$ & $1.84$	  &$-$ & $3.85$	&$-$ & $6.12$   &$-$\\
$h=2^{-8}$ & $2.93$   &$-$ & $8.98$	  &$-$ & $18.74$&$-$ & $32.94$  &$-$  \\
$h=2^{-9}$ & $12.52$  &$-$ & $40.53$  &$-$ & $94.89$&$-$ & $178.58$ &$-$  \\
\end{tabular}
\subcaption{Factorization costs in seconds}
\end{subtable}
\begin{subtable}{1\textwidth}
\centering
\begin{tabular}{c|rr|rr|rr|rr}
     & \multicolumn{2}{c|}{$p=2$}  & \multicolumn{2}{c|}{$p=3$}  & \multicolumn{2}{c|}{$p=4$} & \multicolumn{2}{c}{$p=5$}  \\ 
           & ILUT  & GS  & ILUT   &GS   &  ILUT & GS  & ILUT & GS     \\ \hline 
$h=2^{-6}$ & $0.01$   &$0.07$& $0.02$ &$0.23$ & $0.02$	&$0.98$  &$0.04$&  $4.24$  \\
$h=2^{-7}$ & $0.05$   &$0.18$& $0.06$ &$0.66$ & $0.09$	&$3.12$  &$0.13$&  $13.26$  \\
$h=2^{-8}$ & $0.21$   &$0.70$& $0.21$ &$2.49$ & $0.34$	&$10.85$ &$0.49$&  $45.14$ \\
$h=2^{-9}$ & $0.86$   &$3.01$& $0.87$ &$9.93$ & $1.35$	&$43.12$ &$1.94$&  $168.74$ \\
\end{tabular}
\subcaption{Solver costs in seconds}
\end{subtable}
\begin{subtable}{1\textwidth}
\centering
\begin{tabular}{c|rr|rr|rr|rr}
     & \multicolumn{2}{c|}{$p=2$}  & \multicolumn{2}{c|}{$p=3$}  & \multicolumn{2}{c|}{$p=4$} & \multicolumn{2}{c}{$p=5$}  \\ 
           & ILUT  & GS  & ILUT   &GS   &  ILUT & GS  & ILUT & GS     \\ \hline 
$h=2^{-6}$ & $0.61$   &$0.57$&  $1.48$  &$1.31$ & $3.05$	&$3.14$   &$5.59$&   $8.39$  \\
$h=2^{-7}$ & $2.78$   &$1.98$&  $6.11$  &$4.66$ & $13.29$	&$12.28$  &$22.07$&  $30.54$  \\
$h=2^{-8}$ & $10.82$  &$7.87$&  $26.52$ &$20.13$& $55.25$	&$44.38$  &$97.68$&  $113.71$ \\
$h=2^{-9}$ & $45.15$  &$33.15$& $108.50$&$74.62$& $240.10$	&$188.50$ &$453.20$& $434.22$ \\
\end{tabular}
\subcaption{Total costs in seconds}
\end{subtable}
\caption{CPU timings for the first benchmark using $h$-multigrid with different smoothers.}
\label{tab:CPU_h}
\end{table}
\newpage 

\section{CPU times compared to an alternative smoother}
\label{cpu_scms}

Table \ref{tab:cpu_scms} shows the CPU timings with $h$-multigrid adopting the smoother from \cite{tak1} and a $p$-multigrid method using ILUT. For each configuration, the assembly, factorization and solver costs are shown separately. 

\begin{table}[h!]
\centering
\begin{subtable}{1\textwidth}
\centering
\begin{tabular}{c|rr|rr|rr|rr}
     & \multicolumn{2}{c|}{$p=2$}  & \multicolumn{2}{c|}{$p=3$}  & \multicolumn{2}{c|}{$p=4$} & \multicolumn{2}{c}{$p=5$}  \\ 
           & $p$-MG  & $h$-MG & $p$-MG   &$h$-MG  & $p$-MG  & $h$-MG  & $p$-MG & $h$-MG     \\ \hline 
$h=2^{-6}$ & $0.42$  &$0.30$  & $0.76$   &$0.71$  & $1.46$  &$1.62$   &$2.70$  & $3.19$  \\
$h=2^{-7}$ & $1.62$  &$1.19$  & $3.04$   &$2.86$  & $6.26$  &$6.47$   &$11.74$ & $11.73$ \\
$h=2^{-8}$ & $6.50$  &$4.88$  & $12.47$  &$11.59$ & $24.19$ &$25.86$  &$46.63$ & $49.67$ \\
$h=2^{-9}$ & $26.91$ &$19.44$ & $47.85$  &$46.19$ & $93.50$ &$104.10$ &$184.05$& $202.64$\\
\end{tabular}
\subcaption{Assembly costs}
\end{subtable}
\begin{subtable}{1\textwidth}
\centering
\begin{tabular}{c|rr|rr|rr|rr}
     & \multicolumn{2}{c|}{$p=2$}  & \multicolumn{2}{c|}{$p=3$}  & \multicolumn{2}{c|}{$p=4$} & \multicolumn{2}{c}{$p=5$}  \\ 
           & $p$-MG  & $h$-MG  & $p$-MG   &$h$-MG  & $p$-MG   & $h$-MG  & $p$-MG & $h$-MG     \\ \hline 
$h=2^{-6}$ & $0.06$  &$0.01$   & $0.19$   &$0.01$  & $0.41$  &$0.02$ &$0.78$  &	$0.02$\\
$h=2^{-7}$ & $0.30$  &$0.02$   & $0.96$   &$0.03$  & $2.20$  &$0.06$ &$4.64$  &	$0.08$ \\
$h=2^{-8}$ & $1.25$  &$0.09$   & $4.16$   &$0.12$  & $10.64$ &$0.22$ &$22.53$ & $0.27$\\
$h=2^{-9}$ & $5.26$  &$0.33$   & $17.44$  &$0.43$  & $44.30$ &$0.88$ &$120.64$& $1.04$\\
\end{tabular}
\subcaption{Set-up smoother costs}
\end{subtable}
\begin{subtable}{1\textwidth}
\centering
\begin{tabular}{c|rr|rr|rr|rr}
     & \multicolumn{2}{c|}{$p=2$}  & \multicolumn{2}{c|}{$p=3$}  & \multicolumn{2}{c|}{$p=4$} & \multicolumn{2}{c}{$p=5$}  \\ 
           & $p$-MG  & $h$-MG  & $p$-MG    &$h$-MG  & $p$-MG   & $h$-MG  & $p$-MG & $h$-MG     \\ \hline 
$h=2^{-6}$ & $0.02$  &$0.10$ & $0.02$  &$0.16$ & $0.03$	&$0.26$   &$0.05$  & $0.39$\\
$h=2^{-7}$ & $0.06$  &$0.34$ & $0.09$  &$0.52$ & $0.09$	&$0.86$   &$0.16$  & $1.22$\\
$h=2^{-8}$ & $0.23$  &$1.42$ & $0.26$  &$2.08$ & $0.48$	&$3.47$   &$0.54$  & $4.63$ \\
$h=2^{-9}$ & $0.93$  &$6.65$ & $1.09$  &$8.90$ & $1.99$	&$14.49$  &$2.12$  & $18.87$\\
\end{tabular}
\subcaption{Solver costs}
\end{subtable}
\begin{subtable}{1\textwidth}
\centering
\begin{tabular}{c|rr|rr|rr|rr}
     & \multicolumn{2}{c|}{$p=2$}  & \multicolumn{2}{c|}{$p=3$}  & \multicolumn{2}{c|}{$p=4$} & \multicolumn{2}{c}{$p=5$}  \\ 
           & $p$-MG  & $h$-MG  & $p$-MG    &$h$-MG  & $p$-MG   & $h$-MG  & $p$-MG & $h$-MG     \\ \hline 
$h=2^{-6}$ & $0.50$  &$0.41$ & $0.97$  &$0.88$ & $1.90$	 &$1.90$   &$3.53$   & $3.60$\\
$h=2^{-7}$ & $1.98$  &$1.55$ & $4.09$  &$3.41$ & $8.55$	 &$7.39$   &$16.54$  & $13.03$\\
$h=2^{-8}$ & $7.98$  &$6.39$ & $16.89$ &$13.79$& $35.31$ &$29.55$  &$69.70$  & $54.57$ \\
$h=2^{-9}$ & $33.10$ &$26.42$& $66.38$ &$55.52$& $139.79$&$119.47$ &$306.81$ & $222.55$\\
\end{tabular}
\subcaption{Total costs}
\end{subtable}
\caption{CPU times (in seconds) for convergence with $p$-multigrid and $h$-multigrid.}
\label{tab:cpu_scms}
\end{table}

\end{document}